\UseRawInputEncoding
\documentclass[10pt]{article}

\usepackage{amssymb,amsmath,amsthm,epsfig}
\numberwithin{equation}{section}
\usepackage{latexsym, enumerate}
\usepackage{eepic}
\usepackage{epic}
\usepackage{color}
\usepackage{ifpdf}
\usepackage{subfigure}

\usepackage{dsfont}
\usepackage{multirow}
\usepackage{makecell}
\usepackage{bm}
\usepackage{epstopdf}
\usepackage{graphicx}
\graphicspath{{tuii/CT/}{tuii/Denoise/}{tuii/Deblurring/}{tuii/Heat/}{tuii/Decon/}{tuii/Diagnosis_mh/}{tuii/}}

\usepackage{hyperref}
\hypersetup{hypertex=true,
    colorlinks=true,
    linkcolor=blue,
    filecolor=blue,
    urlcolor=blue,
    citecolor=blue,
    bookmarksopen=true,
}

\usepackage{booktabs}
\usepackage{caption}

\usepackage[linesnumbered,ruled]{algorithm2e}

\theoremstyle{plain}

\theoremstyle{definition}
\newtheorem{defn}{Definition}[section]

\theoremstyle{remark}

\begin{document}

\title{\large\bf Solving linear Bayesian inverse problems using a fractional total variation-Gaussian (FTG) prior and transport map}
\author{
Zejun Sun\thanks
{College of Mathematics, Hunan University, Changsha 410082, China.
Email: sunzejun@hnu.edu.cn}
\and
Guang-Hui Zheng\thanks
{College of Mathematics, Hunan Provincial Key Laboratory of Intelligent Information Processing and Applied Mathematics, Hunan University, Changsha 410082, China.
Email: zhenggh2012@hnu.edu.cn (Corresponding author)}
}
\date{}
\maketitle
\begin{center}{\bf ABSTRACT}
\end{center}\smallskip
The Bayesian inference is widely used in many scientific and engineering problems, especially in the linear inverse problems in infinite-dimensional setting where the unknowns are functions.
In such problems, choosing an appropriate prior distribution is an important task.
In particular, when the function to infer has much detail information, such as many sharp jumps, corners, and the discontinuous and nonsmooth oscillation, the so-called total variation-Gaussian (TG) prior is proposed in function space to address it.
However, the TG prior is easy to lead the blocky (staircase) effect in numerical results.
In this work, we present a fractional order-TG (FTG) hybrid prior to deal with such problems, where the fractional order total variation (FTV) term is used to capture the detail information of the unknowns and simultaneously uses the Gaussian measure to ensure that it results in a well-defined posterior measure.
For the numerical implementations of linear inverse problems in function spaces, we also propose an efficient independence sampler based on a transport map, which uses a proposal distribution derived from a diagonal map, and the acceptance probability associated to the proposal is independent of discretization dimensionality.
And in order to take full advantage of the transport map, the hierarchical Bayesian framework is applied to flexibly determine the regularization parameter.
Finally we provide some numerical examples to demonstrate the performance of the FTG prior and the efficiency and robustness of the proposed independence sampler method.

\smallskip
{\bf keywords}: Bayesian inference, fractional order total variation, transport map, independence sampler, hybrid prior

\section{Introduction}
The Bayesian inference methods \cite{Gelman_2013, 2005_Jari} have been popular among solving inverse problems, where the estimation results uncertainties can be quantified by learning the statistical information such as moments, confidence intervals, and the marginalizing.
A typical Bayesian inference is learning from the measurement data by incorporating some prior information to yield the posterior distribution and has  different features compared with classical deterministic regularization methods.
As we know, the most practical inverse problems are often highly ill-posed due to the limited and noisy measurement data and the prior distribution plays a significant role on the performance of the Bayesian inference.
Recently the Bayesian inverse problems in the infinite dimensional setting have been extensively studied and \cite{Dashti_2015, stuart_2010} privode a Bayesian framework for  the inverse problems in function spaces where the unknowns are of infinite dimension.
The Gaussian measures are widely used as the prior distributions for the infinite dimensional Bayesian inverse problems.
In fact such a choice has many advantages about theories and computation in the infinite dimensional Bayesian inference \cite{Dashti_2015}.

In our work, we focus on the linear inverse problem in infinite-dimensional setting.
Let $\mathbf{X}$ be separable Hilbert space, equipped with the Borel $\sigma$-algebra, and $\mathcal{A} : \mathbf{X} \to \mathbb{R}^n$ a linear forward model (unknows-to-observation).
We wish to solve the inverse problem of finding the unknown function $u\in\mathbf{X}$ (in this paper we shall restrict ourselves to the situation where $u$ is a real-valued function defined in $\mathbb{R}^m,\; m=1,2$) from measurement $\mathbf{y}\in \mathbb{R}^n$, which is usually generated by
\begin{align}
\label{Model}
\mathbf{y}=\mathcal{A}u+\boldsymbol{\eta},
\end{align}
where the noise $\boldsymbol{\eta}$ is assumed to be a $n$-dimensional zero-mean Gaussian random variable with covariance matrix $\Sigma_{\boldsymbol{\eta}}$.
In Bayesian inversion, the unknown model input and the measurement data are usually regarded as random variables.
The posterior distribution can be in a closed form only in some special cases.
For example, if the prior and noise are Gaussian, the posterior distribution is Gaussian for the linear forward model.

However the Gaussian prior distributions is usually applied to recover the unknowns with smooth property and in many practical problems, especially in the image reconstructions, are not suitable for modeling such inversion functions with smooth, sharp jumps or discontinuities \cite{Yao_2016}.
In order to model such functions, the total variation (TV) regularization \cite{RUDIN1992259} was proposed in the deterministic inverse problem context.
The TV prior for the Bayesian inference is applied to a range of inverse problems \cite{babacan2008variational, lassas2004can}.
However, with the discretization dimension increasing, the posterior distribution based on TV  prior may not converge to a well-defined infinite-dimensional measure \cite{lassas2004can}.
And the TV prior distributions differ significantly from the Gaussian measures, which many analysis techniques and numerical implementation methods based on Gaussian priors can not be simply and directly extended to non-Gaussian priors.
Recently, \cite{Yao_2016} has proposed a TV-Gaussian (TG) prior to deal with this issue and the corresponding numerical method is given.
The main idea of the TG prior is that the TV term is effective for identifying the sharp jumps and the Gaussian distribution is mainly used to ensure the posterior measure is well-defined in the function space.

In fact, although the total variational regularization can catch the sharp jumps of the unknowns,  the reconstruction results are easy to have the blocky effect \cite{2002_Vogel}.
And the Riemann-Liouville fractional order derivative based total variational regularization constructed in \cite{zhang2015variational} has a better result than that of the total variational regularization in addressing the image reconstruction problems.
Inspired by the success of this nonlocal regularization, in this paper, we propose a hybrid FTV-Gaussian (FTG) prior that we shall replace the total variational term in the hybrid TG prior with a fractional order total variational (FTV) term.
Namely, it uses a FTV term to capture the detail information in the unknowns and the Gaussian distribution makes sure that the posterior distribution is a well-defined probability measure in the function space.
The FTG prior can effectively reduces blocky effects and have better performance than the TG for the detail information in the unknowns, especially for recovering textures in images. More details are present in Section \ref{Numerical examples} about the numerical examples.

Practically we need to explore the posterior distribution in the Bayesian inference.
Markov chain Monte Carlo (MCMC) simulations are flexible  and widely used to draw samples from the complex posterior probability distribution in Bayesian inferences \cite{tierney1994markov, robert2004monte, Gelman_2013}.
As we all know, the standard MCMC algorithms, such as Metropolis-Hastings algorithm, can become arbitrarily slow as the discretization mesh of the unknown function is refined \cite{Cotter_2013, Roberts_2001}, and in this case the algorithm is said to be \textit{dimension-dependent}.
\cite{Cotter_2013} presented a dimension-independent MCMC algorithm that is obtained by constructing a preconditioned Crank-Nicolson (pCN) discretization of a stochastic partial differential equation that preserves the reference measure.
The sampling efficiency of MCMC algorithms for the infinite-dimensional Bayesian inverse problems can be further improved by incorporating the data information in the proposal design that can reflect the local or global geometry of the target distribution, like doing in finite dimensional counterparts, such as the stochastic Newton MCMC \cite{2012_Martin}, the dimension-independent likelihood-informed MCMC \cite{2016_Tiangang} and the adaptive independence sampler \cite{2018_Feng}.
The proposal distribution for the independence sampler is represented as a mixture of a finite number of specially parametrized Gaussian measures in \cite{2018_Feng}.
\cite{Parno_2018} constructed a more effective and adaptive proposal based transport maps to accelerate MCMC and in \cite{Peherstorfer_2019}, a multifidelity preconditioner is proposed to increase the efficiency of MCMC sampling.
This multifidelity approach exploits low-fidelity models to construct a proposal distribution that approximates the posterior distribution via a transport map, then uses this proposal distribution to perform MCMC sampling of the original (high-fidelity) posterior distribution.
However, in this work, we propose a proposal distribution based on a diagonal transport map to accelerate the independence sampler for the linear inverse problem in infinite-dimensional setting.
And the proposed independence sampler based on a diagonal map is a dimension-independent method, which the acceptance probability associated to the proposal derived from the transport map is independent of discretization dimensionality.
In order to make full use of the transport map, the hierarchical Bayesian framework \cite{JIN20107317} is also applied here, where the regularization parameter can be flexibly determined.
Loosely speaking, this approach has two stages: firstly we construct a diagonal transport map that can approximately push forward the reference measure to the posterior measure (not low-fidelity models in \cite{Peherstorfer_2019}) through solving a optimization problem and at the same time the regularization parameter can be determined, and secondly the posterior measure is explored by an independence sampler using a proposal distribution derived from the diagonal map.
Our numerical results in Section \ref{1d_ex2}, \ref{1d_ex1} illustrate the high efficiency of the  diagonal map-based independence sampler algorithm and its robustness with respect to various parameters, such as the noise level and the parameters in hyper-prior.

In summary, the main contributions of this paper are two-fold: (1) We propose a FTG hybrid  prior to handle unknowns with some complex oscillation features and textural information that cannot be well modeled by TV-Gaussian hybrid prior. In particular, it is able to effectively reduces blocky effects. (2) We also provide an efficient independence sampler with a proposal distribution based on a diagonal transport map, which is more efficient than standard pCN. (3) By using FTG hybrid prior, transport map and hierarchical Bayesian which can choose regularization parameters flexibly in combination, we significantly improve the accuracy and efficiency of Bayesian inference for several linear inverse problems of varying dimension, involving the deconvolution problem, inverse source problem, limited computed tomography and image denoising.
The rest of the paper is organized as follows.
In Section \ref{section2}, we introduce the fundamentals of Bayesian inference and the fractional order total variation Gaussian (FTG) prior.
Section \ref{section3} first describes the hierarchical modeling, and then introduces the transport maps in the context of Bayesian inverse problems, explains how the diagonal transport maps can be constructed from reference measure and discusses a numerical solution for the optimization problem.
And the independence sampler using a proposal distribution derived from a diagonal map is described in this section.
Section \ref{section4} presents a range of numerical examples for the linear inverse problems from the one-dimensional deconvolution problem and inverse source identification problem to the two-dimensional limited computed tomography (CT) reconstruction in medical imaging and the image denoising.
The paper ends with some conclusions in Section \ref{section5}.


\section{The FTG priors}\label{section2}
\subsection{The Bayesian framework and hybrid priors}

We first give a brief introduction to the Bayesian inference framework for the infinite-dimensional linear inverse problems.
According to \eqref{Model}, the likelihood function which is law of $\mathbf{y}$ conditional on $u$ is
\begin{align}
\label{likelihood}
L_{\mathbf{y}}(u)\varpropto \exp(-\Phi_{\mathbf{y}}(u)),
\end{align}
where the data-misfit function
\begin{align}
\label{data-misfit}
\Phi_{\mathbf{y}}(u):=\frac{1}{2}\|\mathcal{A}u-\mathbf{y}\|^2_{\Sigma_{\boldsymbol{\eta}}},
\end{align}
where $\|f\|_{\Sigma}:=\|\Sigma^{-\tfrac{1}{2}}f\|^2_2$, for any positive symmetric matrix $\Sigma$, and the $\|\cdot\|_2$ is the Euclidean norm.
We then choose a prior probabilistic measure of $u$ denoted by $\mu_{pr}$ which is related to the posterior measure of $u$,denoted by $\mu^{\mathbf{y}}$, through the Radon-Nikodym (R-N) derivative \cite{stuart_2010}
\begin{align}
\label{RN}
\frac{\mathrm{d}\mu^{\mathbf{y}}}{\mathrm{d}\mu_{pr}}(u)=\frac{1}{Z}\exp(-\Phi_{\mathbf{y}}(u)),
\end{align}
where $Z=\int_{\mathbf{X}}\exp(-\Phi_{\mathbf{y}}(u))\mathrm{d}\mu_{pr}(u)$ is a normalization constant, which can be interpreted as the Bayes' rule in the infinite-dimensional setting. The most popular prior $\mu_{pr}$ is chosen as the Gaussian measure $\mu_0$, i.e. $\mu_{pr}=\mu_0 =\mathcal{N}(0,C_0)$, a zero mean and covariance operator $C_0$ Gaussian measure on $\mathbf{X}$.
Note that the $C_0$ is symmetric positive and of trace class.

Howerer, to better reconstruct the detail informations of unknows, the hybrid prior \cite{Yao_2016} is necessary which takes the form
\begin{align}
\label{add_prior}
\frac{\mathrm{d}\mu_{pr}}{\mathrm{d}\mu_0}(u)\varpropto \exp(-J(u)),
\end{align}
where $J(u)$ is the additional prior (or the regularization) information on $u$.
As a result,  the R-N derivative of $\mu^{\mathbf{y}}$ with respect to $\mu_0$ is
\begin{align}
\label{RN-addprior}
\frac{\mathrm{d}\mu^{\mathbf{y}}}{\mathrm{d}\mu_0}(u)\varpropto \exp(-\Phi_{\mathbf{y}}(u)-J(u)).
\end{align}
In this case, it is well-defined on $\mathbf{X}$ for the posterior measure $\mu^{\mathbf{y}}$ associated with the linear forward model satisfied with certain conditions \cite{Dashti_2015, stuart_2010}.
And in \cite{Yao_2016}, a TV-Gaussian (TG) prior is presented that  total variation (TV) term is used to deal with sharp jumps of the function, and the Gaussian measure is used as a reference measure so that it results in a well-defined posterior measure in the function space.
However, in this paper we introduce the fractional order TG (FTG) prior to better address unknown function in the linear forward, especially in the image reconstruction.

\subsection{The FTG prior}
Here we provide the formulation of the FTG prior and we start with briefly reviewing some basic knowledge on fractional derivatives \cite{Samko1993FractionalIA, Kilbas2006TheoryAA}.
\begin{defn}[Riemann-Liouville (RL) fractional derivative]
Let $n-1<\alpha \le n,\, n\in \mathbb{R}^+$, and $[a,b]\subset \mathbb{R}$, then Riemann-Liouville fractional derivative of a function $f$ are defined as follows:
the left derivative
\[
D^{\alpha}_{[a,x]}f(x):=\frac{1}{\Gamma(n-\alpha)}\left(\frac{\mathrm{d}}{\mathrm{d}x}\right)^n\int^x_{\alpha}\frac{f(\tau)}{(x-\tau)^{\alpha-n+1}}\mathrm{d}\tau,
\]
the right derivative
\[
D^{\alpha}_{[x,b]}f(x):=\frac{1}{\Gamma(n-\alpha)}\left(-\frac{\mathrm{d}}{\mathrm{d}x}\right)^n\int^b_{x}\frac{f(\tau)}{(\tau-x)^{\alpha-n+1}}\mathrm{d}\tau,
\]
where the $\varGamma(\cdot)$ is the standard Gamma function and the corresponding Riesz-RL (center) derivative is defined as
\[
D^{\alpha}_{[a,b]}f(x):=\frac{1}{2}\Bigl(D^{\alpha}_{[a,x]}f(x)+(-1)^nD^{\alpha}_{[x,b]}f(x)\Bigr).
\]
\end{defn}
Then the fractional order $\alpha$ total variation of a function $\varphi$ is given by
\[
TV^{\alpha}(\varphi):=\sup_{\phi\in K}\int_{\Omega}-\varphi \mathrm{div}^{\alpha}\phi\mathrm{d}x,
\]
where $\mathrm{div}^{\alpha}\phi=\sum^d_{i=1}\frac{\partial^{\alpha}\phi_i}{\partial x^{\alpha}_i}$ for $\phi=(\phi_1,\cdots,\phi_d)$, and the $\frac{\partial^{\alpha}\phi_i}{\partial x^{\alpha}_i}$ represents  fractional $\alpha$ order derivative $D^{\alpha}_{[a,b]}\phi_i$ of $\phi_i$ along the $x_i$ direction.
We define the fractional Sobolev space as
\[
W_p^{\alpha}(\Omega)=\Bigl\{u\in L^p(\Omega)\bigl| \; \|u\|_{W_p^{\alpha}(\Omega)}<+\infty\Bigr\},
\]
and the corresponding norm:
\[
\|u\|_{W_p^{\alpha}(\Omega)}=\Bigl(\int_{\Omega}|u|^p\mathrm{d}x+\int_{\Omega}|\nabla^{\alpha}u|^p\mathrm{d}x\Bigr)^{\frac{1}{p}},
\]
where $p$ is any positive integer, $\nabla^{\alpha}=\Bigl(\frac{\partial^{\alpha}}{\partial x^{\alpha}_1},\dots,\frac{\partial^{\alpha}}{\partial x^{\alpha}_m}\Bigr)$, and $\Omega\subseteq\mathbb{R}^m$.
In \cite{2015_Zhang}, it has been proved that the fractional sobolev space $W_p^{\alpha}(\Omega)$ is the Banach space, and has the following embedding relation
\begin{align}
\label{embed}
W_2^{\alpha}(\Omega)\subseteq W_1^{\alpha}(\Omega),
\end{align}
that is, for any $u\in W_1^{\alpha}(\Omega)$, there is a constant $C>0$ such that
\begin{align}
\|u\|_{W_1^{\alpha}(\Omega)}\le C\|u\|_{W_2^{\alpha}(\Omega)}.
\end{align}

In addition, \cite{2015_Zhang} has showed that for any $u\in W_1^{\alpha}(\Omega)$, the fractional order total variation $TV_{\alpha}(u)=\int_{\Omega} |\nabla^{\alpha}u|\mathrm{d}x$.
And for any $1\le p<\infty$ fractional order Sobolev space $W_{p}^{\alpha}(\Omega)$ is separable \cite{Wang_2021_1,Wang_2021_2}.
According to the definition of fractional Sobolev space, the norm of $W_2^{\alpha}(\Omega)$ space is induced by the inner product below
\[
\langle u,v\rangle_{W_2^{\alpha}(\Omega)}=\int_{\Omega}uv\mathrm{d}x+\int_{\Omega}\nabla^{\alpha}u\cdot\nabla^{\alpha}v\mathrm{d}x,\quad u,v\in W_2^{\alpha}(\Omega).
\]
We can obtain that $W_2^{\alpha}(\Omega)$ is the separable Hilbert space.
Thus we choose the $\mathbf{X}=W_2^{\alpha}(\Omega)$ and consider the embedding relation \eqref{embed}.
Then we can get the following $\alpha-$order FTV regularization term or FTG prior
\begin{align}
J(u;\lambda)=\frac{\lambda}{2}\| u\|_{TV^{\alpha}},
\end{align}
where $\lambda$ is a regularization parameter, and
$
\| u\|_{TV^{\alpha}}=\int_{\Omega} |\nabla^{\alpha}u|\mathrm{d}x.
$

\cite{Wang_2021_1,Wang_2021_2} have show that the equation \eqref{add_prior} is a well-behaved prior under certain assumptions on the FTG prior $J$, and the posterior $\mu^{\mathbf{y}}$ is a well-defined probability measure on the $\mathbf{X}$ and it is also Lipschitz in the data $\mathbf{y}$.
Moreover, under some additional assumptions, the posterior measure can be well approximated by a measure defined in a finite-dimensional space.

\section{Diagonal map-based independence sampler}\label{section3}
As we know, the commonly used approach to exploring the posterior is the Markov chain Monte Carlo \cite{gamerman2006markov}, which can generate a stream of samples to estimate some statistical information of complex target probability distributions.
However, samples generated by MCMC are necessarily correlated among successive samples, which means smaller effective sample sizes and larger errors in posterior estimates.
Thus the algorithm usually needs large number of samples to make accurate estimates, in particular, in the infinite-dimensional setting, which leads to bring a heavy computation burden.
In this section, in order to increase the efficiency of independence sampler, we construct a fixed proposal distribution that approximates the posterior distribution at hand through a diagonal transport map that pushes forward a simple reference distribution to the target probability distribution.
And, simultaneously the regularization parameters for the independence sampler can be determined automatically via \eqref{reg_para_MAP} after constructing a transport map at hand.
Moreover, we describe that diagonal map based independence sampler is well-defined in the function space.

We start by briefly reviewing the independence sampler MCMC.
To this end, we define the measures
\begin{align*}
\nu(du,dv)=\mu^{\mathbf{y}}(du)\mu(dv),\\
\nu^{\mathrm{T}}(du,dv)=\mu^{\mathbf{y}}(dv)\mu(du)
\end{align*}
on the product space $\mathbf{X}\times \mathbf{X}$, where the $\mu$ is given a proposal distribution.
When the $\nu$ is absolute continuous with respect to $\nu^{\mathrm{T}}$ \cite{tierney1998, Dashti_2015}, we then can define the acceptance probability
\begin{align}
\label{acc}
a(u, v)=\min\Bigl\{1,\frac{\mathrm{d}\nu}{\mathrm{d}\nu^{\mathrm{T}}}(u,v)\Bigr\} ,
\end{align}
where
\[
\frac{\mathrm{d}\nu}{\mathrm{d}\nu^{\mathrm{T}}}(u,v)=\frac{\mathrm{d}\mu^{\mathbf{y}}}{\mathrm{d}\mu}(v)\frac{\mathrm{d}\mu}{\mathrm{d}\mu^{\mathbf{y}}}(u).
\]
It can be seen that the acceptance probability of the independence sampler is well-defined if and only if the $\nu$ is absolute continuous with respect to $\nu^T$, which obviously need that $\mu^{\mathbf{y}}$ and $\mu$ are equivalent each other.
And it suffices to require that $\mu$ and $\mu_0$ are equivalent since $\mu^{\mathbf{y}}$ and $\mu_0$ are equivalent.
The dimension-indepandent pCN algorithm and independence sampler can be obtained by simply choosing the Gaussian prior as proposal distribution which is equivalent to the $\mu_0$ and invariant \cite{Dashti_2015}.

The standard pCN algorithm used the proposal as follow
\[
v=\sqrt{1-\beta^2}u+\beta \omega,
\]
where $\omega \sim \mathcal{N}(0,C_0)$ and the $\beta$ plays the role of the locality parameter, and then $v \sim \mathcal{N}(\sqrt{1-\beta^2}u, \beta^2C_0)$.
And the associated acceptance probability is
\begin{align}
a(u,v) = \min\Bigl\{1, \exp\Bigl(\Phi_{\mathbf{y}}(u)+J(u;\lambda)-\Phi_{\mathbf{y}}(v)-J(v;\lambda)\Bigr)\Bigr\},
\label{pCN_acc}
\end{align}
where the $\lambda$ is regularization parameter.
And the Algorithm \ref{alg:pCN} describes the standard pCN approach, which is used to compare the sampling with our independence sampler.

However, this simply proposed distribution according to the prior works well only when the data and likelihood is not too informative.
When the posterior is far away from the prior, and the $\Phi_{\mathbf{y}}$ varies extremely depending on where it is evaluated, the independence sampler will not work well.
Next we introduce an efficient proposal measure derived from a diagonal map for the independence sampler.

\begin{algorithm}[htbp]
\KwIn{Date-misfit function $\Phi_{\mathbf{y}}(u)$, prior $J(u;\lambda)$ and $\mu_0$, regularization parameter $\lambda$, total number of steps $K$.}
\KwOut{pCN samples of the target distribution, $\{u^{(i)}\}_1^K$.}
Choose the initial state $u^{(1)}\sim \mu_0$\;
\For{$i \gets 1 \ldots K-1$}{
Computer the propose $v=\sqrt{1-\beta^2}u^{(i)}+\beta \omega,\quad \omega \sim \mu_0$\;
Calculate the acceptance probability $a(u^{(i)},v)$ via the \eqref{pCN_acc} \;
Set $u^{(i+1)}=v$ with probability $a(u^{(i)},v)$; else set $u^{(i+1)}=u^{(i)}$\;
}
\Return{\textrm{Target samples} $\{u^{(i)}\}_1^K$.}
\caption{The standard pCN}
\label{alg:pCN}
\end{algorithm}

\subsection{Construction of transport maps}
Before describing the transport maps, we first introduce the hierarchical Bayesian formulation in order to make the most of it, which we can flexibly determine the regularization parameter.
\subsubsection{Hierarchical Bayesian formulation}
The probability density function corresponding to the target measure $\mu^{\mathbf{y}}$ is the posterior $\pi^{\mathbf{y}}$, and the probability density function corresponding to Gaussian measure $\mu_0$ is denoted by $\pi_0$. Then according to \eqref{RN-addprior}, we can obtain that
\begin{align}
\pi^{\mathbf{y}}(u)\varpropto \exp(-\Phi_{\mathbf{y}}(u)-J(u;\lambda))\pi_0(u),
\end{align}
We rewrite $\mathbf{u}\in\mathbb{R}^d$ to represent the unknowns and let $\Sigma_{\boldsymbol{\eta}}=\sigma^2I$ and the covariance operator  $C_0\in\mathbb{R}^{d\times d}$, where $\sigma$ is the standard deviation and $I$ is a $n$-dimensional identity matrix. Then the posterior density $\pi^{\mathbf{y}}$ can be obtained as
\begin{align}
\pi^{\mathbf{y}}(\mathbf{u})\varpropto \exp\Bigl(-\frac{1}{2\sigma^2}\|\mathcal{A}\mathbf{u}-\mathbf{y}\|^2_2-\frac{\lambda}{2}\|\mathbf{u}\|_{TV^{\alpha}}-\frac{1}{2}\|\mathbf{u}\|^2_{C_0}\Bigr).
\end{align}
The posterior density $\pi^{\mathbf{y}}$ provides the complete distribution of $u$ relying on the observations $\mathbf{y}$. Following the same steps as in\cite{Dashti_2013}, we can compute the MAP point with the FTG prior, $\mathbf{u}_{MAP}:=arg\;\max_{\mathbf{u}}\pi^{\mathbf{y}}(\mathbf{u})$, which is equivalent to the following minimization problem
\begin{align}
\min_{\mathbf{u}}\Bigl\{\frac{1}{2\sigma^2}\|\mathcal{A}\mathbf{u}-\mathbf{y}\|^2_2+\frac{\lambda}{2}\|\mathbf{u}\|_{TV^{\alpha}}+\frac{1}{2}\|\mathbf{u}\|^2_{C_0}\Bigr\},
\end{align}
where the regularization parameter $\lambda$ plays a critical role in classical inverse problems.
It is essential to choose a suitable $\lambda$ for Bayesian inverse problems.
However, as we know, it is nontrivial task to ascertain the regularization parameter in almost inverse problems.
Thanks to hierarchical Bayesian modeling \cite{JIN20107317, gamerman2006markov}, we can overcome the difficulty flexibly.
The idea is to let the data $\mathbf{y}$ determine the parameters in the hope of effectively diminishing the effect of the initial (prior) assumptions of their values on the inverse solution.
Then, the unknown function $u$ and the regularization parameter $\lambda$ can be identified at the same time.
In hierarchical Bayesian framework, $\lambda$ (a hyper-parameter) can be regarded as a random variable.
If we choose Gamma distribution $G(\lambda ; k, \vartheta)$ as the hyper-prior \cite{Gelman_2013} for $\lambda$, i.e.,
 \begin{align}
 G(\lambda ; k, \vartheta)=\frac{\vartheta^{k}}{\varGamma(k)}\lambda^{k-1}\exp(-\vartheta\lambda),
 \end{align}
 where the positive constants $k$ and $\vartheta$ is called shape parameter and rate parameter respectively.
Then the posterior density can be written as
\begin{align}
\label{H_postprior}
\pi^{\mathbf{y}}(\mathbf{u},\lambda)\propto \lambda^{k-1}\exp\Bigl(-\frac{1}{2\sigma^2}\|\mathcal{A}\mathbf{u}-\mathbf{y}\|^2_2-\frac{1}{2}\|\mathbf{u}\|^2_{C_0}-\frac{\lambda}{2}\|\mathbf{u}\|_{TV^{\alpha}}-\vartheta\lambda\Bigr),
\end{align}
and the MAP estimate for the posterior density \eqref{H_postprior} can be easily obtained by minimizing the following functional
\begin{align}
\label{H_MAP}
\mathcal{R}(\mathbf{u},\lambda)=\frac{1}{2\sigma^2}\|\mathcal{A}\mathbf{u}-\mathbf{y}\|^2_2+\frac{1}{2}\|\mathbf{u}\|^2_{C_0}+\vartheta\lambda+\frac{\lambda}{2}\|\mathbf{u}\|_{TV^{\alpha}}-(k-1)\mathrm{ln}(\lambda).
\end{align}
One salient feature of the functional is that the regularization parameter $\lambda$ can be computed by data-driven procedure.
And the minimization problem of $\eqref{H_MAP}$ can be solved by the iteratively reweighted approach \cite{Li_2010}, which is usually applied in compressed sensing.
However, we focus on the relationship between $\lambda$ and $\mathbf{u}$ from the \eqref{H_MAP} rather than its minimizer.
Now taking the partial derivative of $\mathcal{R}(\mathbf{u},\lambda)$ with respect to $\lambda$, we can obtain
\begin{align}
\label{partial_lambda}
\frac{\partial{\mathcal{R}(\mathbf{u},\lambda)}}{\partial{\lambda}}=\vartheta+\frac{1}{2}\|\mathbf{u}\|_{TV^{\alpha}}-\frac{k-1}{\lambda},
\end{align}
its second order partial derivative is $\tfrac{\partial^2{\mathcal{R}(\mathbf{u},\lambda)}}{\partial{\lambda^2}} =\tfrac{k-1}{\lambda^2}>0$ for $k>1$, and let $\frac{\partial{\mathcal{R}(\mathbf{u},\lambda)}}{\partial{\lambda}}=0$, we can get
\begin{align}
\label{reg_para_MAP}
\lambda=\frac{2(k-1)}{\|\mathbf{u}\|_{TV^{\alpha}}+2\vartheta}.
\end{align}
After construction of a transport map, the above equation can be used to calculate the regularization parameter for our independence sampler.

\subsubsection{Optimal transport}
In \cite{ELMOSELHY20127815}, the measure-preserving transport maps constructed via the solution of an optimization problem between continuous probability measures is first applied to Bayesian inference.
\cite{ELMOSELHY20127815} presented a variational approach to the construction of transport maps explicitly that pushes forward the prior measure to the posterior measure, which entirely avoid Markov chain simulation.
The transport approach of Benjamin et al. \cite{Peherstorfer_2019} and Parno et al. \cite{Parno_2018} instead follows a precondition MCMC sampling.
\cite{Parno_2018} used transport map to obtain a proposal distributions that can more effectively explore the target density and is adapted as the MCMC sampling proceeds.
And \cite{Peherstorfer_2019} proposed a multifidelity approach that the low-fidelity model is used to construct a transport map and the high-fidelity posterior distribution is explored using a non-Gaussian proposal distribution derived from the transport map.
Below, we follow \cite{Marzouk_2016, Peherstorfer_2019} to denote some notions of transport maps.

We will refer to the posterior measure $\mu^{\mathbf{y}}$ and Gaussian measure $\mu_{ref}$ as the target and reference measures on $\mathbb{R}^d$, respectively.
A transport map $\mathnormal{T}: \mathbb{R}^d \to \mathbb{R}^d$ is a deterministic coupling that pushes forward $\mu_{ref}$ to $\mu^{\mathbf{y}}$, satisfying
\begin{align}
\label{Td}
\mathnormal{T}_{\sharp}\mu_{ref}=\mu^{\mathbf{y}}.
\end{align}
In other words, $\mu^{\mathbf{y}}(B)=\mu_{ref}(T^{-1}(B))$ for any Borel set $B\in\mathbb{R}^d$.
Both the reference measure $\mu_{ref}$ and the target measure $\mu^{\mathbf{y}}$ are absolutely continuous with respect to the Lebesgue measure on $\mathbb{R}^d$ that assure the existence of transport maps satisfying \eqref{Td}.
Of course, there may be infinitely many such transport maps between the reference measure and the target measure.
One way of guaranteeing the uniqueness of map is to introduce a transport cost function and minimizes it simultaneously satisfying the constraint \eqref{Td}.
This minimization problem is called the Monge problem \cite{Vershik_2013, villani2003topics, villani2009optimal}, and its solution is the optimal transport map.
If the cost is taken to be a quadratic form in \cite{Bonnotte2013FromKR, Carlier_2010}, the optimal transport map is exactly the Knothe-Rosenblatt rearrangement \cite{Bonnotte2013FromKR, Carlier_2010, Rosenblatt1952RemarksOA}.
However, in this paper we directly assume that the transport map is a triangular diffeomorphism $T$ such that $\nabla T\succ0$ (i.e. monotone increasing) as in \cite{Parno_2018, ELMOSELHY20127815, Marzouk_2016, Peherstorfer_2019}, instead of being particularly concerned with the optimality aspect of the transport.

As noted above, the lower triangular maps take the form
\begin{align}
\label{T}
T(\mathbf{u})=
\begin{bmatrix}
\!\!\!T_1(u_1) \qquad\qquad\quad\\
\!\!\!\!\!\!T_2(u_1,u_2) \qquad\quad\\
\vdots \qquad\qquad\qquad\\
T_d(u_1, u_2, \ldots ,u_d) \\
\end{bmatrix}
\end{align}
where $\mathbf{u} \in \mathbb{R}^d$ and $T_k : \mathbb{R}^k \to \mathbb{R} \; (k=1,2,\dots,d)$ is $k$th component function of the transport map $T$.
In this setting, the $T_k$ of the map depends only on the first $i$ input variables, and it holds that $\det \nabla T>0$ (see \cite{Marzouk_2016} for more details).
Since the reference and the target measures are absolutely continuous, existence and uniqueness of such a lower triangular transport map (i.e. Knothe-Rosenblatt rearrangement) are guaranteed \cite{Bonnotte2013FromKR, Carlier_2010, Rosenblatt1952RemarksOA}.
To obtain a useful approximation of the transport map, we will define a map-induced density $\widetilde{\pi}_{ref}(u)$ and minimize the distance between the density $\pi_{ref}(u)$ of Gaussian measure $\mu_{ref}$ and  this map-induced density.
The next subsections describes the setup of this optimization problem.

\subsubsection{Optimization problems}
To set up the optimization problem, following \cite{Marzouk_2016}, the density form of \eqref{Td} can be written as $\mathnormal{T}_{\sharp}\pi_{ref}=\tfrac{1}{\tau}\pi^{\mathbf{y}}$, where $\tau$ is the normalization constant of posterior density and the transport map $T$ only move mass of $\mathbf{u}$ but not including the $\lambda$.
Now consider the pushforward of the reference density under the map $T$, and it is defined as
\begin{align}
\mathnormal{T}_{\sharp}\pi_{ref}(\mathbf{u}) := \pi_{ref}(T^{-1}(\mathbf{u}))|\det \nabla T^{-1}(\mathbf{u})|,
\end{align}
where $\det \nabla T^{-1}(\mathbf{u})$ denotes  the determinant of the Jacobian of the inverse of the map at $u$.
Then we can obtain the map-induced density
\begin{align}
\label{map_T}
\widetilde{\pi}_{ref}(\mathbf{u})=\mathnormal{\widetilde{T}}_{\sharp}^{-1}\Bigr[\frac{1}{\tau}\pi^{\mathbf{y}}(\mathbf{u},\lambda)\Bigl]=\frac{1}{\tau}\pi^{\mathbf{y}}(\widetilde{T}(\mathbf{u}),\lambda)|\det \nabla \widetilde{T}(\mathbf{u})|,
\end{align}
where $\widetilde{T}$ is an approximation of a transport map $T$, which will be obtained via numerical optimization.
If the reference density and the map-induced density are equal, i.e. $\pi_{ref}=\widetilde{\pi}_{ref}$, the $\widetilde{T}$ can exactly satisfy $\mathnormal{\widetilde{T}}_{\sharp}\mu_{ref}=\mu^{\mathbf{y}}$.
Thus we can minimize a distance between $\pi_{ref}$ and $\widetilde{\pi}_{ref}$ to obtain the $\widetilde{T}$.
In this paper, we use the Kullback-Leibler (KL) divergence to measure the distance between distributions as in \cite{Parno_2018, Peherstorfer_2019}.
Let the $\mathcal{D}_{KL}(\cdot||\cdot)$ and $\mathbb{E}_{\pi_{ref}}[\cdot]$ denote the KL divergence and integration with respect to the reference measure, respectively.
Then a minimizer of the optimization problem:
\begin{align}
\label{opt_problem}
\min \mathcal{D}_{KL}(\pi_{ref}||\widetilde{\pi}_{ref}) & =\mathbb{E}_{\pi_{ref}}\Bigl[\log\bigl(\frac{\pi_{ref}}{\widetilde{\pi}_{ref}}\bigr)\Bigr], \\
s.t. \nabla \widetilde{T} & \succ 0, \notag \\
\widetilde{T} & \in \mathcal{T}, \notag
\end{align}
is  a valid approximation of a transport map \cite{ELMOSELHY20127815}, where the constraint $\nabla \widetilde{T} \succ 0$ suffices to enforce mmonotonicity of a triangular map
and $\mathcal{T}$ is some space of smooth lower triangular functions from $\mathbb{R}^d$ to $\mathbb{R}^d$.
If $\mathcal{T}$ is rich enough, we will obtain $\mathcal{D}_{KL}(\pi_{ref}||\widetilde{\pi}_{ref}) =0$.
Then the solution $\widetilde{T}$ of this optimization problem will satisfy \eqref{Td} \cite{Marzouk_2016, ELMOSELHY20127815}.

Furthermore, form the \eqref{map_T}, the objective function of this optimization problem can be written as
\begin{align}
\mathcal{D}_{KL}(\pi_{ref}||\widetilde{\pi}_{ref})=\mathbb{E}_{\pi_{ref}}\Bigl[\log\pi_{ref}+\log\tau-\log\pi^{\mathbf{y}}(\widetilde{T},\lambda)-\log|\det \nabla \widetilde{T}|\Bigr],
\end{align}
where the $\mathbb{E}_{\pi_{ref}}[\log\pi_{ref}+\log\tau]$ is independent of the map $\widetilde{T}$ and thus a constant that can be ignored for the purposes of optimization.

\subsubsection{Diagonal approximation of triangular map}\label{approx_map}
To obtain the numerical solution of this optimization problem \eqref{opt_problem}, the infinite-dimensional function space $\mathcal{T}$ must be replaced with a finite-dimensional subspace $\mathcal{\widetilde{T}} \subset \mathcal{T}$.
In \cite{Peherstorfer_2019}, each component function of the approximation map $\widetilde{T}_k$ is parameterized with the integrated-squared ansatz, which can enforce the monotonicity constraints explicitly and capture nonlinear dependencies in the target measure.
However, this parameterization of the map is computationally expensive in the infinite dimensional setting.
In this paper, the $\widetilde{T}_k$ is parameterized by expanding it in a diagonal basis of univariate polynomials for the linear inverse problems.

Let each component of the map be written as $\widetilde{T}_k(\mathbf{a}_k; \mathbf{x}),\; k=1,2,  \dots, d$, where $\mathbf{a}_k \in \mathbb{R}^{q+1}$ is a column vector coefficients.
Then, we can express each component of the transport map $\widetilde{T} \in \mathcal{\widetilde{T}}$ as
\begin{align}
\label{map_pare}
\widetilde{T}_k(\mathbf{a}_k; \mathbf{x})=P_q(\mathbf{a}_k; x_k)=a_{k,0}+a_{k,1}x_k+a_{k,2}x_k^2+\dots+a_{k, q}x_k^q,
\end{align}
where $x_k$ is the $k$th component of $\mathbf{x}$, and $P_q(\mathbf{a}_k; x_k)$ is a $q$ degree univariate polynomial with respect to $x_k$, which forces $\widetilde{T}$ to be lower triangular.
%
It is easy to see that the number of coefficients is proportional to that of the mode parameters, and as a result the diagonal approximation of transport map $\widetilde{T}_k$ can be rapidly calculated by numerical optimization of \eqref{opt_problem}.
In fact, the diagonal parameterization \eqref{map_pare} is the same to that using the multi-index sets $\mathcal{J}^D_k$ with $p=q$ in \cite{Marzouk_2016}.

\subsubsection{Numerical optimization}\label{numercal opt}
There is need to approximate the expectation with respect to reference measure in the objective of \eqref{opt_problem} in the process of numerical optimization.
We approximate the expectation $\mathbb{E}_{\pi_{ref}}[\cdot]$ by its sample-average approximation (SAA) \cite{Kleywegt_2002}, i.e., a Monte Carlo estimator with $M$ independent samples, denoted by $\{\mathbf{x}^{(1)},\mathbf{x}^{(2)},\dots,\mathbf{x}^{(M)}\}$, from the reference measure $\mu_{ref}$.
Obviously, the $\mathbb{E}_{\pi_{ref}}[\cdot]$ can be calculated more accurately, as the cardinality of the sample set grows.
Then, the coefficients $(\mathbf{a}_1,\mathbf{a}_2,\dots,\mathbf{a}_d)$ is denoted by $F$ and we can obtain the optimization problem
\begin{align}
\label{d_opt_problem}
\min_{F} \; & \frac{1}{M}\sum^M_{i=1}\Biggl[-\log\pi^{\mathbf{y}}\biggl(\widetilde{T}(F; \mathbf{x}^{(i)}),\lambda\biggr)-\sum^d_{k=1}\log\frac{\partial\widetilde{T}_k}{\partial x_k}\biggl|_{\mathbf{x}^{(i)}}\Biggr], \\
s.t.& \; \frac{\partial\widetilde{T}_k}{\partial x_k}\biggl|_{\mathbf{x}^{(i)}} > 0, \quad k=1,2,\dots, d, \text{and}\; i=1,2,\dots, M, \notag
\end{align}
where we have simplified the monotonicity constraint $\nabla \widetilde{T} \succ 0$ by using the fact that  $\nabla \widetilde{T}$ is lower triangular and according to the \eqref{reg_para_MAP} we can explicitly obtain the $\lambda$ in $\pi^{\mathbf{y}}$,
\begin{align}
\label{de_reg}
\lambda=\frac{2(k-1)}{\bigl\|\mathbb{E}_{\pi_{ref}}[\widetilde{T}(\mathbf{x})]\bigr\|_{TV^{\alpha}}+2\vartheta},
\end{align}
where $\mathbb{E}_{\pi_{ref}}[\widetilde{T}(\mathbf{x})]=\tfrac{1}{M}\sum^M_{i=1}\widetilde{T}(F; \mathbf{x}^{(i)})$ via the SAA.
In our work, we use the alternating direction algorithm to solve this optimization problem, which is summarized in Algorithm \ref{alg:ad}.
Then, the above numerical optimization problem \eqref{d_opt_problem}, in fact, just require to optimize the parameter $F$ in $\widetilde{T}$.

If the unknowns are only endowed with a Gaussian prior $\mu_0$ for the linear forward model inverse problems,
then the posterior is the Gaussian, i.e. $\mathbf{u} \sim \mathcal{N}(\mu_p,\Sigma_p)$, where
\[
\mu_p =\Sigma_pA^T\Sigma_{\boldsymbol{\eta}}^{-1}\mathbf{y},\quad \Sigma_p =(A^T\Sigma_{\boldsymbol{\eta}}A+C_0^{-1})^{-1},
\]
where $A\in \mathbb{R}^{n\times d}$ is the discretization of the forward model $\mathcal{A}$, and the transport map is linear and available in closed form:
\begin{align}
\label{initial map}
T(\mathbf{x})=\mathbf{z}_0+Z_0\mathbf{x},
\end{align}
where $Z_0C_0Z_0^{T}=\Sigma_p$ and $\mathbf{z}_0=\mu_p$.
The initial point of the optimization problem \eqref{d_opt_problem} is chosen as this linear map \eqref{initial map}, i.e. $F_0=[\mathbf{z}_0\;\; Z_0]^T$ for the alternating direction method in Algorithm \ref{alg:ad}.
In fact, the optimization problem can become unconstrained depending on the parameterization of the transport map (e.g. integrated exponential parametrization or squared-integrated parameterization, see \cite{Marzouk_2016, Peherstorfer_2019} for more details) that the constraint can be automatically satisfied.

\begin{algorithm}[htbp]
\KwIn{Target density $\pi^{\mathbf{y}}(\mathbf{x},\lambda)$, Gamma distribution parameters $k$ and $\vartheta$, samples $\{\mathbf{x}^{(1)},\mathbf{x}^{(2)},\dots,\mathbf{x}^{(M)}\}$ from reference Gaussian measure $\mu_{ref}$, total number of iterations $K$.}
\KwOut{Direct transport map $\widetilde{T}$.}
Choose the initial coefficients $F_0$\;
\For{$k \gets 0 \ldots K$}{
Compute the parameter $\lambda_{k+1}$ by
$$
\lambda_{k+1}=\frac{2(k-1)}{\bigl\|\tfrac{1}{M}\sum^M_{i=1}\widetilde{T}(F_k; \mathbf{x}^{(i)})\bigr\|_{TV^{\alpha}}+2\vartheta};
$$

Update the coefficients $F_{k+1}$ by solving the optimization problem \eqref{d_opt_problem} for fixed $\lambda_{k+1}$\;
}
\Return{Direct transport map $\widetilde{T}$.}
\caption{Alternating direction method}
\label{alg:ad}
\end{algorithm}

\subsection{Linear diagonal map-based independence sampler}
We expect that the diagonal transport map approximately pushforward the reference onto the posterior, it is reasonable to consider an independence sampler for the diagonal transport map together with the reference distribution serves as a proposal distribution.
And the independence sampler with the proposal that is close to target space through a transport map can greatly reduce integrated autocorrelation time during sampling, which can clearly improve sampling efficiency.

In this paper, the approximation space for the map $\widetilde{T}$ is taken by the first-order polynomials (linear), i.e. $q=1$ in \eqref{map_pare}, which is able to catch the main information of posterior and be rapidly calculated for the numerical optimization.
Once we have a valid approximation of a transport map $\widetilde{T}$ between reference measure and target measure at hand, then the regularization parameters  is determined by the \eqref{de_reg}.
Moreover, the proposed distribution is set as the $\widetilde{T}_{\sharp}\mu_{ref}$, which is also a Gaussian measure, dented by $\mathcal{N}(m_1,C_1)$, because of the linear diagonal map, and the $\mu_{ref}$ can be taken by any Gaussian measures in the sampling process.
Obviously, the $\widetilde{T}_{\sharp}\mu_{ref}$ and $\mu_0$ are equivalent (see \cite{2018_Feng} for more details).
Further, if we also take the Gaussian measure $\mathcal{N}(m_1,C_1)$ as our Gaussian prior $\mu_0$ in the sampling process, then the proposal distribution $\widetilde{T}_{\sharp}\mu_{ref}$ is apparently reversible, and hence invariant, with respect to this $\mu_0$ \cite{Dashti_2015}.
Then, we can get
\begin{align}
\label{acc_is}
\frac{\mathrm{d}\mu^{\mathbf{y}}}{\mathrm{d}\mu}(v)\frac{\mathrm{d}\mu}{\mathrm{d}\mu^{\mathbf{y}}}(u)=\Phi_{\mathbf{y}}(u) + J(u;\lambda) - \Phi_{\mathbf{y}}(v) - J(v;\lambda),
\end{align}
for the independence sampler.

In this work, the complete linear diagonal map-based independence sampler is summarized in Algorithm \ref{alg:diagonal_map}.
And in Algorithm \ref{alg:diagonal_map}, we set the initial state $u^{(1)}=\mathbb{E}_{\pi_{ref}}[\widetilde{T}(\mathbf{x})]$ that belongs to  approximation target space, which can accelerate the convergence of the Markov chains.
Although the regularization parameter $\lambda$ is determined in advance, as will be demonstrated with our numerical results in Section \ref{Numerical examples}, the additional FTG prior still has an outstanding performance in terms of handling the complex structure and detail information of model inputs.

\begin{algorithm}[htbp]
\KwIn{Direct transport map $\widetilde{T}$ obtained by the alternating direction method in Algorithm \ref{alg:ad}, reference distribution $\mu_{ref}$, the data-misfit function $\Phi_{\mathbf{y}}(u)$, the prior $J(u;\lambda)$, total number of steps $K$.}
\KwOut{MCMC samples of the target distribution, $\{u^{(i)}\}_1^K$.}
Compute the initial state $u^{(1)}=\mathbb{E}_{\pi_{ref}}[\widetilde{T}(\mathbf{x})]$ using the SAA\;
Compute the regularization parameter $\lambda$ via the $\eqref{de_reg}$\;
\For{$i \gets 1 \ldots K-1$}{
Propose $\hat{r}$ from $\mu_{ref}$\;
Pushforward proposed state $\hat{r}$ onto target, $v=\widetilde{T}(\hat{r})$\;
Calculate the acceptance probability $a(u^{(i)},v)$ using the $\eqref{acc}$:
\[
a(u^{(i)}, v)=\min \Bigl\{1, \Phi_{\mathbf{y}}(u^{(i)}) + J(u^{(i)};\lambda) - \Phi_{\mathbf{y}}(v) - J(v;\lambda) \Bigr\};
\]

Set $u^{(i+1)}=v$ with probability $a(u^{(i)},v)$; else set $u^{(i+1)}=u^{(i)}$\;
}
\Return{\textrm{Target samples} $\{u^{(i)}\}_1^K$.}
\caption{Linear diagonal map-based independence sampler}
\label{alg:diagonal_map}
\end{algorithm}

\section{Numerical examples}\label{Numerical examples}\label{section4}
In this section, we present some examples and applications to demonstrates the performance of our FTG prior and the linear diagonal map-based independence sampler.
These applications and examples for the linear inverse problems range from the one-dimensional deconvolution problem and inverse source identification problem to the two-dimensional limited computed tomography reconstruction in medical imaging and the image denoising.

We use the $Gr\ddot{u}nwald$ formula to discretize the Riemann-Liouville fractional derivative.
If $0<\alpha\le1$, i.e. $n=1$, the Riemann-Liouville fractional derivative is approximated by the following standard $Gr\ddot{u}nwald$ formula
\begin{align}
\nabla^{\alpha}u(x_l)=\frac{1}{2h^{\alpha}}\Biggl(\sum^l_{j=0}\varpi^{\alpha}_ju_{l-j}-\sum^{d-l}_{j=0}\varpi^{\alpha}_ju_{l+j}\Biggr),
\end{align}
and if $1<\alpha\le2$, i.e. $n=2$, the shifted $Gr\ddot{u}nwald$ formula
is used
\begin{align}
\nabla^{\alpha}u(x_l)=\frac{1}{2h^{\alpha}}\Biggl(\sum^{l+1}_{j=0}\varpi^{\alpha}_ju_{l-j+1}-\sum^{d-l+1}_{j=0}\varpi^{\alpha}_ju_{l+j-1}\Biggr),
\end{align}
where $l=1,2,\dots,d-1$, and $\varpi^{\alpha}_0=1, \varpi^{\alpha}_j=(1-\tfrac{\alpha+1}{j})\varpi^{\alpha}_{j-1},\, j=1,2,\dots,d$.
In addition, if no special explanation is given, we choose the following exponential covariance for the Gaussian prior $\mu_0$,
\begin{align}
C_0(x_1,x_2)=\gamma\exp\Biggl[-\frac{1}{2}\Bigl(\frac{x_1-x_2}{\nu}\Bigr)^2\Biggr],
\end{align}
where $\gamma$ and $\nu$ are the parameters.

To obtain the transport map, in this paper, we use the interior point algorithm to solve the \eqref{d_opt_problem} for fixed $\lambda$ in Algorithm \ref{alg:ad}.
If no special explanation is given, the numerical optimization problem is performed with MATLAB's fmincon optimizer, where the step tolerance (StepTolerance) is set to $10^{-6}$; the SpecifyConstraintGradient and SpecifyOdjectiveGradient are set to true (we provide the gradients for the solver in the linear inverse problems, see the \cite{ELMOSELHY20127815} for more detail).
And the reference distribution $\mu_{ref}$ is set the same as the Gaussian distribution $\mu_0$ in solving the numerical optimization process, and the chain was run for $1\times10^5$ steps for our diagonal map-based independent sampler in the Algorithm \ref{alg:diagonal_map} and then all sample is used to calculate the posterior mean.

We calculate the average relative error (RelErr) to measure the difference between the reconstructed result (the posterior mean) and the target.
For the posterior mean $\mathbf{x}$ and the target $\mathbf{x}_0$, RelErr is defined as
\begin{align}
RelErr=\frac{\|\mathbf{x}-\mathbf{x}_0\|_2}{\|\mathbf{x}_0\|_2}.
\end{align}
\subsection{Deconvolution problem}\label{1d_ex1}

\subsubsection{Problem setup}
Consider the Fredholm first kind integral equation \cite{2002_Vogel} of convolution type:
\begin{align}
\label{decon}
(\mathcal{A}f)(x):=\int_{\Omega}A(x-y)f(y)\mathrm{d}y=g(x),\quad x\in\Omega.
\end{align}
This two-dimensional version of a model that occurs in optical imaging.
In this application, $g$ stands for the blurred image intensity, $f$ represents light source intensity and the kernel $A$ characterizes blurring effects that occur during image formation.
In this example, we consider the one-dimensional version of above model, i.e. $\Omega=[l,r]$, and its  kernel is
\begin{align}
A(x)=\xi\exp\Bigl(-\frac{x^2}{2\delta}\Bigr),
\end{align}
where the $\xi$ and $\delta$ are positive parameters.
Obviously, we can use the standard numerical quadrature to get the accurate approximation of $\mathcal{A}f=g$.
Then the midpoint quadrature is applied to discretize the equation \eqref{decon} and a discrete linear system $A\mathbf{f}=\mathbf{g}$ is obtained, where the $A=(a_{ij})_{d\times d}$ and
\begin{align}
a_{ij}=h\,\xi\exp\Bigl[-\frac{\bigl((i-j)h\bigr)^2}{2\delta^2}\Bigr],
\end{align}
where $h=(r-l)/d$.
And the observations are obtained by the $\mathbf{y}=A\mathbf{f}+\boldsymbol{\eta}$, where the $\boldsymbol{\eta}$ is the Gaussian white noise.

\subsubsection{Set up of inverse problems}
In order to illustrate get the numerical results, we take $\xi=1/(\delta\sqrt{2\pi}),\,\delta=0.02$ and $\Omega=[0,1]$.
We consider the light source
 \[ f(x)=
\begin{cases}
0.5 &  0.1\le x<0.25,\\
0.25 &  0.35\le x<0.4,\\
\sin^4(2\pi x) & 0.5\le x<1,\\
0 &  otherwise,
\end{cases}
\]
and set $d=120$.
The measurement data are obtained by $\mathbf{y}=A\mathbf{f}+\boldsymbol{\eta}$ and $\boldsymbol{\eta}$ is assumed to be the Gaussian noise with zero mean and standard deviation $4.84\times10^{-2},\,9.7\times10^{-3},\,4.8\times10^{-3}$, which corresponds to $5\%,\,1\%,\,0.05\%$ noise, respectively, with respect to the maximum norm of the output $A\mathbf{f}$.
Note that the measurement data $\mathbf{y}$ is computed from a twice finer grid.
The Gaussian prior $\mu_0$ is taken by zero mean and set $\gamma=0.016, \nu=0.0003$ and the shape and rate parameter of Gamma distribution are set to $k=2\times10^3,\;\vartheta=1$, respectively.

We then construct a linear diagonal transport map from the Gaussian distribution $\mu_0$ to the posterior \eqref{H_postprior} using the alternating direction method in Algorithm \ref{alg:ad}.
We use $M=1000$ samples of the Gaussian distribution $\mu_0$ to approximate the expected value via the SAA in the objective function (see Section \ref{numercal opt}).
The transport map $\widetilde{T}$ and the regularization parameter $\lambda$ are then used as precondition for the independent sampler as shown in Algorithm \ref{alg:diagonal_map}.
And the proposal $\mu_{ref}$ is the Gaussian distribution with zeros mean and the standard deviation $2\times10^{-3}$ in Algorithm \ref{alg:diagonal_map}.

\subsubsection{Result}

\begin{figure}[htbp]
 \centering
 \begin{tabular}{@{}c@{}c@{}}
\includegraphics[width=0.45\textwidth, height=0.25\textwidth]{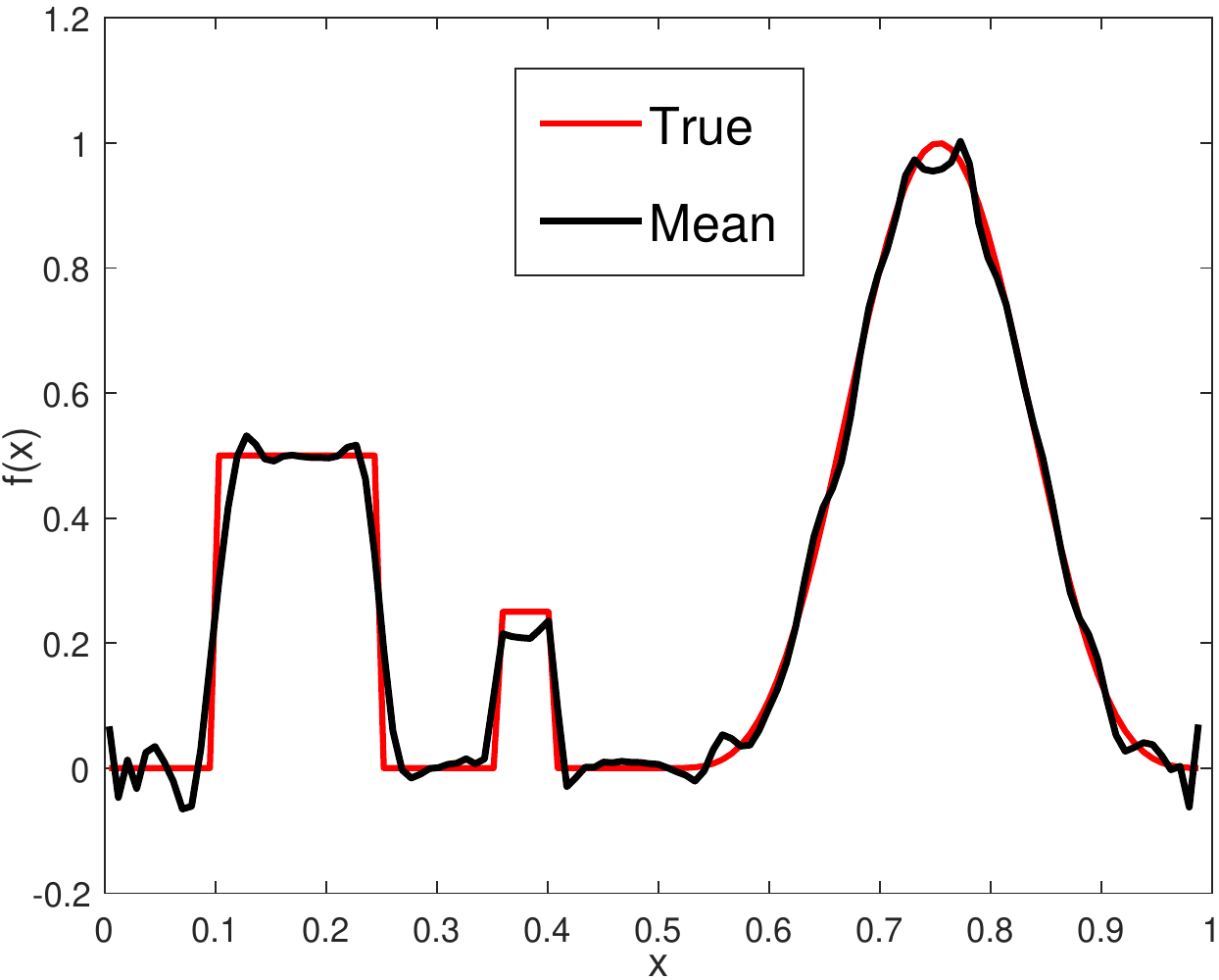}&
\includegraphics[width=0.45\textwidth, height=0.25\textwidth]{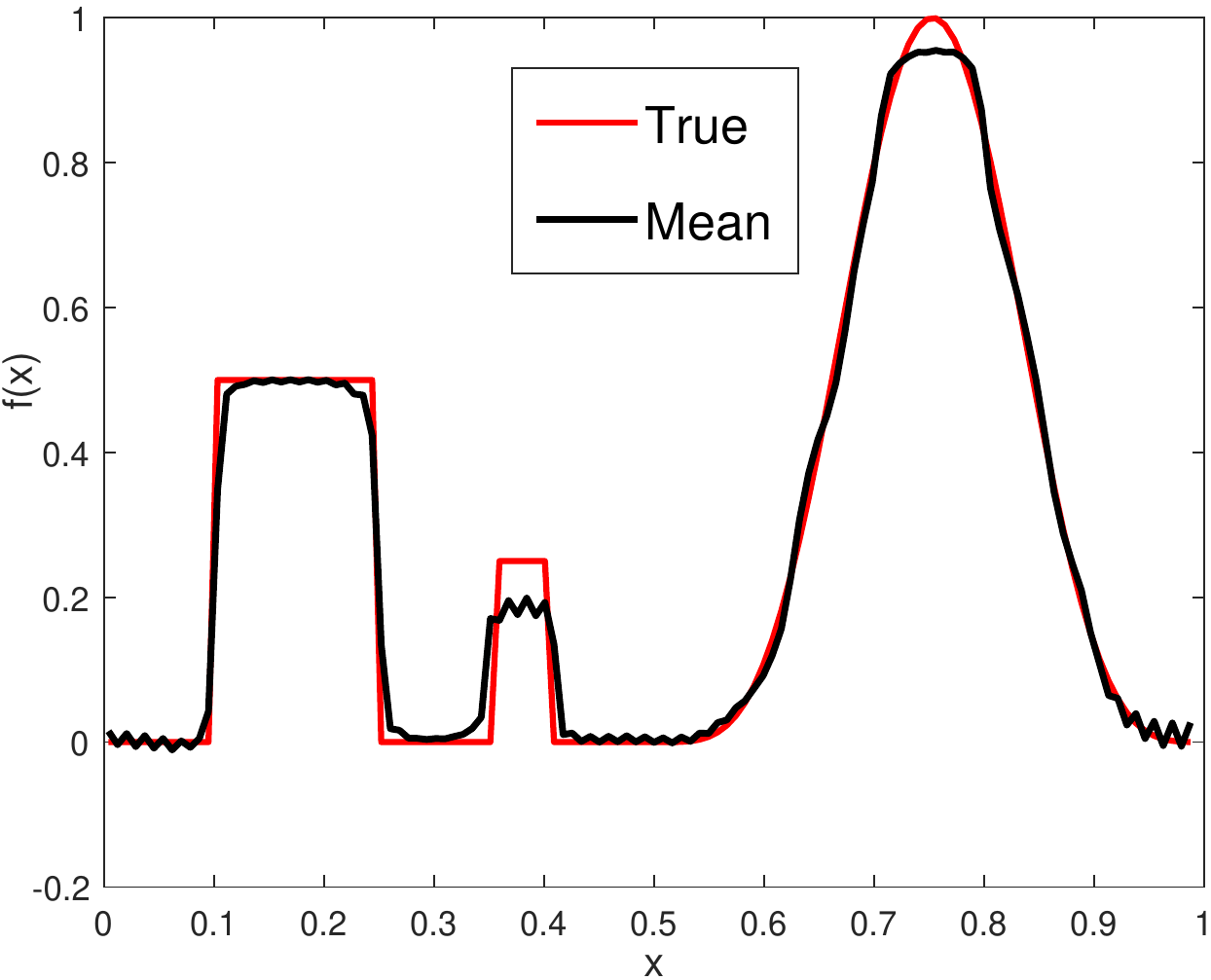} \\
\includegraphics[width=0.45\textwidth, height=0.12\textwidth]{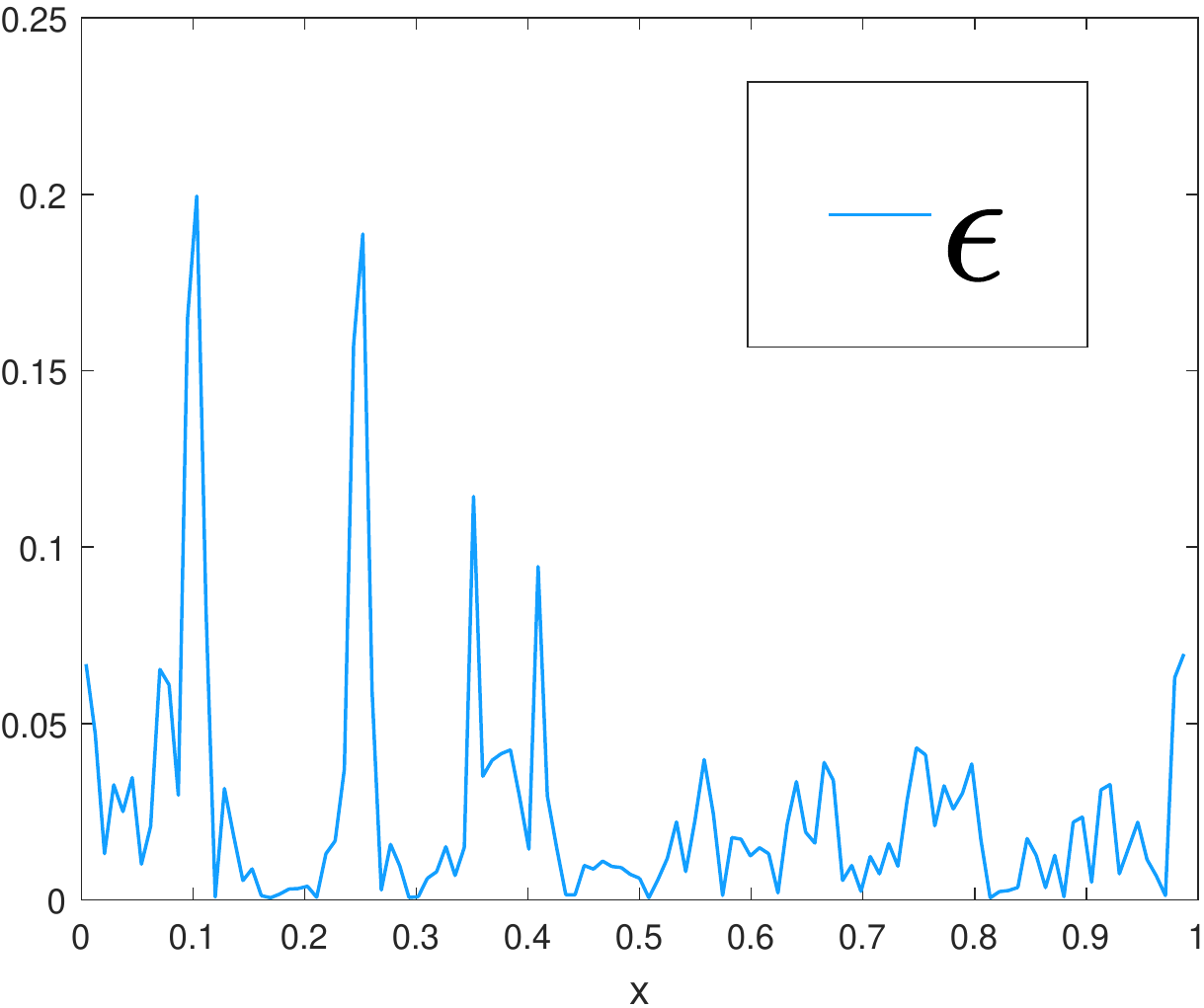}&
\includegraphics[width=0.45\textwidth, height=0.12\textwidth]{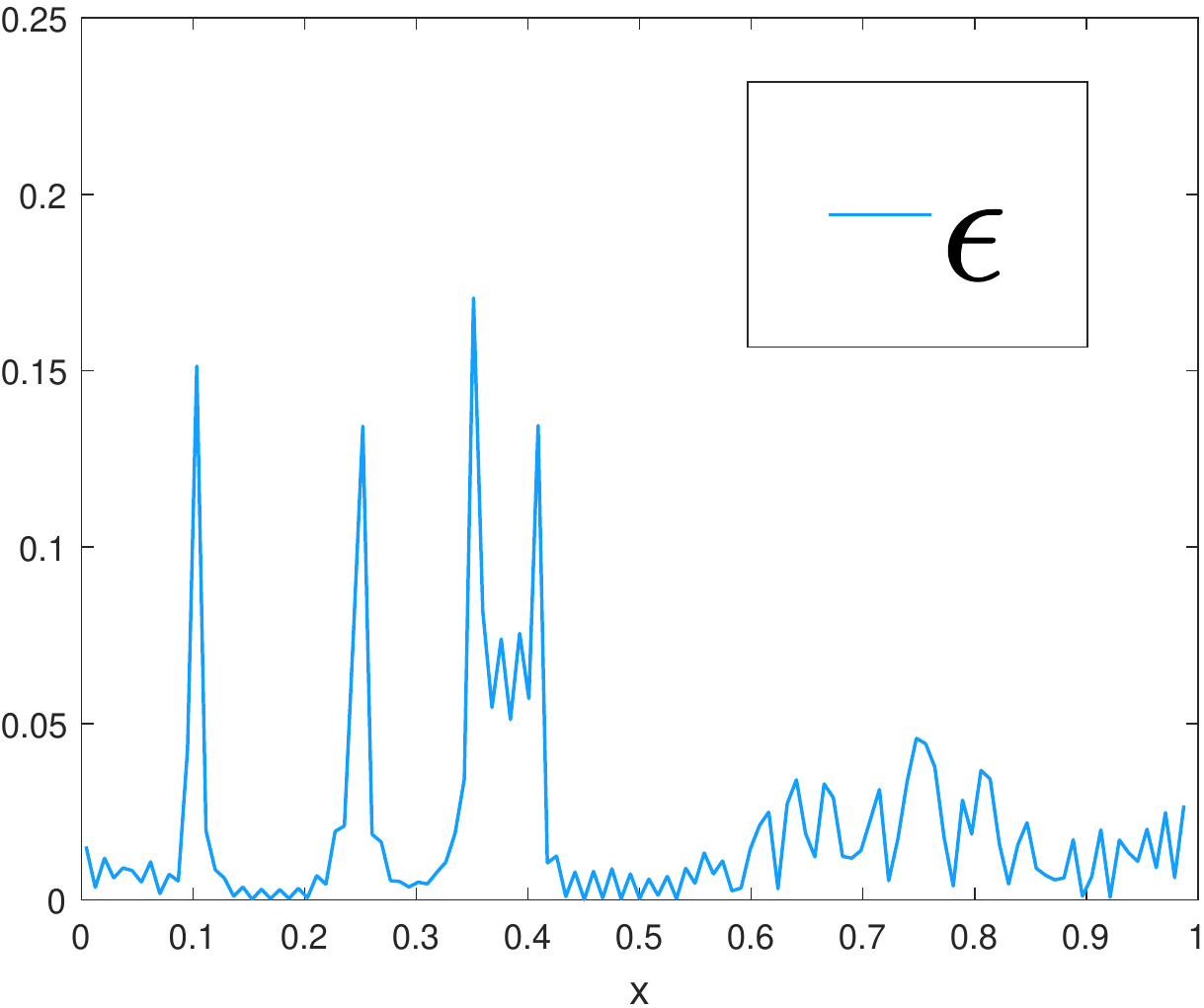}  \vspace{-5pt}\\
FTG with $\alpha=0.5$&FTG with $\alpha=0.9$\\  \hline
\includegraphics[width=0.45\textwidth, height=0.25\textwidth]{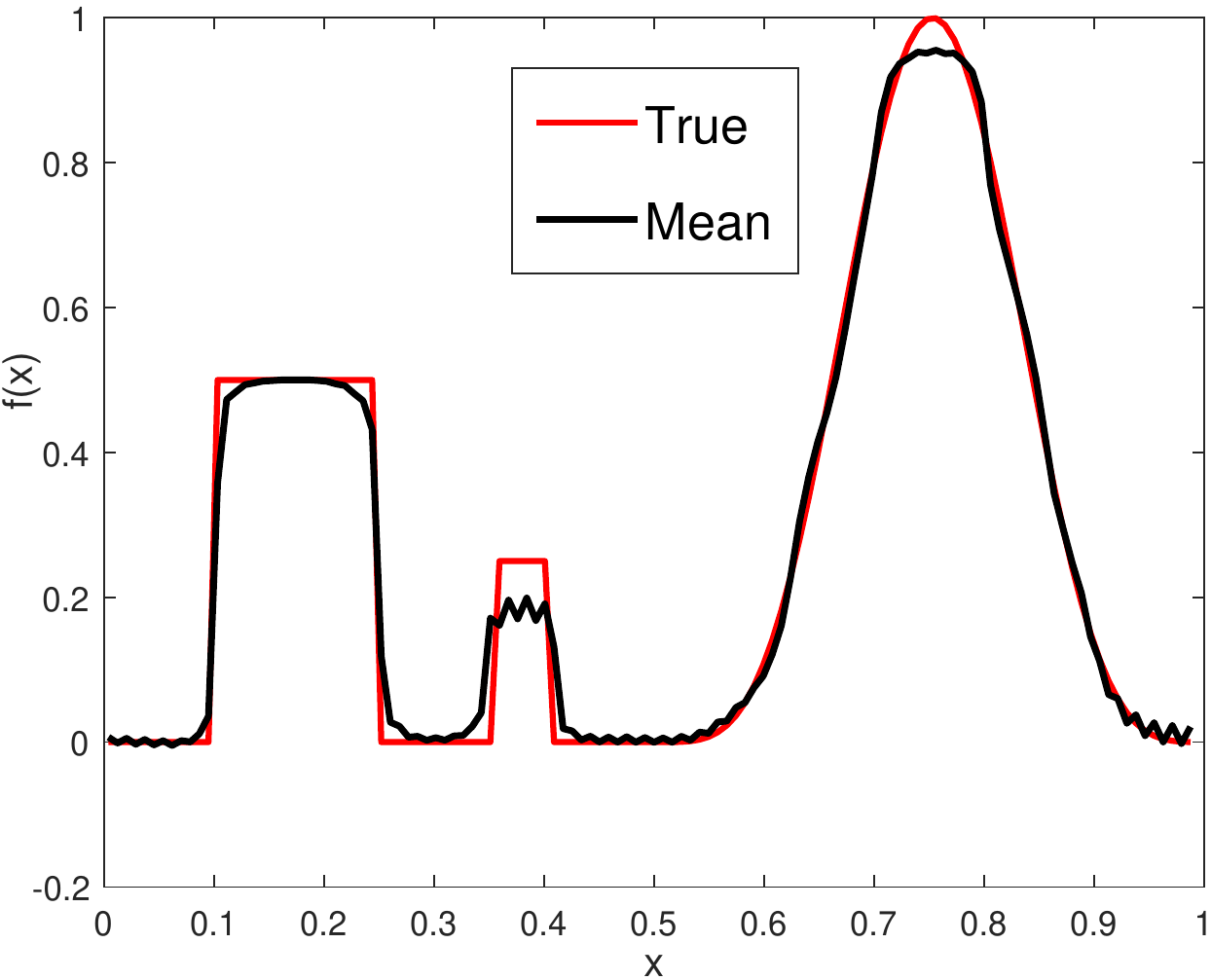} \hspace{2pt}&
\includegraphics[width=0.45\textwidth, height=0.25\textwidth]{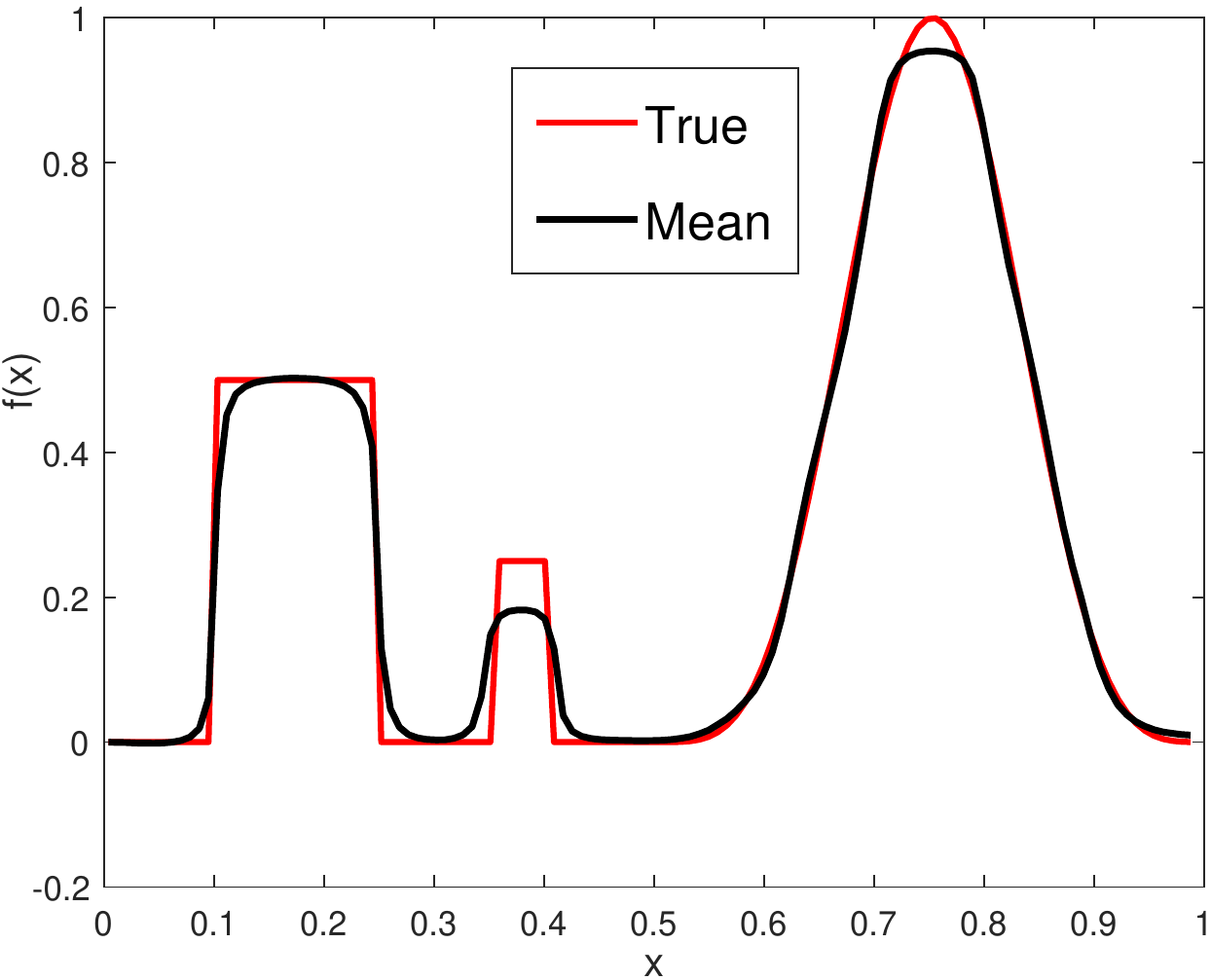}\\
\includegraphics[width=0.45\textwidth, height=0.12\textwidth]{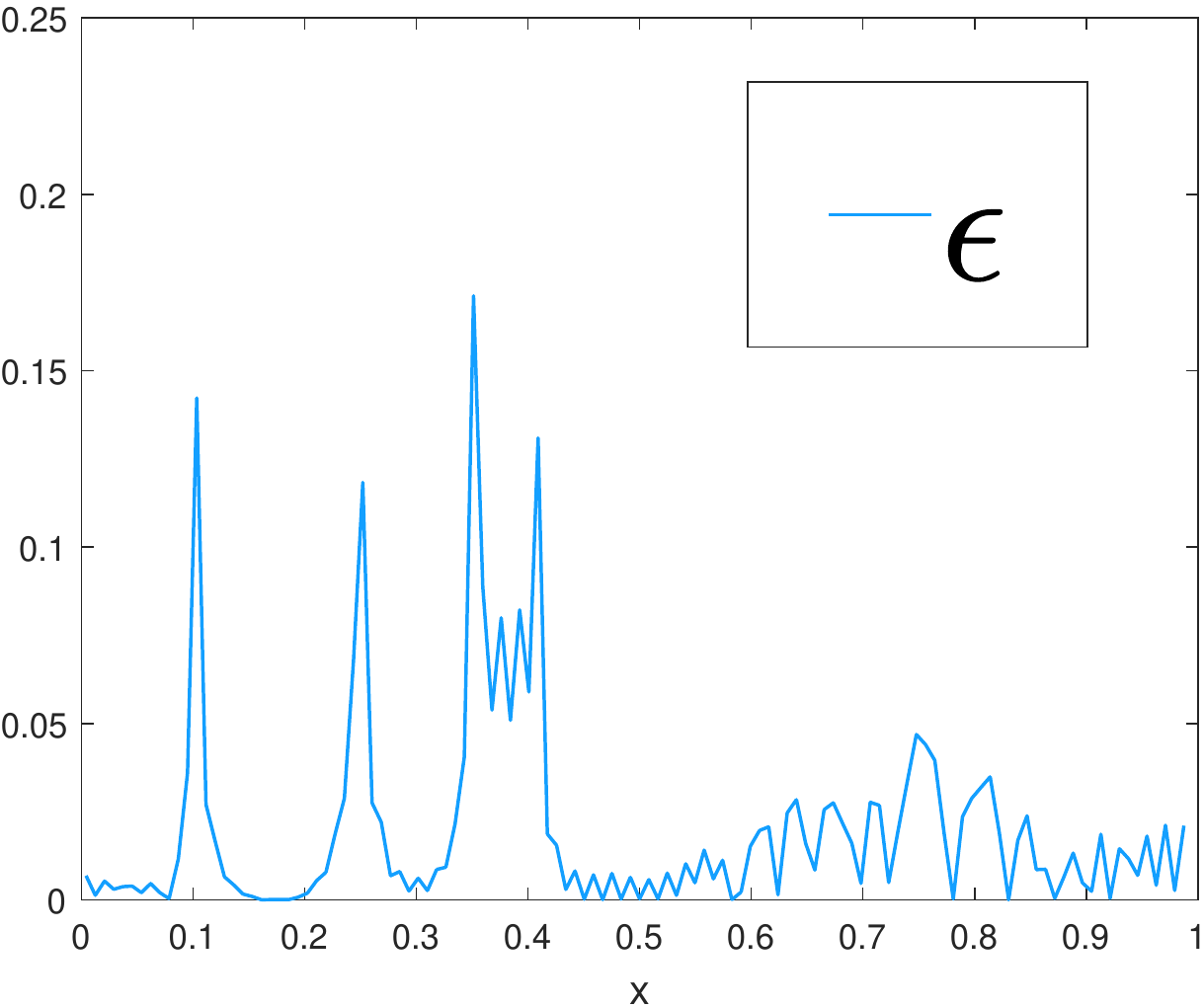}&
\includegraphics[width=0.45\textwidth, height=0.12\textwidth]{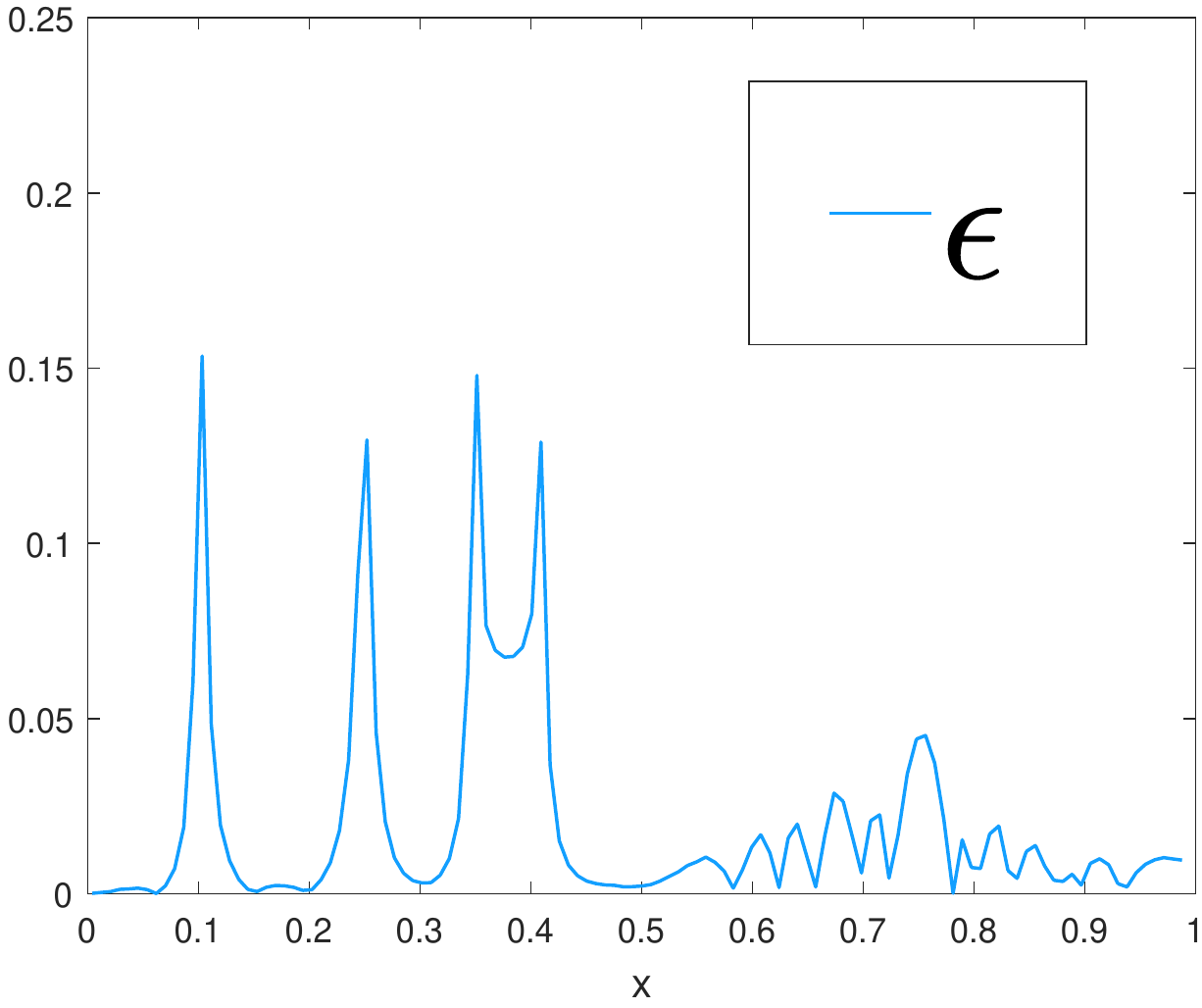}  \vspace{-5pt}\\
 FTG with $\alpha=0.95$ &TG \\ \hline
 \end{tabular}
\caption{Reconstruction results for the deconvolution $f(x)$.
The posterior mean and absolute error $\epsilon$ using FTG prior with $\alpha=0.5,\;0.9,\;0.95$ and TG prior for $1\%$ noise level.
}
\label{decon_result_01}
\end{figure}

We first compare the performance of FTG prior to that of TG prior in two cases.
For  $0<\alpha\le1$, the posterior mean and absolute error $\epsilon$ with TG prior and with FTG prior are plotted in Figure \ref{decon_result_01} and the relative error $RelErr$ are listed in Table \ref{table_decon_01}.
From this table and figure, in general, we can conclude that the reconstruction results using FTG prior gradually converges to that using TG prior if the fractional order $\alpha \to 1^-$.
However, from the Table \ref{table_decon_01}, our FTG prior with $\alpha=0.95,\;0.99$ have better numerical results compared with the TG prior.
And the reconstructed results with FTG prior for $\alpha=0.9,\;0.95$ outperforms that with TG prior and can catch some details such as the corner points, seeing the absolute error curve in Figure \ref{decon_result_01}.
Therefore if the reconstructed target has much this details information, as discussed the example CT reconstruction in Section \ref{2d_ex1}, the FTG prior with the $\alpha \to 1^-$ can obtain better recovery result, compared with the TG prior.
In addition, we can find that the FTG with $\alpha=0.2,\,0.5$ have similar the $RelErr$ value in Table \ref{table_decon_01}.
It seems that our Gaussian prior $\mu_0$ is proper for this example in the process of solving numerical optimization.
This is because the Gaussian prior $\mu_0$ will play main role in the numerical results for the a small fractional order FTG prior.

\begin{table}[htbp]
\centering
\caption{The $RelErr$ values of the deconvolution results using FTG prior with $0<\alpha\le1$ and TG prior for $1\%$ noise level.}
\begin{tabular}{cccccccccc}
  \toprule
          & TG & $\alpha=0.2$ & $\alpha=0.5$ & $\alpha=0.8$ & $\alpha=0.9$ & $\alpha=0.95$ & $\alpha=0.99$ \\
  \midrule
  $RelErr$ & 0.0836&0.1020& 0.1020& 0.0930 & 0.0837 & \textbf{0.0822} &\textbf{0.0828}\\
    \bottomrule
  \end{tabular}
  \label{table_decon_01}
\end{table}

\begin{figure}[htbp]
 \centering
 \begin{tabular}{@{}c@{}c@{}}
\includegraphics[width=0.45\textwidth, height=0.25\textwidth]{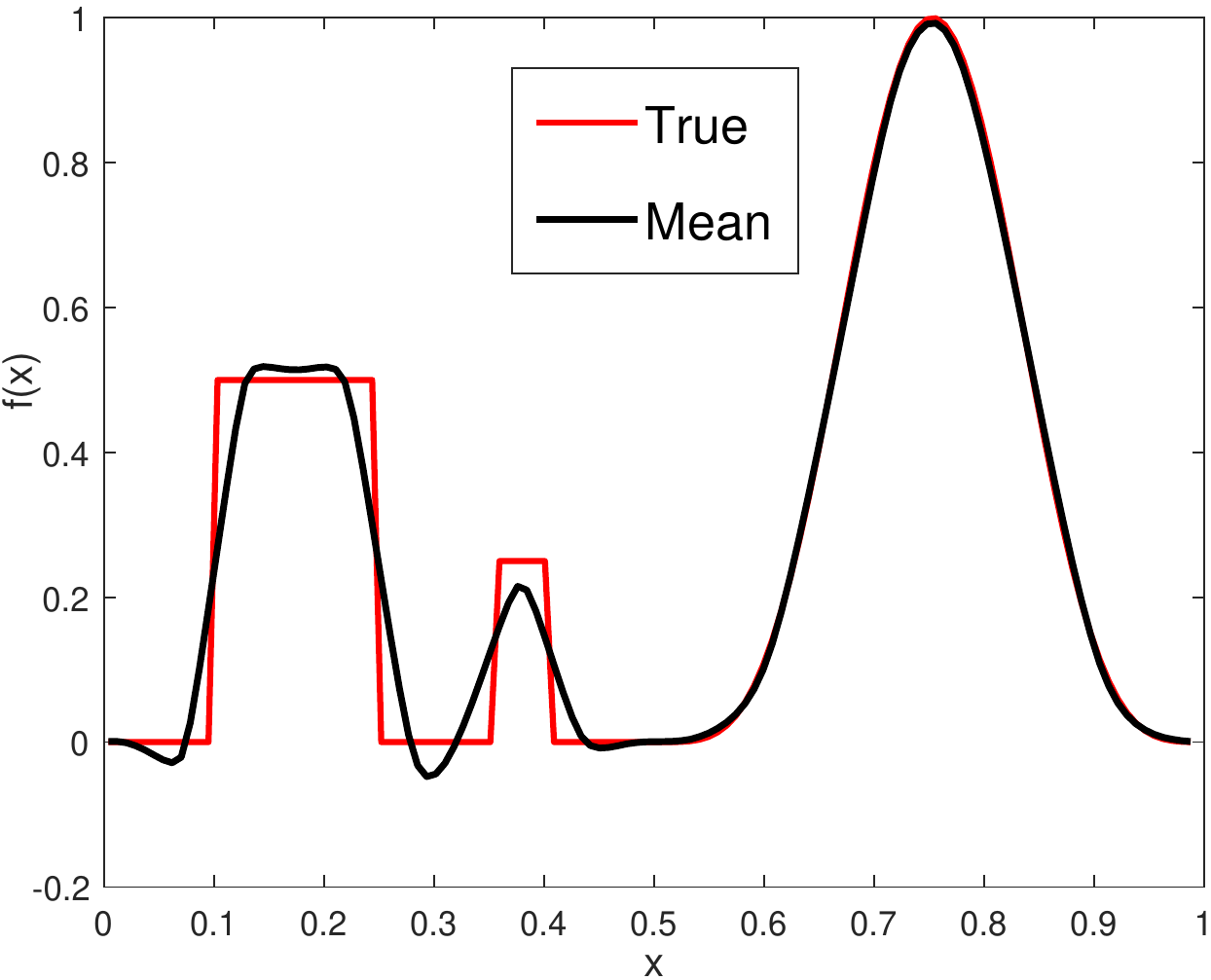}&
\includegraphics[width=0.45\textwidth, height=0.25\textwidth]{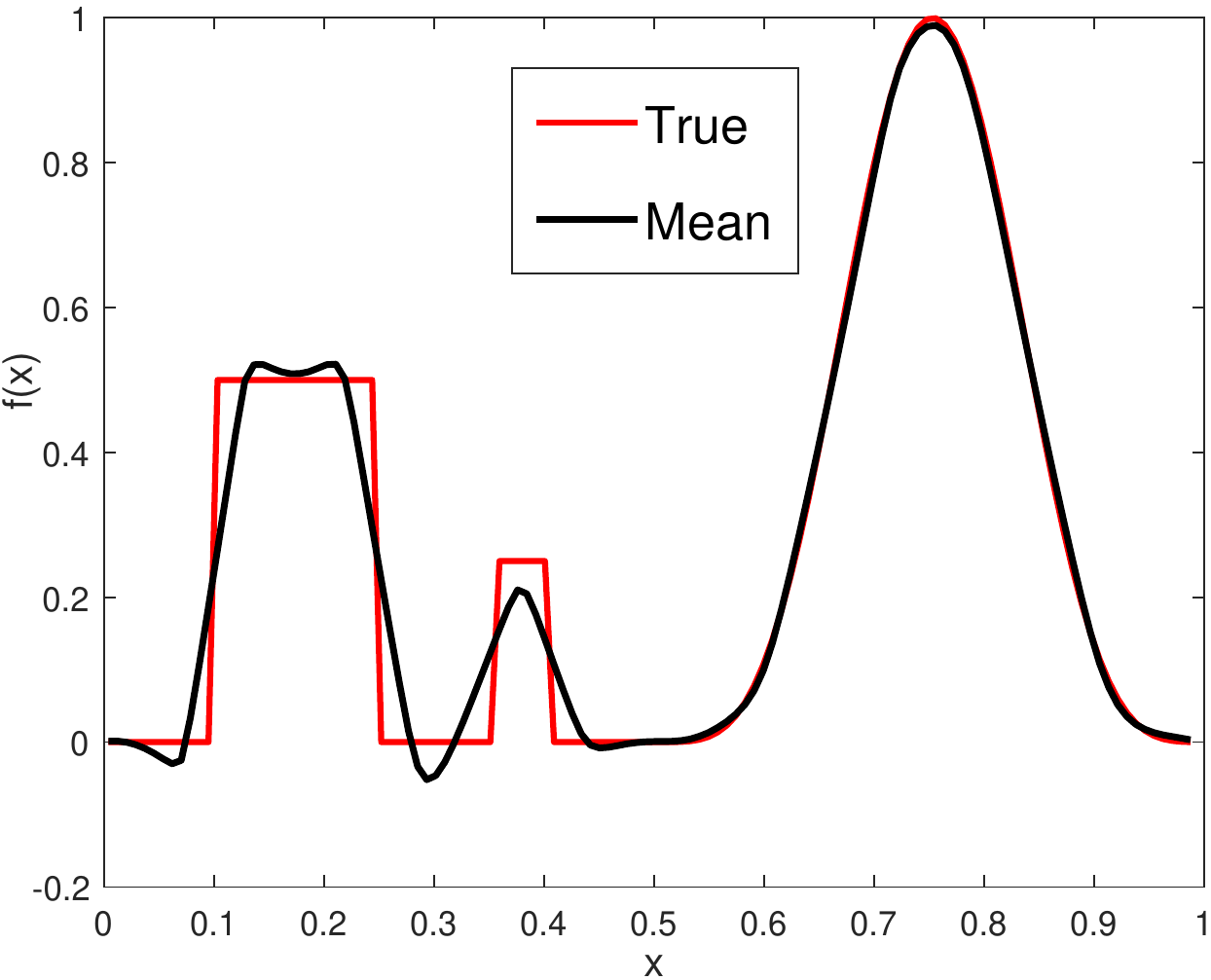} \\
\includegraphics[width=0.45\textwidth, height=0.12\textwidth]{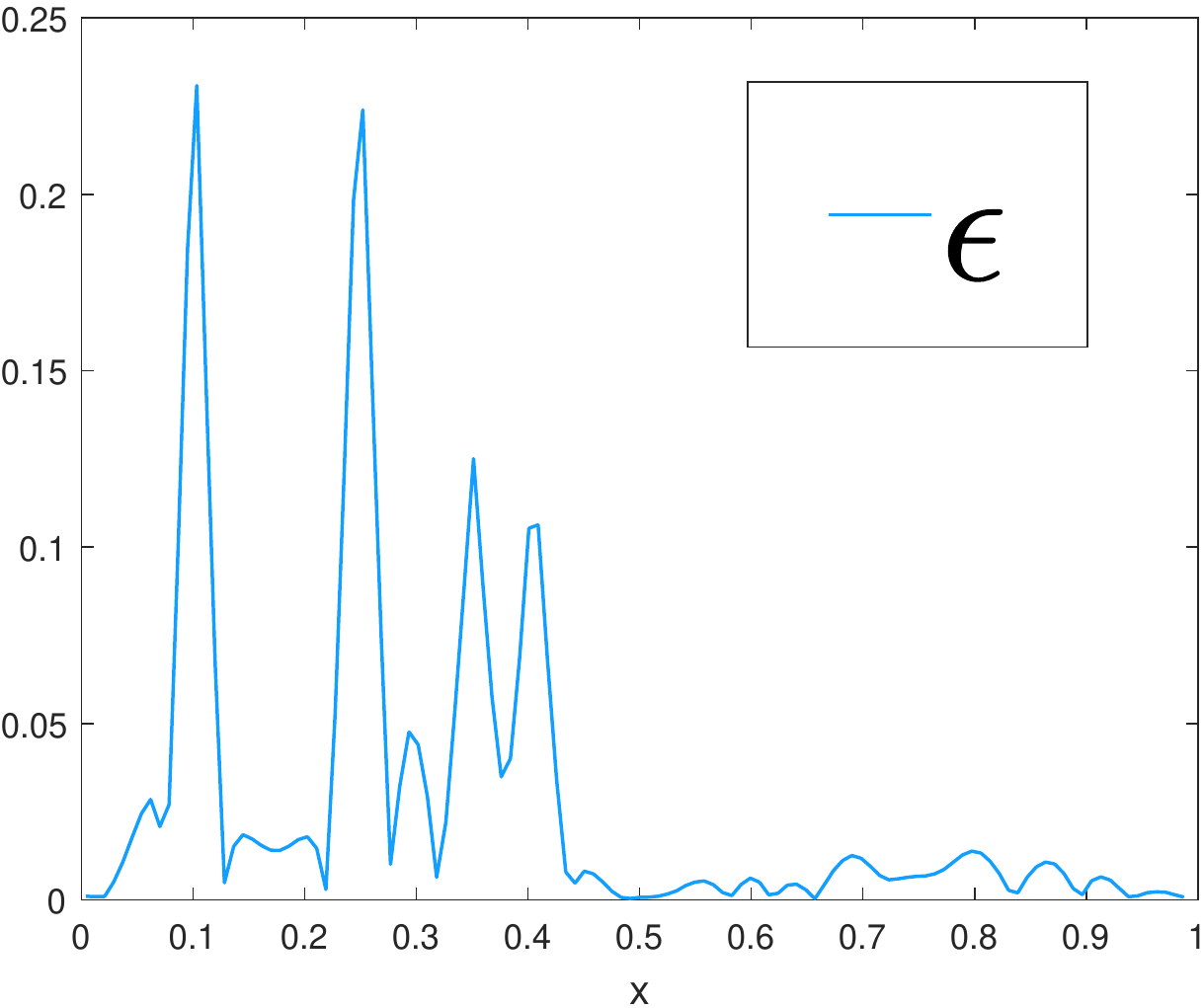}&
\includegraphics[width=0.45\textwidth, height=0.12\textwidth]{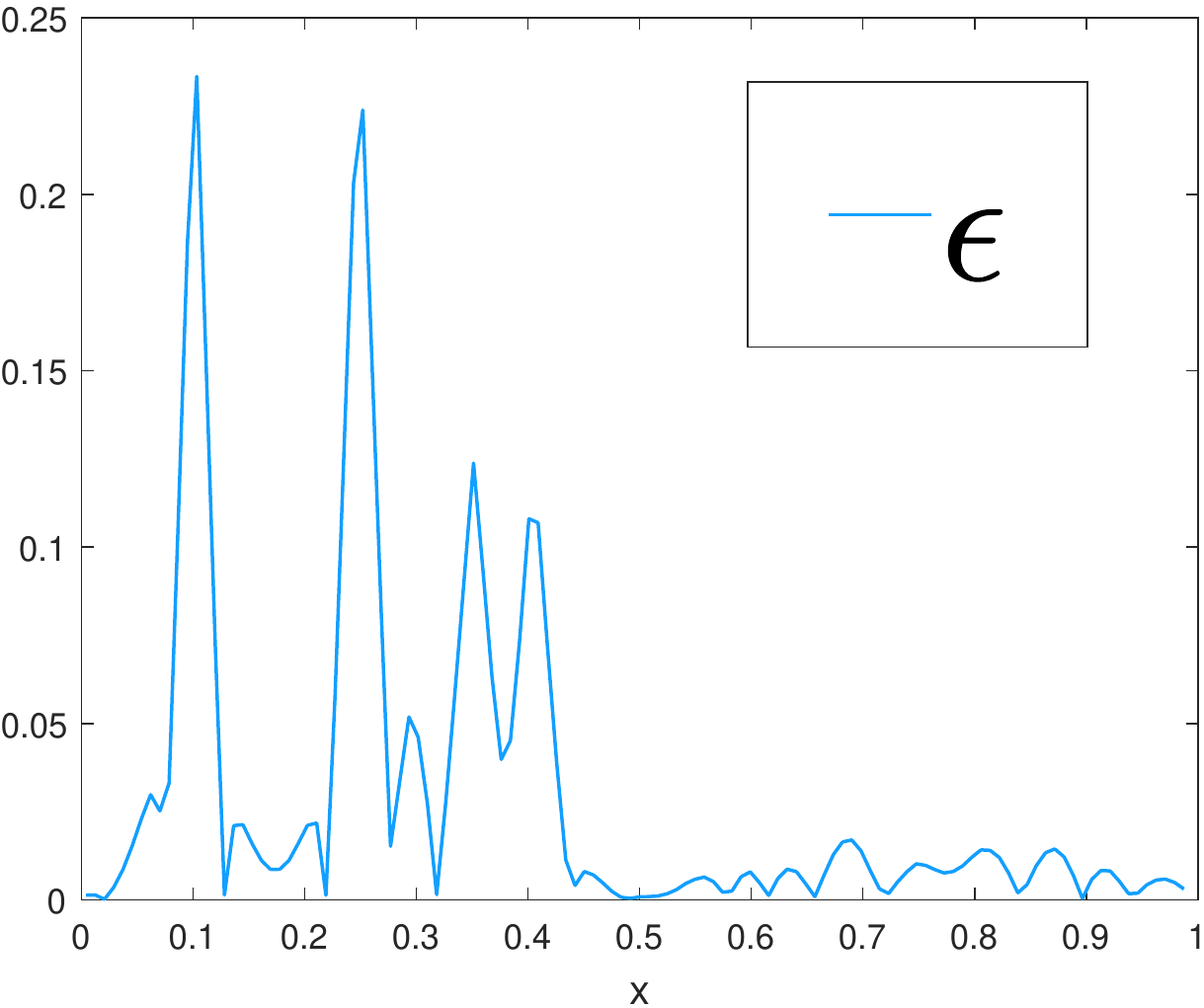}  \vspace{-5pt}\\
FTG with $\alpha=1.01$ &FTG with $\alpha=1.05$\\  \hline
\includegraphics[width=0.45\textwidth, height=0.25\textwidth]{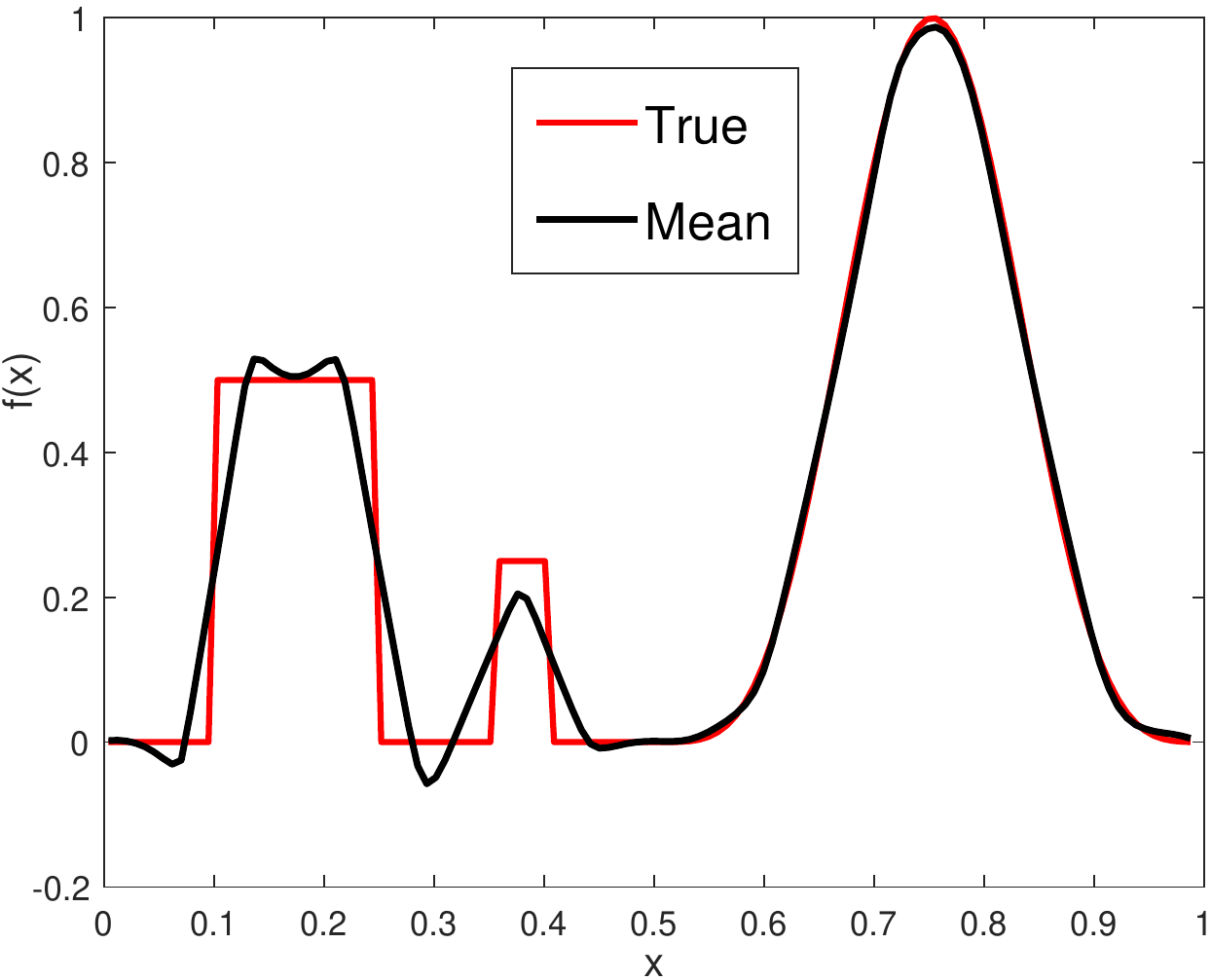} \hspace{2pt}&
\includegraphics[width=0.45\textwidth, height=0.25\textwidth]{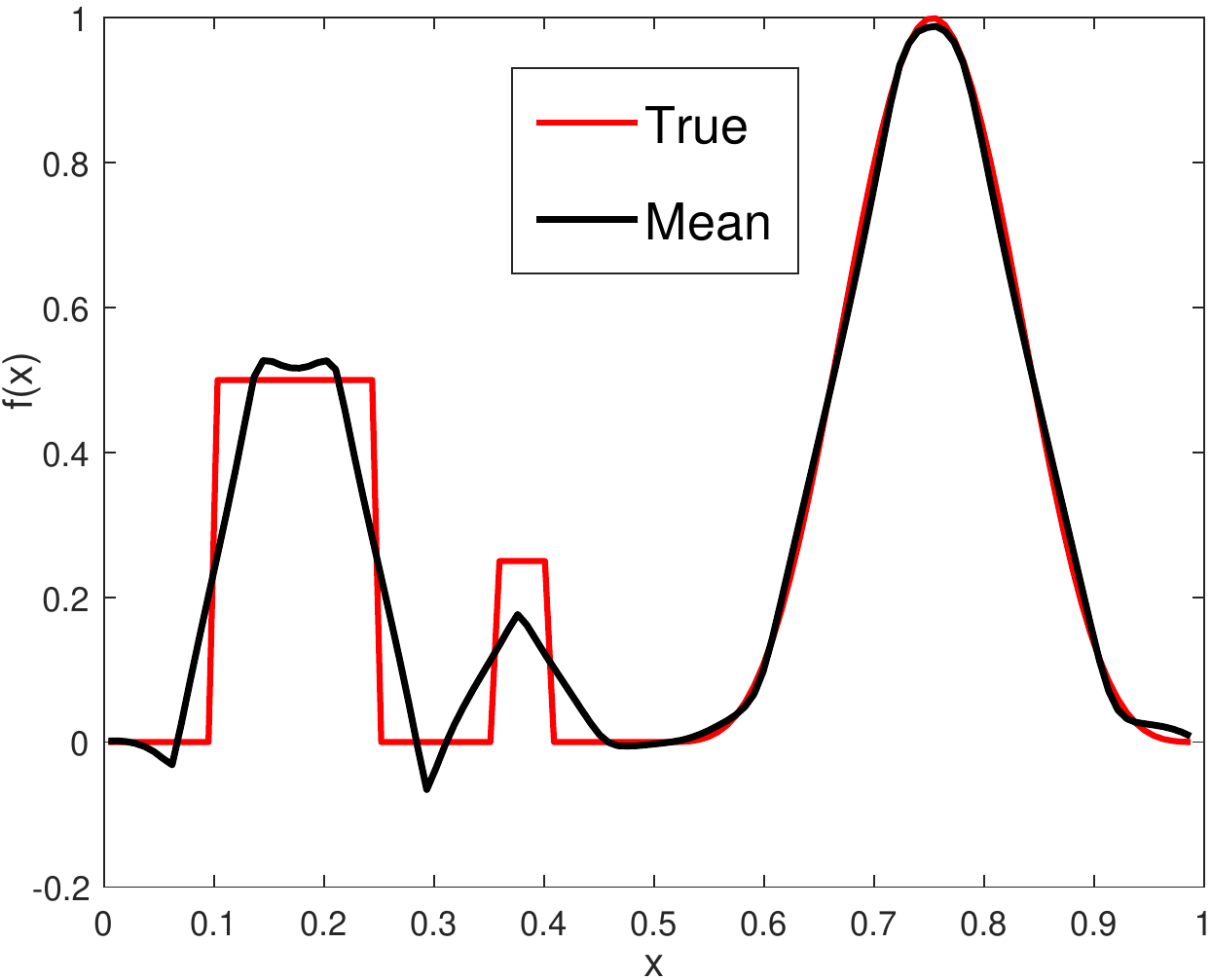}\\
\includegraphics[width=0.45\textwidth, height=0.12\textwidth]{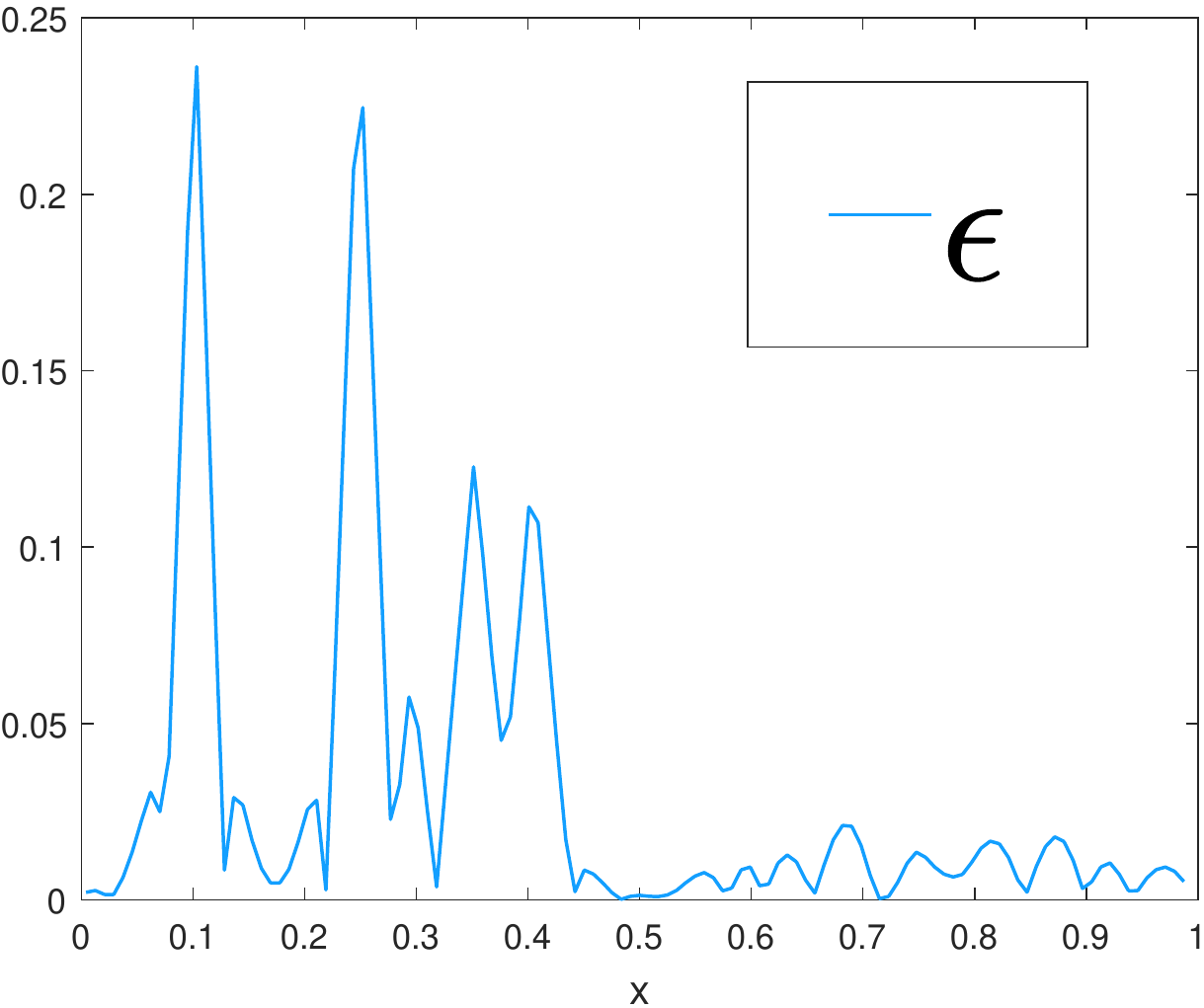}&
\includegraphics[width=0.45\textwidth, height=0.12\textwidth]{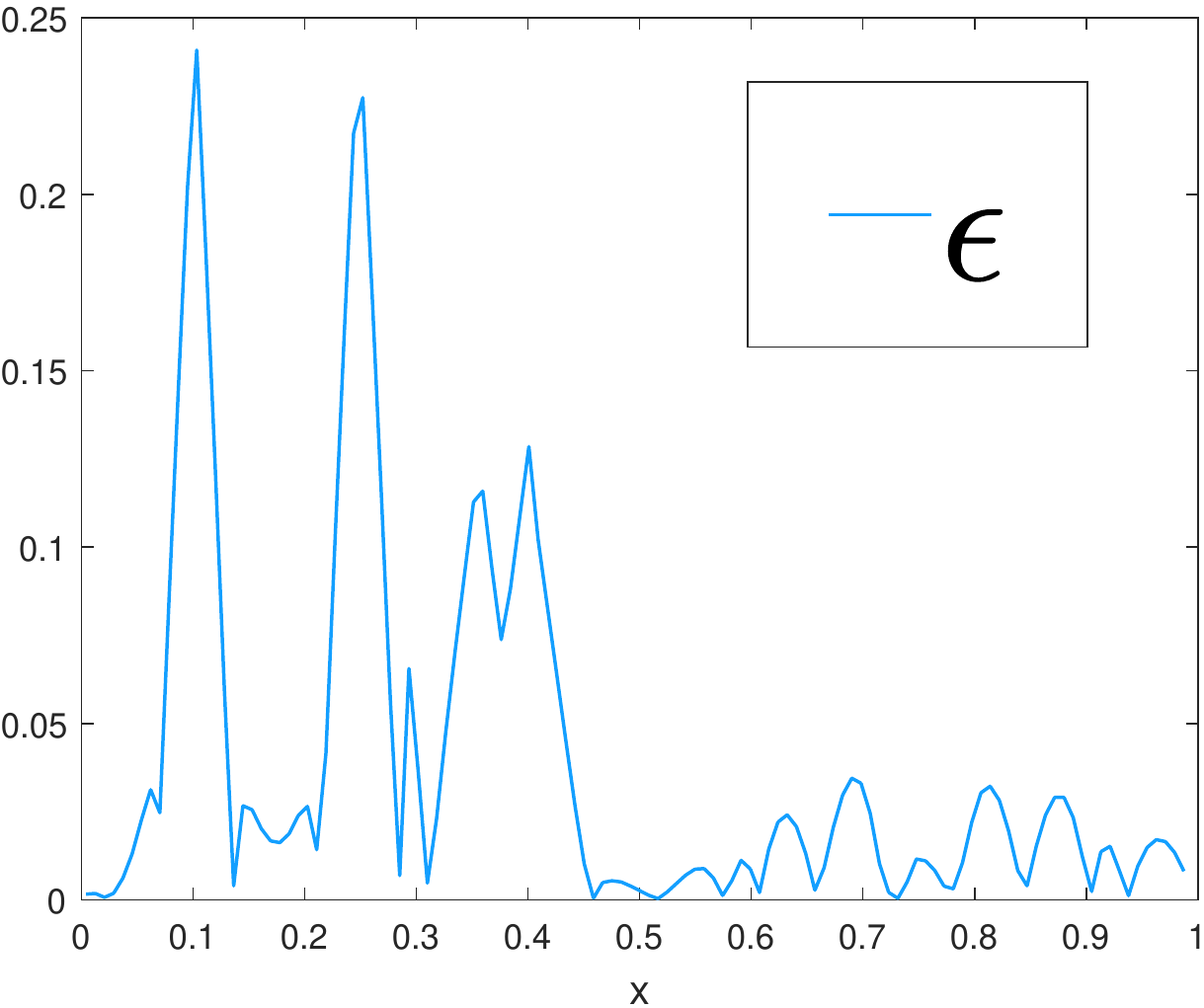}  \vspace{-5pt}\\
 FTG with $\alpha=1.1$ &FTG with $\alpha=1.5$ \\ \hline
 \end{tabular}
\caption{Reconstruction results for the deconvolution $f(x)$.
The posterior mean and absolute error $\epsilon$ using FTG prior with $\alpha=1.01,\;1.05,\;,1.1,\;1.5$ and TG prior for $1\%$ noise level.
}
\label{decon_result_12}
\end{figure}

The posterior mean and absolute error $\epsilon$ with FTG prior for $1<\alpha\le2$ are plotted in Figure \ref{decon_result_12} and the relative error $RelErr$ are listed in Table \ref{table_decon_12}.
As we can see, from the Table \ref{table_decon_12}, the TG prior yield a lower $RelErr$ values compared to the FTG results.
However, in the Figure \ref{decon_result_12}, our FTG prior with $\alpha=1.01,\; 1.05,\; 1.1$ can eliminate well the staircase effect in term of a smooth part in the reconstruction target $f(x)$ and also catch the piecewise constant structure.
Thus, if the reconstructed target has much oscillation information, our FTG prior with $\alpha \to 1^+$ can not only outperforms the TG, but also maintain an acceptable relative error.

\begin{table}[htbp]
\centering
\caption{The $RelErr$ values of the deconvolution results using FTG prior with $1<\alpha\le2$ and TG prior for $1\%$ noise level.}
\begin{tabular}{cccccccccc}
  \toprule
          & TG & $\alpha=1.01$ & $\alpha=1.05$ & $\alpha=1.1$ & $\alpha=1.2$ & $\alpha=1.5$ & $\alpha=1.8$ \\
  \midrule
  $RelErr$ & 0.0836&0.1271& 0.1309& 0.1355 &0.1427& 0.1537 & 0.1597\\
    \bottomrule
  \end{tabular}
  \label{table_decon_12}
\end{table}

\begin{figure}[htbp]
 \centering
 \begin{tabular}{@{}c@{}c@{}c@{}}
\includegraphics[scale=0.5]{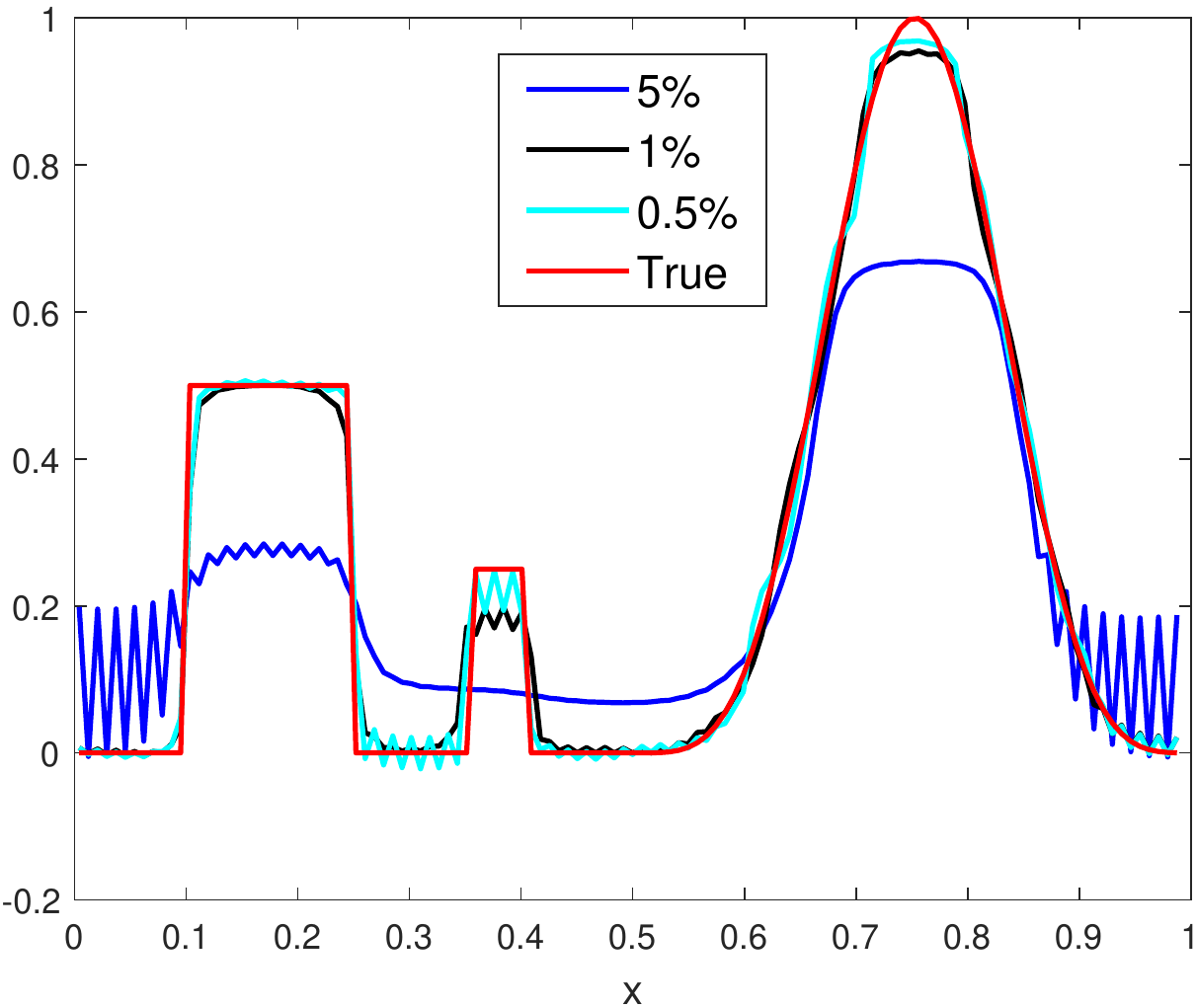}& \hspace{20pt}
\includegraphics[scale=0.5]{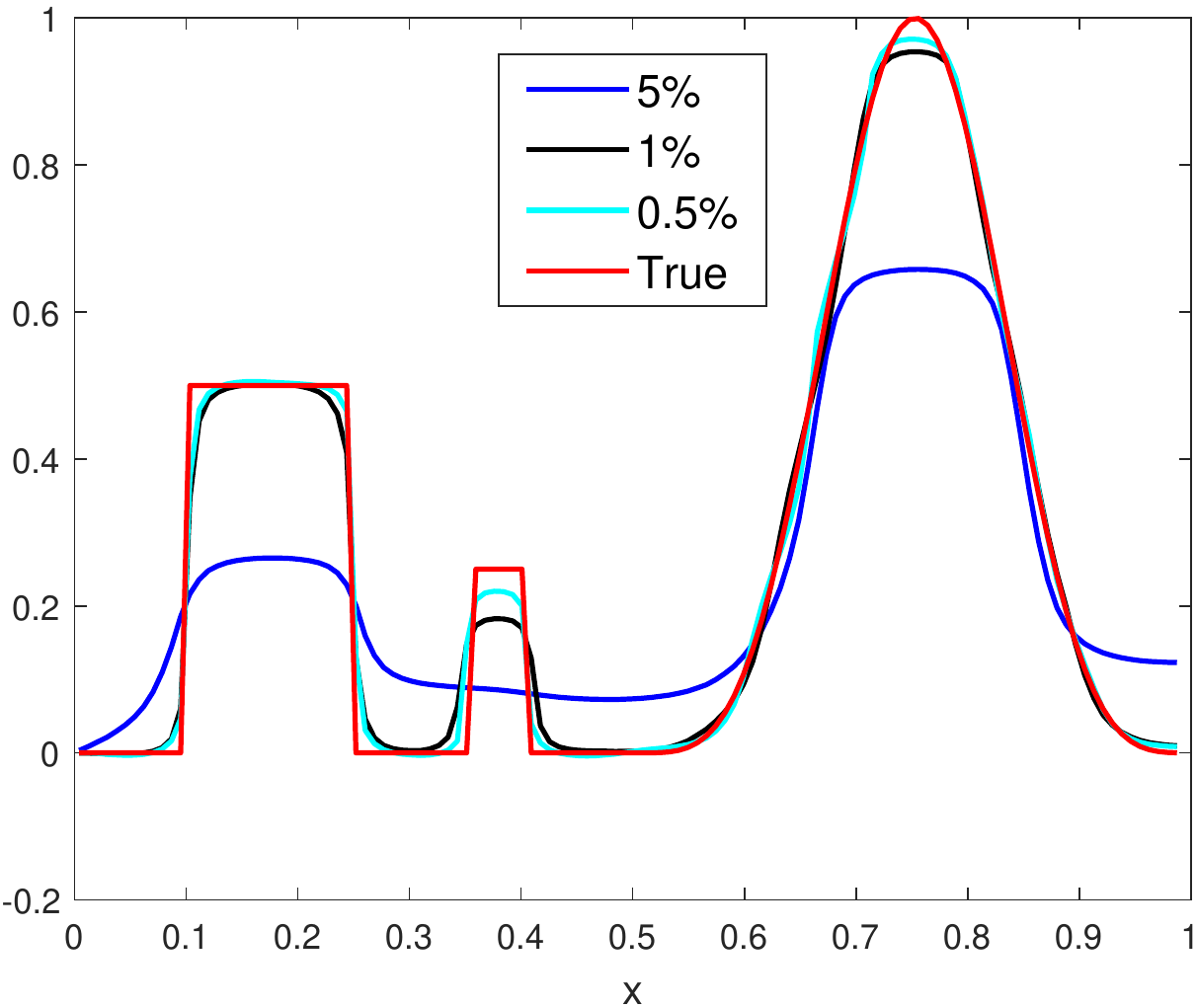} \\
\includegraphics[scale=0.5]{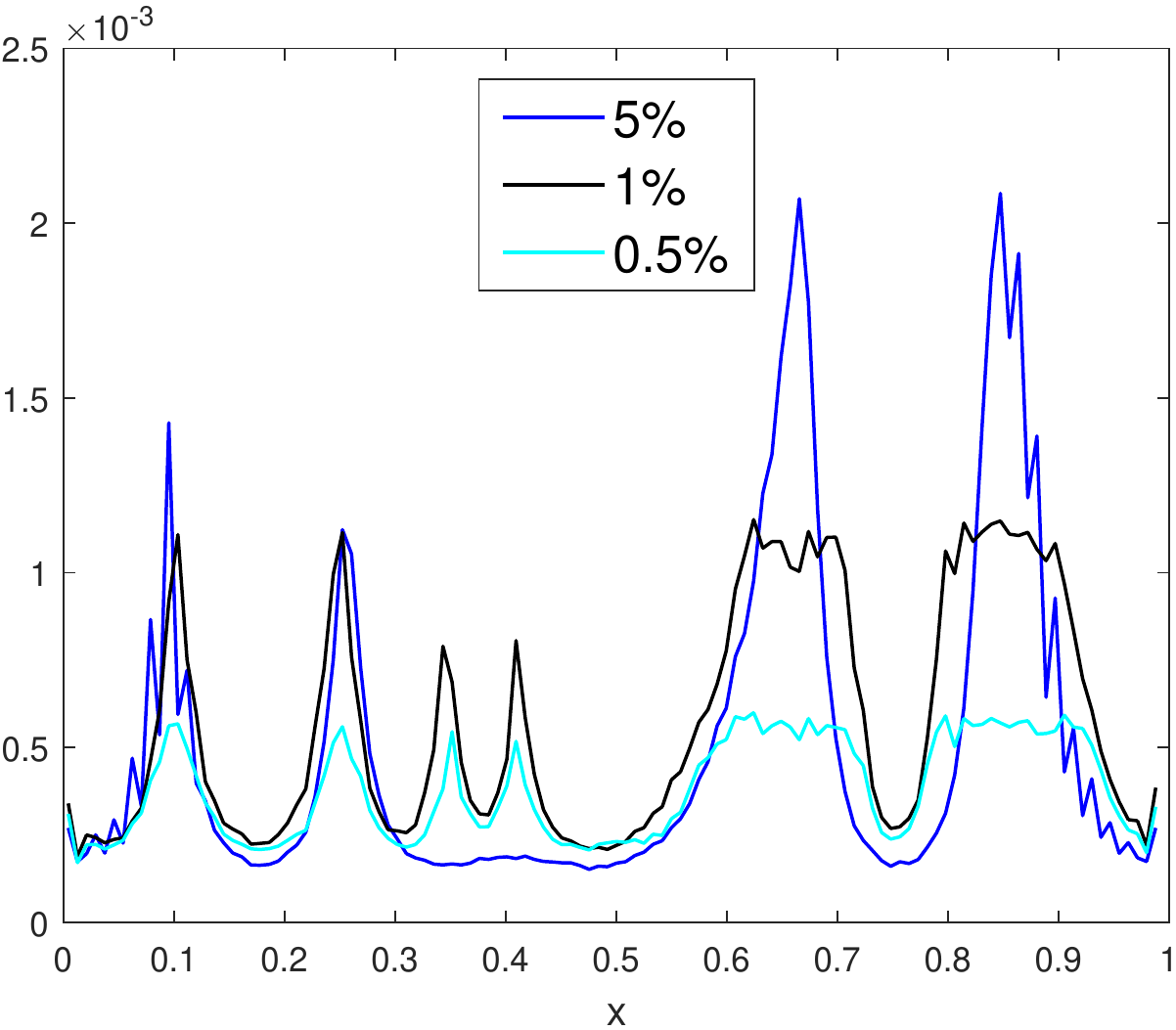}& \hspace{20pt}
\includegraphics[scale=0.5]{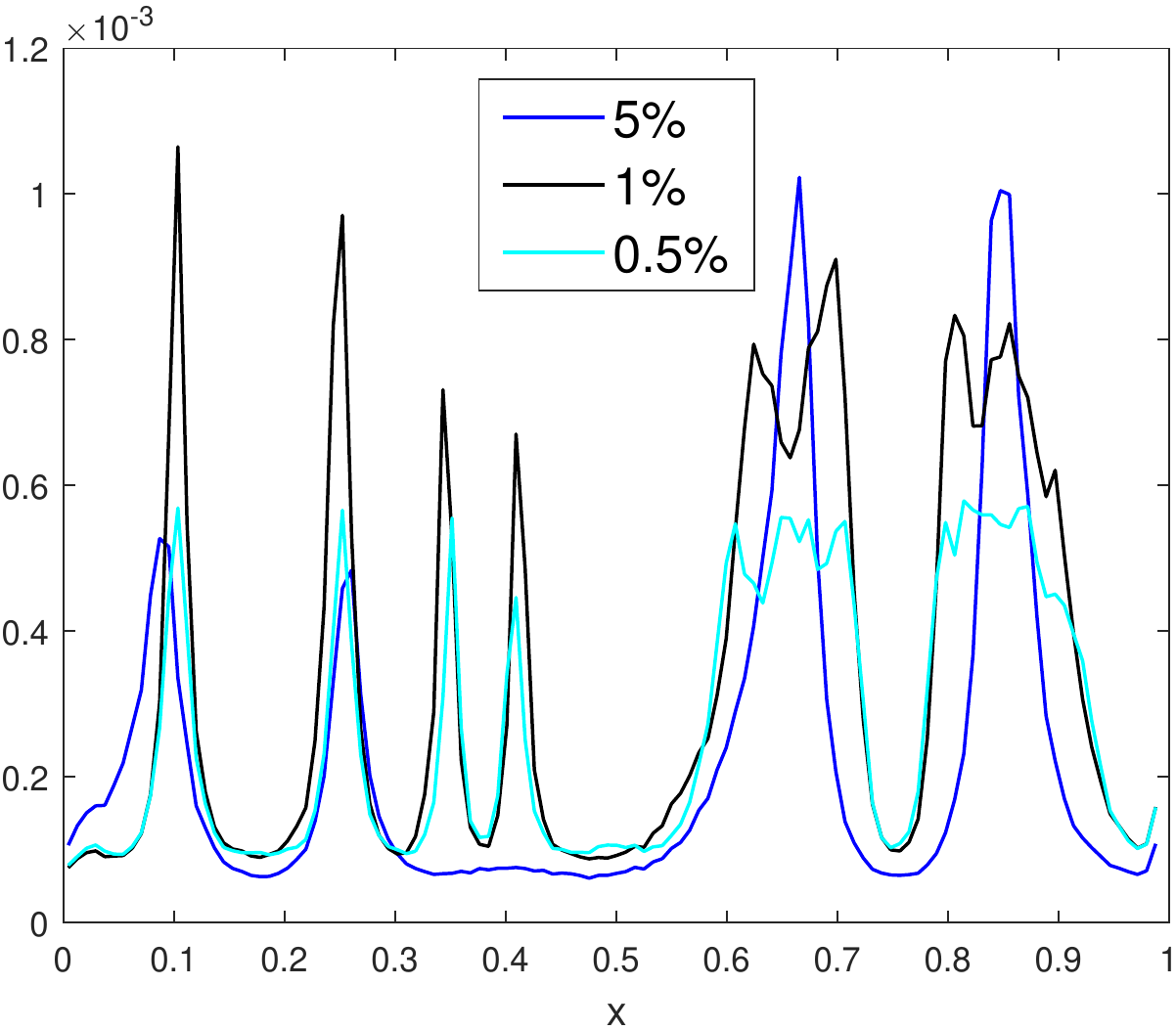} \\
FTG with $\alpha=0.95$ & \hspace{20pt} TG \\
 \end{tabular}
\caption{Deconvolution problem. The posterior mean (top row) and the corresponding posterior standard deviation (bottom row) for different noise level $5\%,\,1\%,\,0.5\%$, using the linear diagonal map-based independence sampler in Algorithm \ref{alg:diagonal_map} for $\alpha=0.95$ FTG and TG prior.
}
\label{decon_result_noise}
\end{figure}

Next, we study the influence of noise level on the numerical results.
The posterior mean and posterior standard deviation for different noise level $5\%,\,1\%,\,0.5\%$, using the linear diagonal map-based independence sampler for FTG prior with $\alpha=0.95$ and TG prior, are plotted in Figure \ref{decon_result_noise}.
Comparing the numerical results, it is observed that the posterior mean is able to approximate the exact $f(x)$ well as the noise level decreases and at the same time, the posterior standard derivation gradually becomes smaller regardless of the TG or FTG prior.


Finally, we study the sensitivity of the inversion results with respect to the parameter $(k,\theta)$ of Gamma distribution and list the corresponding numerical results in Table \ref{decon_table_sen}.
We find that when fixing one of the parameter, changing the other parameter had little effect on the results.
And according to the relationship \eqref{reg_para_MAP}, the parameter $(k,\theta)$ are depended on the value of the $\|u\|_{TV^{\alpha}}$ in some extent.
\begin{table}[htbp]
\centering
\caption{The $RelErr$ values of the deconvolution results with respect to the parameters $(k,\theta)$ of Gamma distribution.}
\begin{tabular}{cccccc}
  \toprule
   & \multicolumn{5}{c}{$(k,\theta)$} \\
  \cmidrule{2-6}
           & $(3000,1)$ & $(2500,1)$ & $(2000,1)$ &$(2000,0.1)$ &$(2000,0.01)$ \\
  \midrule
  $\alpha=0.95$ & 0.0870 & 0.0846 & \textbf{0.0822}&0.0881&0.0891  \\
  TG    & 0.0891 & 0.0865& \textbf{0.0836}&0.0893&0.0902 \\
  $\alpha=1.05$ & 0.1373 & 0.1341 & \textbf{0.1309}&0.1507&0.1574 \\
  \bottomrule
\end{tabular}
\label{decon_table_sen}
\end{table}

\subsection{Inverse source identification problems}\label{1d_ex2}

\subsubsection{Problem setup}
We consider the following initial-boundary value problem for the non-homogeneous heat equation.

\begin{align}
\label{heat_eq1}
\frac{\partial v(\mathbf{x},t)}{\partial t} &=\Delta v(\mathbf{x},t)+f(\mathbf{x}), \quad (\mathbf{x},t)\in \Omega \times (0,\mathsf{T}]\\
\label{heat_eq2}
v(\mathbf{x},0) &=\psi(\mathbf{x}), \quad  (\mathbf{x},t)\in \Omega \\
v(\mathbf{x},t) &=0, \quad (\mathbf{x},t)\in \partial \Omega \times  (0,\mathsf{T}],
\label{heat_eq3}
\end{align}
where the $\Delta$ is the Laplace operator, $\Omega$ is a bounded domain of $\mathbb{R}^m$, $\mathsf{T}>0$, and $\psi$ and the source term $f$ are given functions.

Our task is to determine the heat source on the usual initial-boundary conditions with the assistance of the final temperature data.
In fact, this is inversely determined and usually ill-posed.
In this paper the inverse problem of determining $f(\mathbf{x})$ from the knowledge of $\psi(\mathbf{x})$ and the final temperature measurement
\begin{align}
\label{Tc}
v(\mathbf{x},\mathsf{T})=\phi(\mathbf{x}), \quad \mathbf{x}\in \Omega.
\end{align}
As shown in \cite{JOHANSSON200766, Johansson_2007}, when the data $\phi$ and $\psi$ satisfy suitable conditions, the above linear inverse problem has a unique solution.
Analogy to \cite{yan2010computational}, we use the finite difference method (FDM) to solve the one-dimensional heat equation \eqref{heat_eq1}$-$\eqref{heat_eq3}, and the equation \eqref{heat_eq1} is discretized by using Crank-Nicolson method, i.e. using the forward-difference for the time derivative and a weighted combination of backward-difference and forward-difference approximations for the remainder of the equation.
Let $\Omega$ be the interval $(0,r)$, and $\Delta t$ is the time step size and use the notation, $v_n(x,t_n)$ where $t_n=n\Delta t,\; n=1,2,\dots,N,\;t_0=0$ and $t_N=\mathsf{T}$,  then we can get the time-discrete equation
\begin{align}
\Bigl(\frac{1}{\Delta t}-w\Delta\Bigr)v_{n+1}=\Bigl(\frac{1}{\Delta t}+(1-w)\Delta\Bigr)v_{n}+f,
\end{align}
where $0\le w\le 1$ is the weights.
We use the uniform grid and let $\Delta x$ be the space step size,  $x_j=j\Delta x, \; j=0,1,2,\dots,d+1$ be a set of discrete points that defines the grid.
The operator $\Delta$ is discretized by using the second-order central difference and apply the initial conditions, and then we can obtain the fully discrete equation
\begin{align}
\label{de_eq}
D_+V_{n+1}=D_-V_{n}+\mathbf{f},\quad n=0,1,2,\dots,N-1,
\end{align}
where $V_n=\bigl(v_n(x_1),v_n(x_2),\dots,v_n(x_d)\bigr)^T$, $\mathbf{f}=\bigl(f(x_1),f(x_2),\dots,f(x_d)\bigr)^T$, and $D_+,\;D_-$ are the discretization matrices of the $1/\Delta t-w\Delta,\; 1/\Delta t-(1-w)\Delta$, respectively.

Both sides of the equation \eqref{de_eq} multiply the inverse matrix of $D_+$ and recurse, which can yield
\begin{align}
V_N=D^NV_0+\sum_{i=0}^{N-1}D^iD_+^{-1}\mathbf{f},
\end{align}
where $D=D_+^{-1}D_-$, $D^0=I$, the condition \eqref{Tc}  $V_N=(\phi(x_1),\phi(x_2),\dots,\phi(x_d))^T$ and the initial condition \eqref{heat_eq2} $V_0=(\psi(x_1),\psi(x_2),\dots,\psi(x_d))^T$.
Let the $\mathbf{b}=V_N-D^NV_0,\; H=\sum_{i=0}^{N-1}D^iD_+^{-1}$ and then we can get the matrix equation $\mathbf{b}=H\mathbf{f}$.
Thus, the solution of heat equation can be obtained by solve the matrix equation and at the same time  we can get the observations through $\mathbf{y}=H\mathbf{f}+D^NV_0+\boldsymbol{\eta}$ where the $D^NV_0+H\cdot$ is the discretization of our linear forward model $\mathcal{A}$ and the $\boldsymbol{\eta}$ is the noise.

\subsubsection{Set up of inverse problems}
In order to illustrate the numerical results, we take $\mathsf{T}=1,\,r=12$, $w=0.5$ and the initial temperature
\[
\psi(x)=sin(\pi x), \quad x\in [0,12].
\]
In this example, we consider the heat source:
\[ f(x)=
\begin{cases}
0.5 &  0.75\le x<2,\\
-(x-3)(x-5) &  3\le x<5,\\
x-5 & 5\le x<6,\\
-x+7 &  6\le x<7,\\
-(x-7)(x-9) &  7\le x<9,\\
0.5 &  10\le x<11.25,\\
0 &  otherwise,
\end{cases}
\]
and set $d=150,\, N=120$.
The measurement data are obtained by $\mathbf{y}=H\mathbf{f}+D^NV_0+\boldsymbol{\eta}$ and $\boldsymbol{\eta}$ is assumed to be the Gaussian noise with zero mean and standard deviation $6.185\times10^{-3},6.185\times10^{-4},6.185\times10^{-5}$, which corresponds to $1\%,0.1\%,0.01\%$ noise, respectively, with respect to the maximum norm of the output $H\mathbf{f}+D^NV_0$.
Note that the measurement data $\mathbf{y}$ is computed from a twice finer grid.
The Gaussian prior $\mu_0$ is taken by zero mean and set $\gamma=0.03, \nu=0.0009$ and the shape and rate parameter about Gamma distribution are set to $k=1\times10^4,\;\vartheta=1$, respectively.

Then using the alternating direction method, we construct a linear diagonal transport map between the Gaussian distribution $\mu_0$ and the posterior \eqref{H_postprior}.
Analogy to above deconvolution problem, we use $M=1000$ samples of the Gaussian distribution $\mu_0$ to approximate the expected value via the SAA in the objective function (see Section \ref{numercal opt}).
And then the regularization parameter $\lambda$ can be determined by the formula \eqref{reg_para_MAP}, which together with the diagonal map $\widetilde{T}$ is used as the precondition for the independence sampler as shown in Algorithm \ref{alg:diagonal_map}.
The proposal $\mu_{ref}$ is the Gaussian distribution with zero mean and the standard deviation $4\times10^{-3}$ in the sampling process.
\subsubsection{Result}

\begin{figure}[htbp]
 \centering
 \begin{tabular}{@{}c@{}c@{}c@{}}
\includegraphics[scale=0.5]{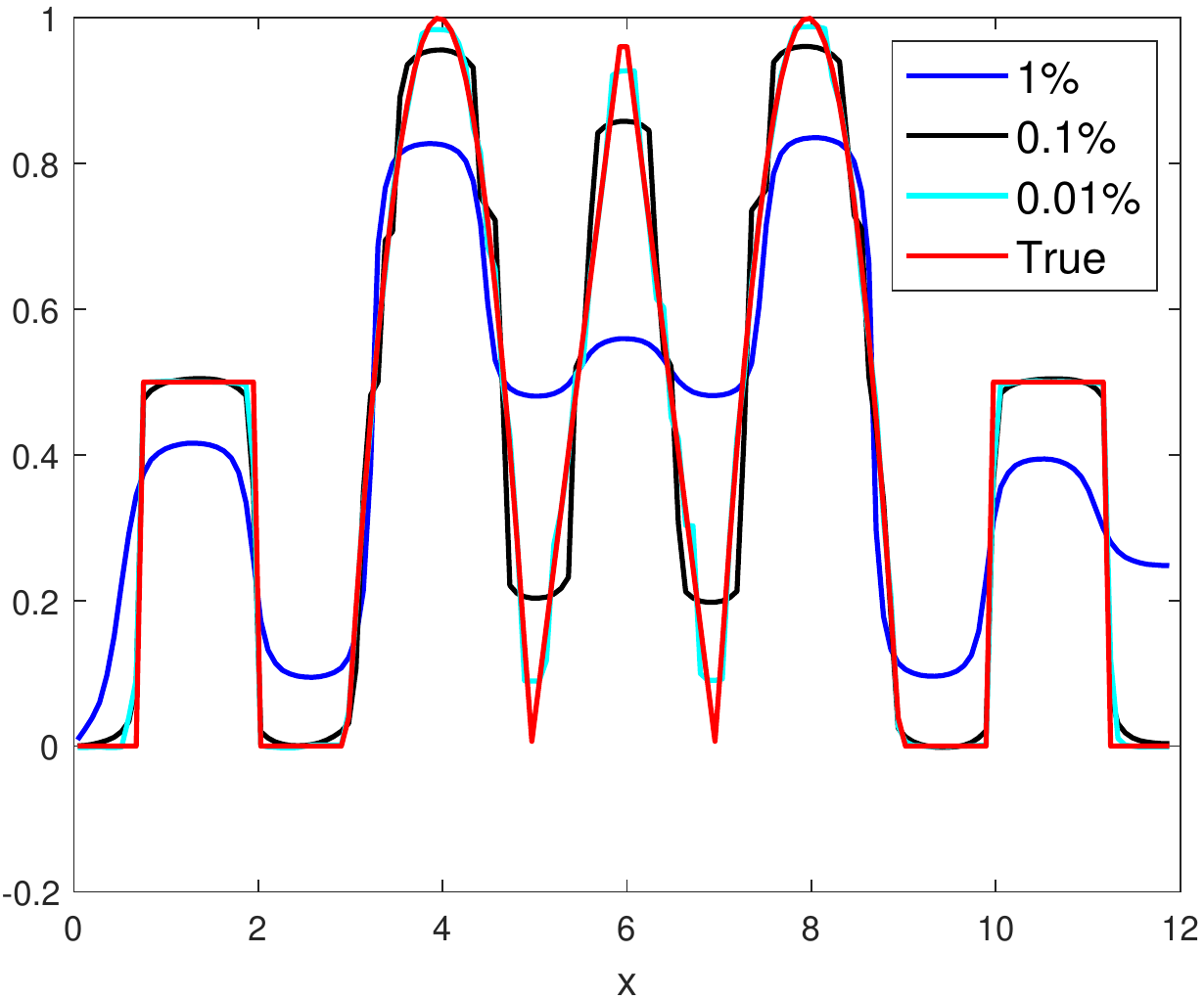}& \hspace{20pt}
\includegraphics[scale=0.5]{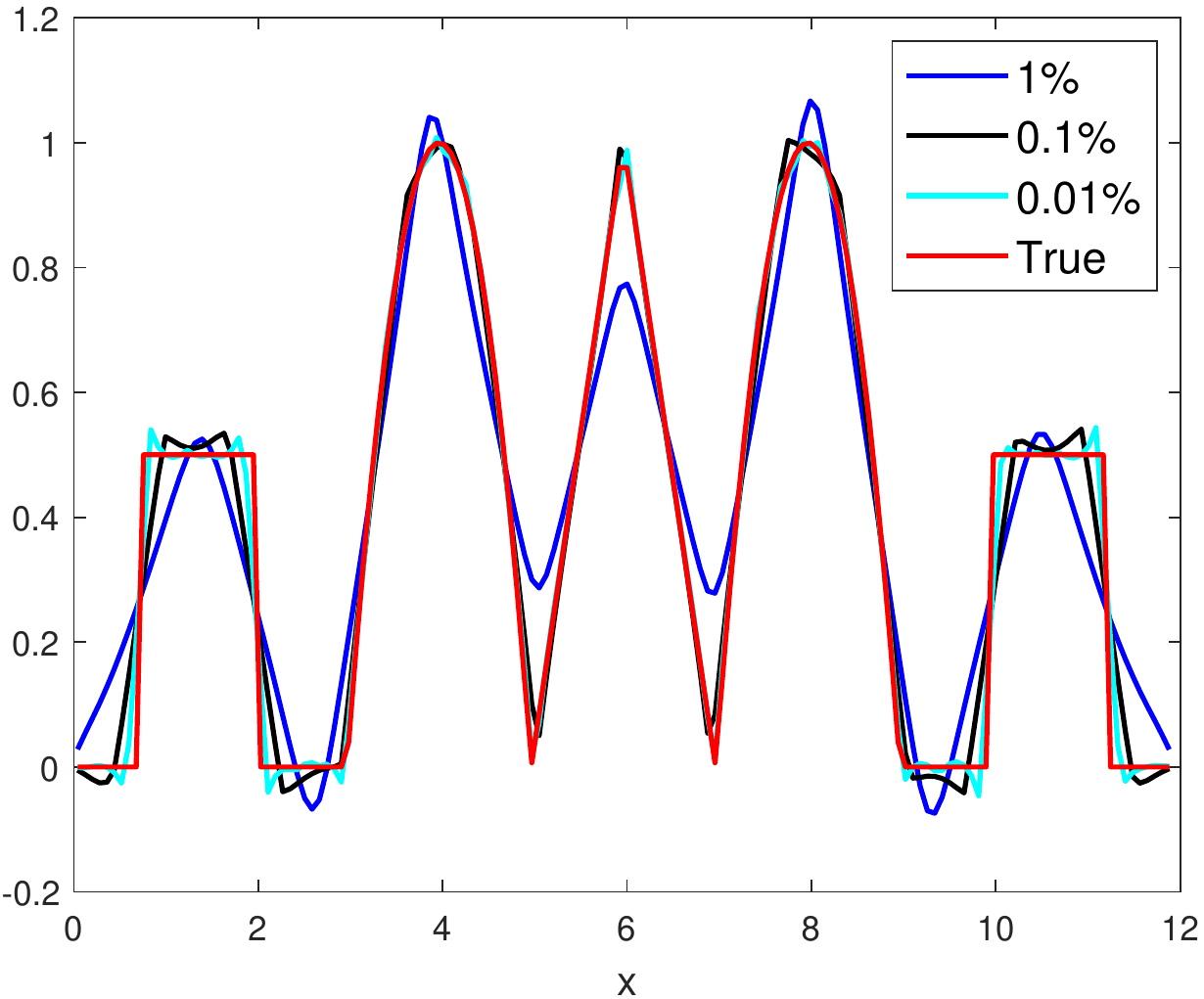}\\
\includegraphics[scale=0.5]{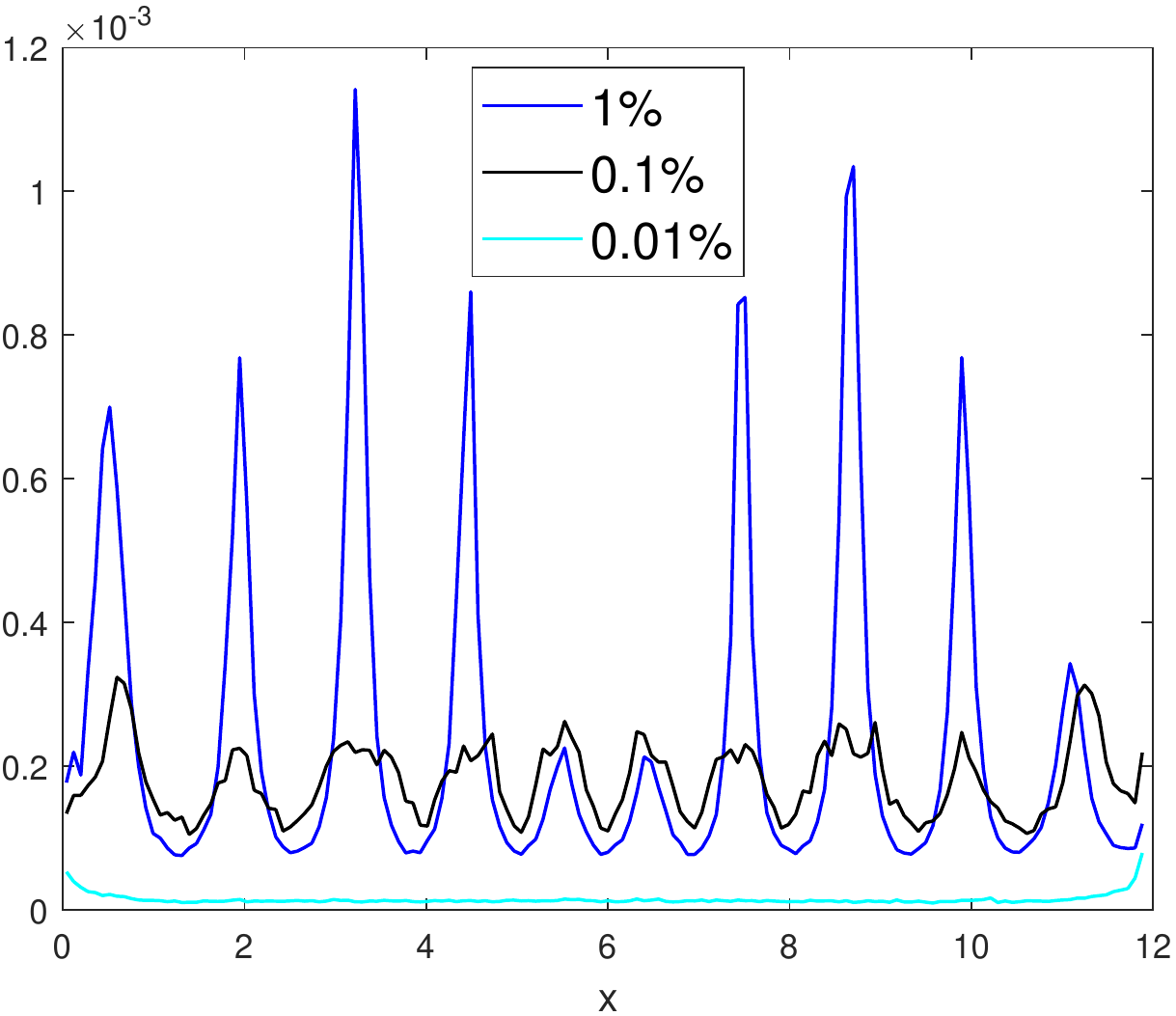}& \hspace{20pt}
\includegraphics[scale=0.5]{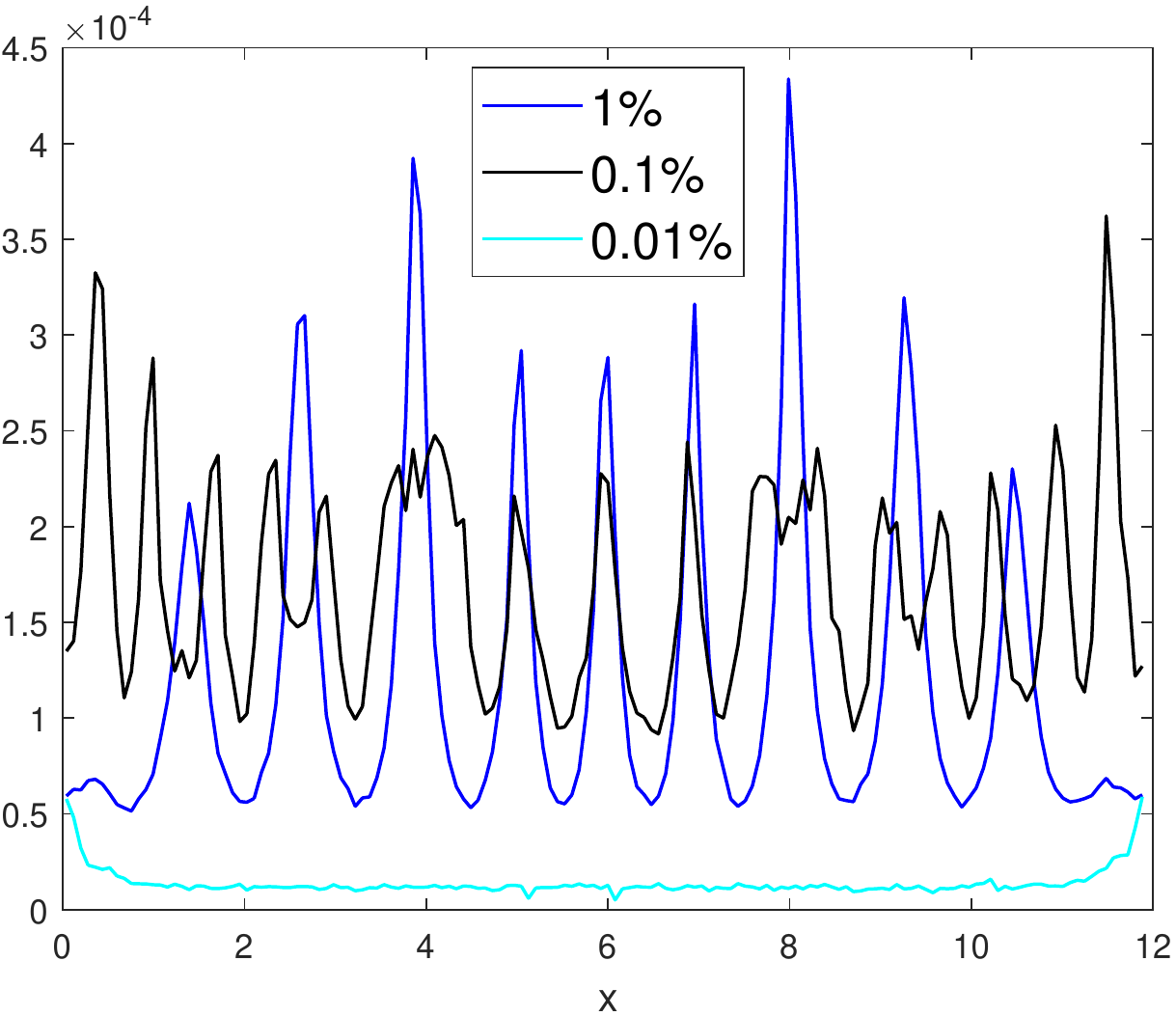}\\
TG&FTG with $\alpha=1.1$ \\
 \end{tabular}
\caption{Inverse source identification problems. The posterior mean (top row) and the corresponding posterior standard deviation (bottom row) for different noise level $1\%,\,0.1\%,\,0.01\%$, using the linear diagonal map-based independence sampler Algorithm \ref{alg:diagonal_map} with TG and $\alpha=1.1$ FTG prior.
}
\label{ex2_result_noise}
\end{figure}

At first, we consider the robustness of the diagonal map-based independence sampler Algorithm \ref{alg:diagonal_map} with respect to the noise level.
The posterior mean and posterior standard deviation for different noise level $1\%,\,0.1\%,\,0.01\%$, using our independence sampler with FTG for $\alpha=1.1$ and TG prior, are plotted in Figure \ref{ex2_result_noise}.
Comparing the numerical results, it is observed that the posterior mean is more and more consistent with the exact $f(x)$ as the noise level becomes smaller and simultaneously, the posterior standard derivation presents a stable trend and becomes smaller gradually whether it is for TG or FTG.
In the following numerical results for this example, we fix the noise level as $0.1\%$.

\begin{figure}[htbp]
 \centering
 \begin{tabular}{@{}c@{}c@{}}
\includegraphics[width=0.45\textwidth, height=0.25\textwidth]{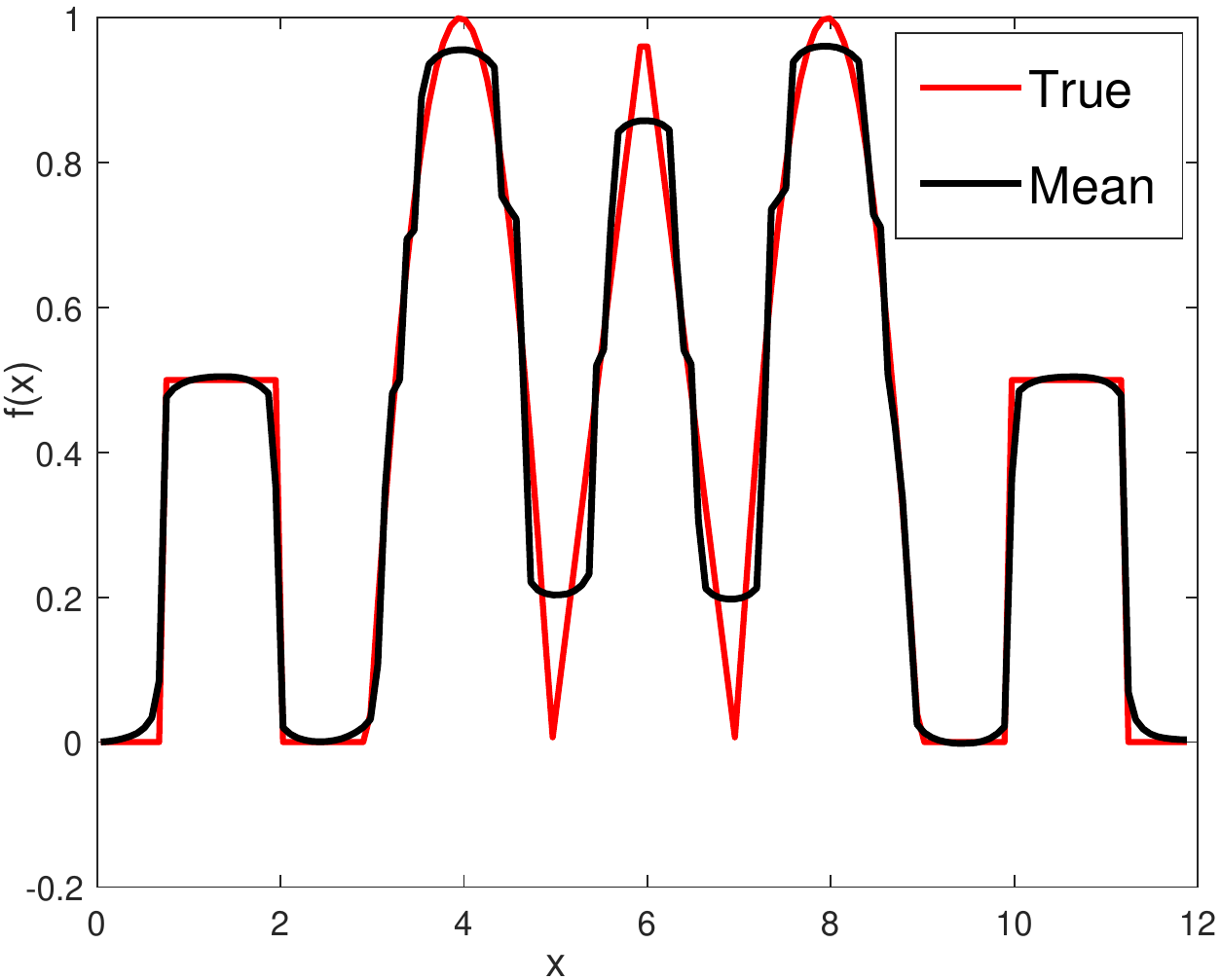}&
\includegraphics[width=0.45\textwidth, height=0.25\textwidth]{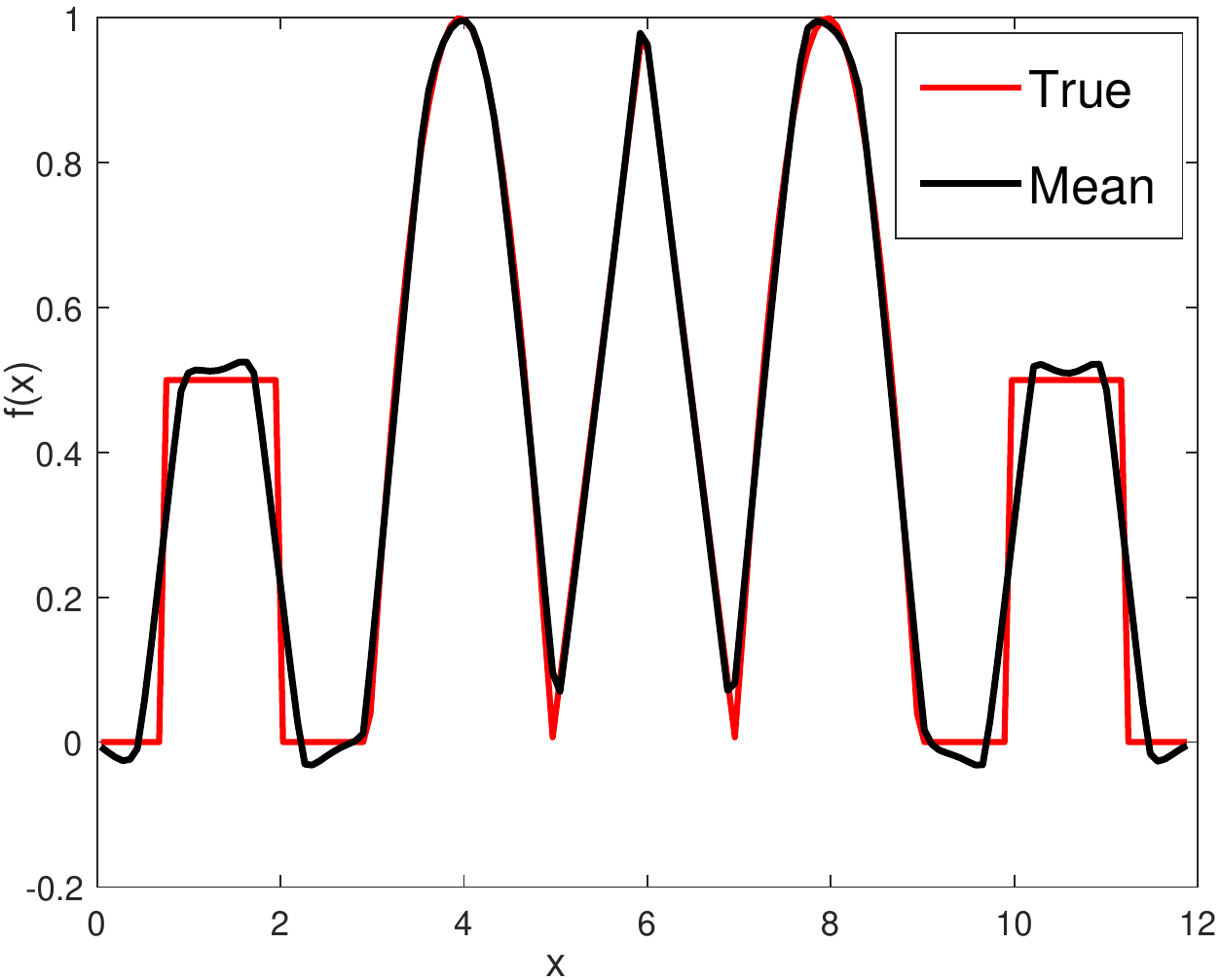} \\
\includegraphics[width=0.45\textwidth, height=0.12\textwidth]{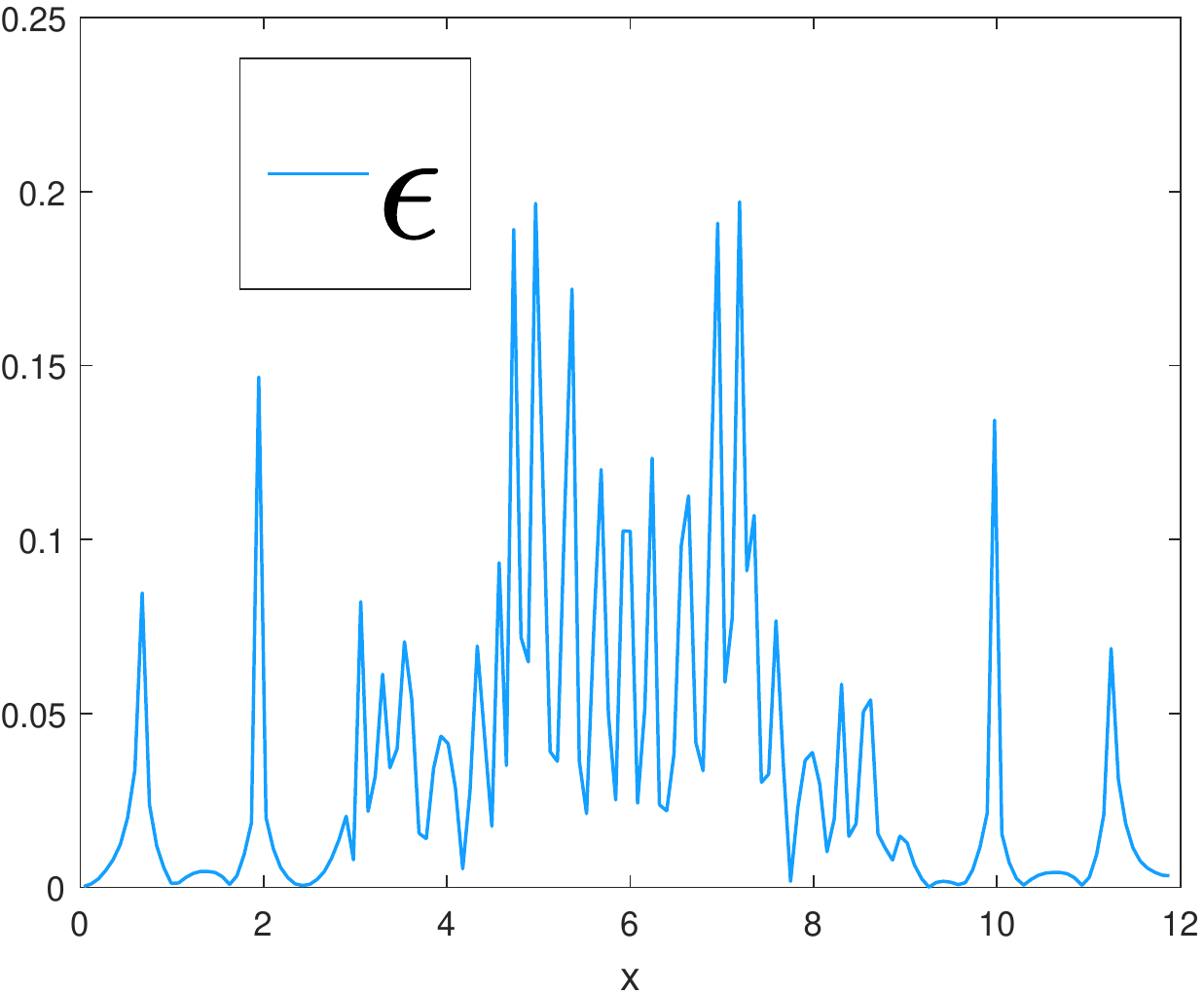}&
\includegraphics[width=0.45\textwidth, height=0.12\textwidth]{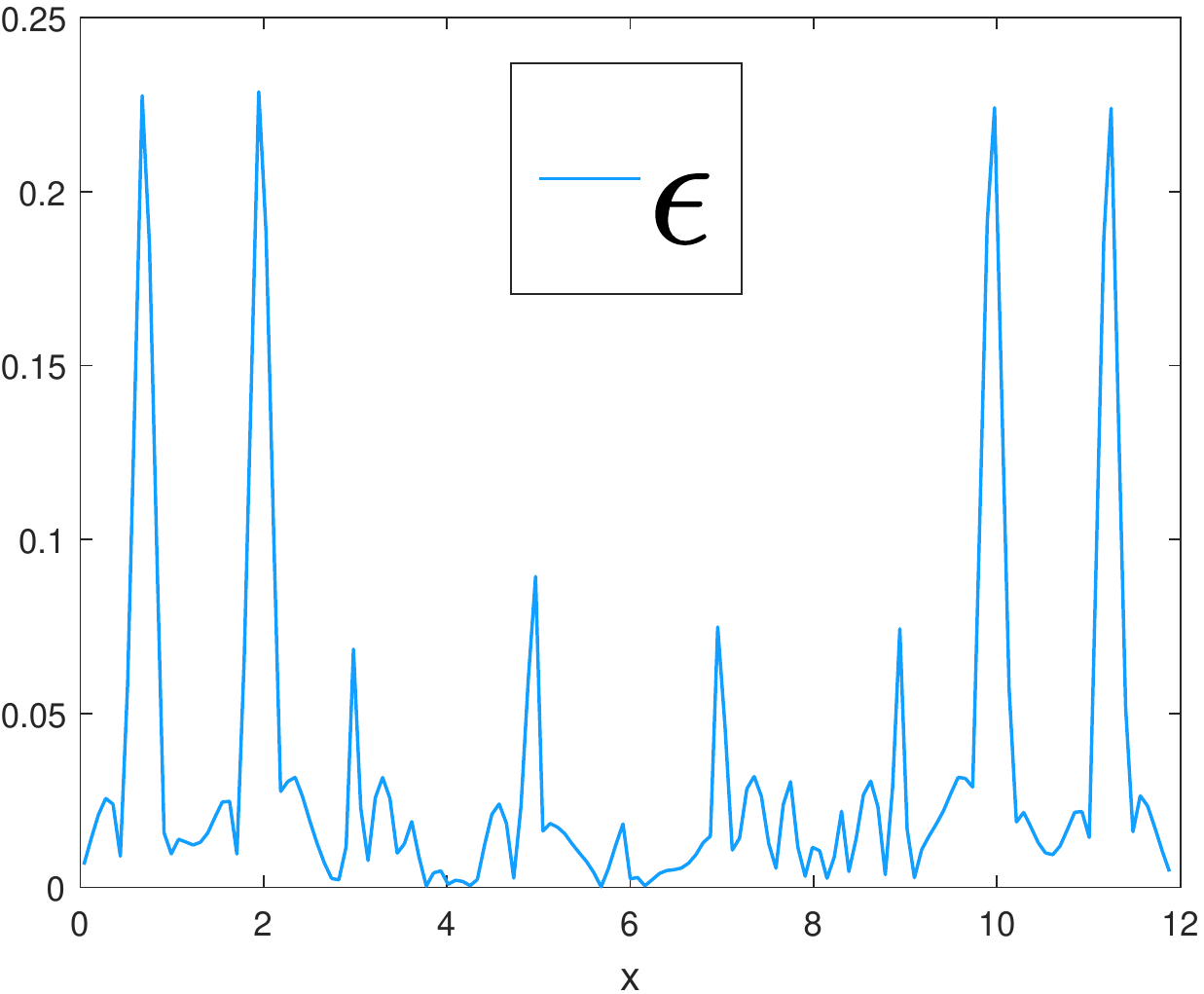}  \vspace{-5pt}\\
TG&FTG with $\alpha=1.05$\\  \hline
\includegraphics[width=0.45\textwidth, height=0.25\textwidth]{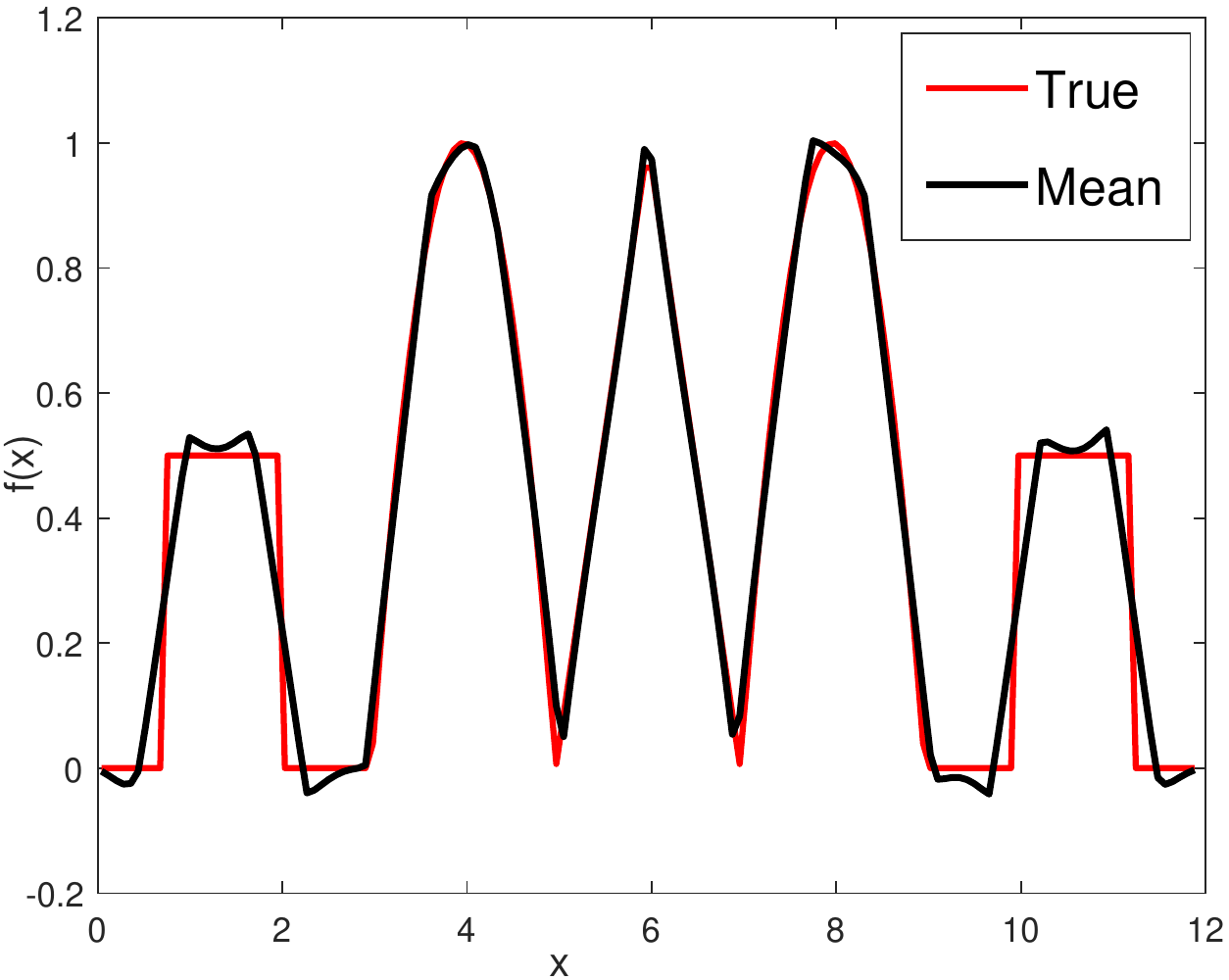} \hspace{2pt}&
\includegraphics[width=0.45\textwidth, height=0.25\textwidth]{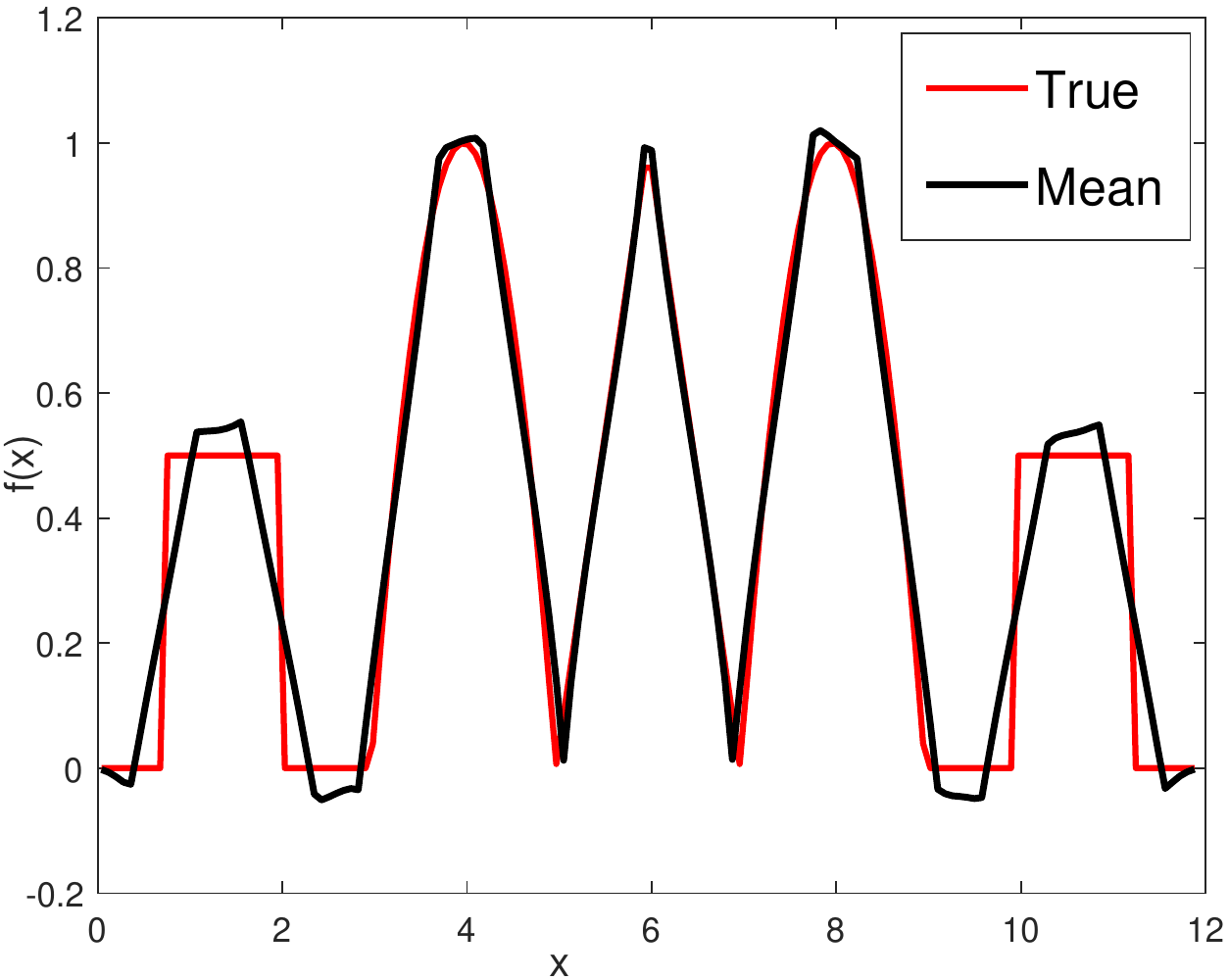}\\
\includegraphics[width=0.45\textwidth, height=0.12\textwidth]{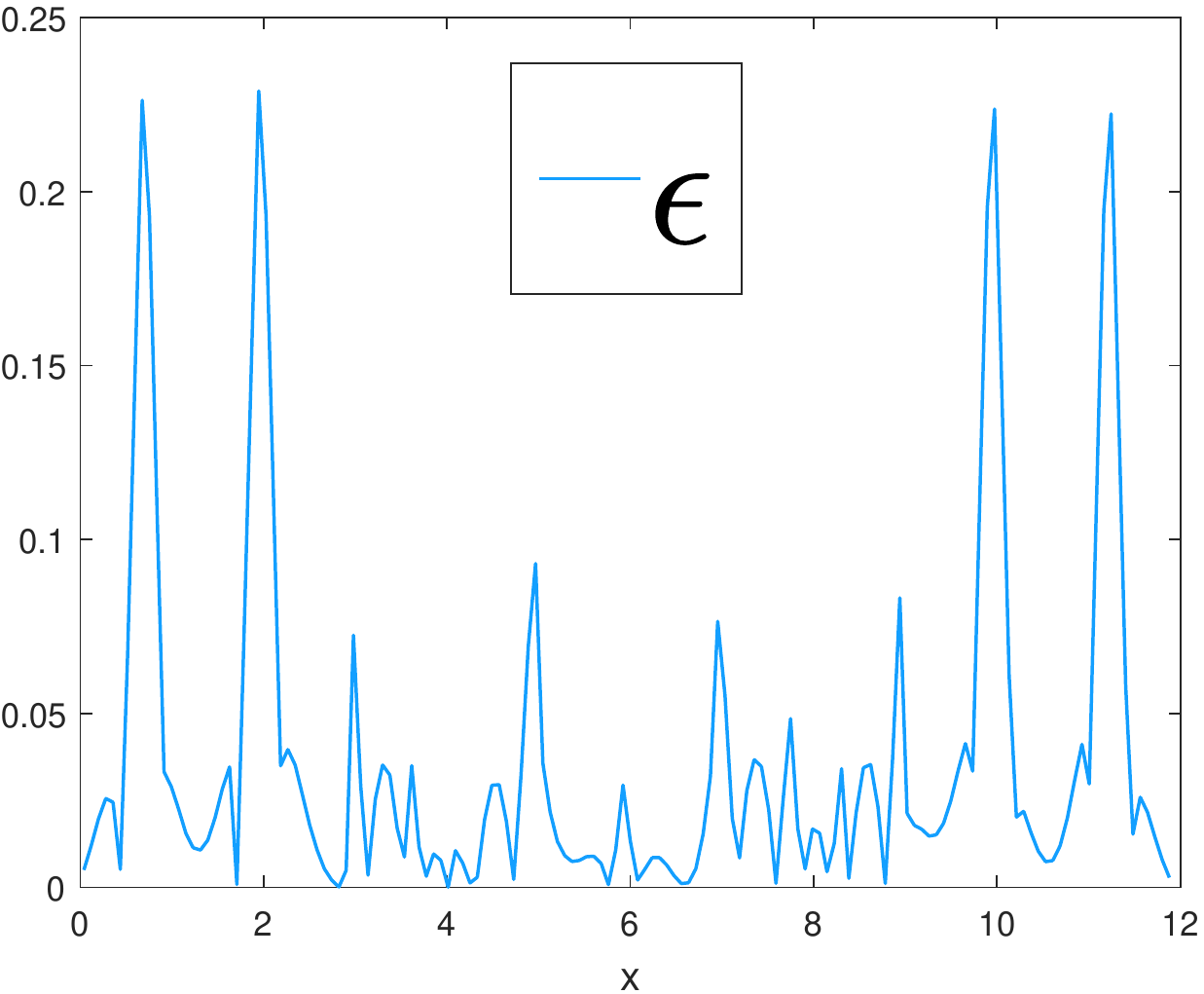}&
\includegraphics[width=0.45\textwidth, height=0.12\textwidth]{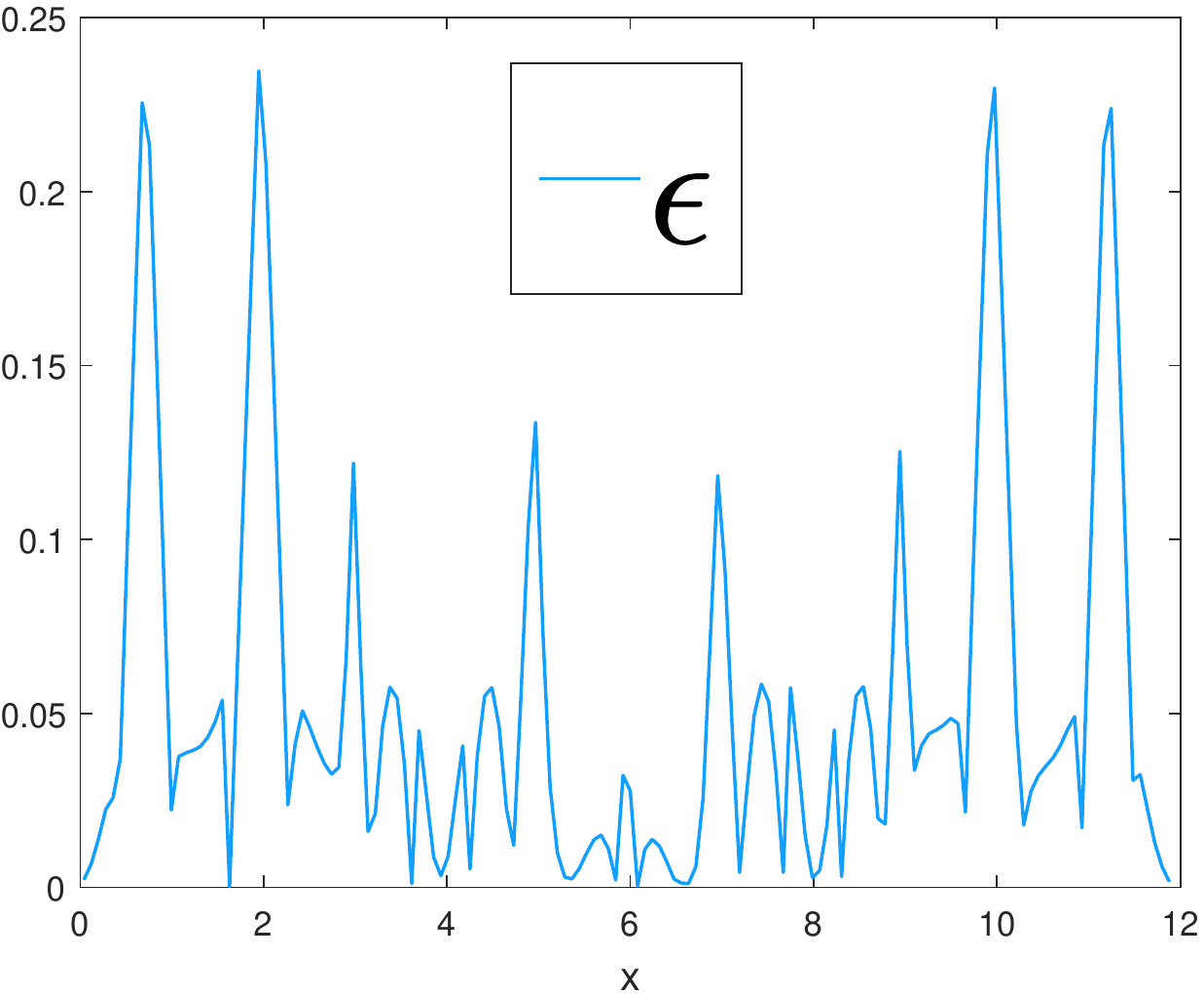}  \vspace{-5pt}\\
 FTG with $\alpha=1.1$ &FTG with $\alpha=1.5$\\ \hline
 \end{tabular}
\caption{Reconstruction results for the heat source $f(\mathbf{x})$.
The posterior mean and absolute error $\epsilon$ using the linear diagonal map-based independence sampler Algorithm \ref{alg:diagonal_map} for FTG with $\alpha=1.05,\;1.1,\;1.5$ and TG prior.
}
\label{ex2_result_12}
\end{figure}

Similar to the above deconvolution example, below we compare the reconstruction performance in two cases.
One situation is for the FTG prior with $1<\alpha\le2$.
The posterior mean and absolute error $\epsilon$ using FTG with $\alpha=1.05,\;1.1,\;1.5$ and TG prior are plotted in Figure \ref{ex2_result_12} and the relative error $RelErr$ are listed in Table \ref{table_heat_12}.
As we can see from this tables, the $RelErr$ values using our FTG prior for $1<\alpha\le2$ are consistent with the TG prior and are getting smaller and smaller with the fractional order $\alpha$ approaching to $1^+$.
And the FTG prior with $\alpha=1.05,\;1.1$ has better performs than the TG prior in term of the complex oscillation structure of the heat source from the Figure \ref{ex2_result_12}.
In this figure and the corresponding absolute error curve, our FTG prior with $\alpha=1.05,\;1.1$ can eliminate well the staircase effect whether it is sharp points or a smooth part in the reconstruction target $f(x)$ and at the same time also catch the piecewise constant part in it, compared with the numerical result using the TG prior.
With the fractional order $\alpha \to 1^+$, FTG prior can not only reconstruct the complex oscillation structures of target well, but also maintain a relative error similar to that of TG prior, cf. the Table \ref{table_heat_12} and Figure \ref{ex2_result_12}.
Therefore, obviously, the FTG prior with some appropriate $\alpha$ in $(1,2]$ can obtain better recovery result than that using the TG prior, when the reconstructed target has much complex oscillation information, as discussed the image denoising example in Section \ref{2d_ex2}.

\begin{table}[htbp]
\centering
\caption{The $RelErr$ values of the inverse source results using FTG with $1<\alpha\le2$ and TG prior.}
\begin{tabular}{cccccccccc}
  \toprule
          & TG & $\alpha=1.01$ & $\alpha=1.05$ & $\alpha=1.1$ & $\alpha=1.2$ & $\alpha=1.5$ & $\alpha=1.8$ \\
  \midrule
  $RelErr$ & 0.1054&0.1096& 0.1124& 0.1173 &0.1266& 0.1446 & 0.1571\\
    \bottomrule
  \end{tabular}
  \label{table_heat_12}
\end{table}

\begin{figure}[htbp]
 \centering
 \begin{tabular}{@{}c@{}c@{}}
\includegraphics[width=0.45\textwidth, height=0.25\textwidth]{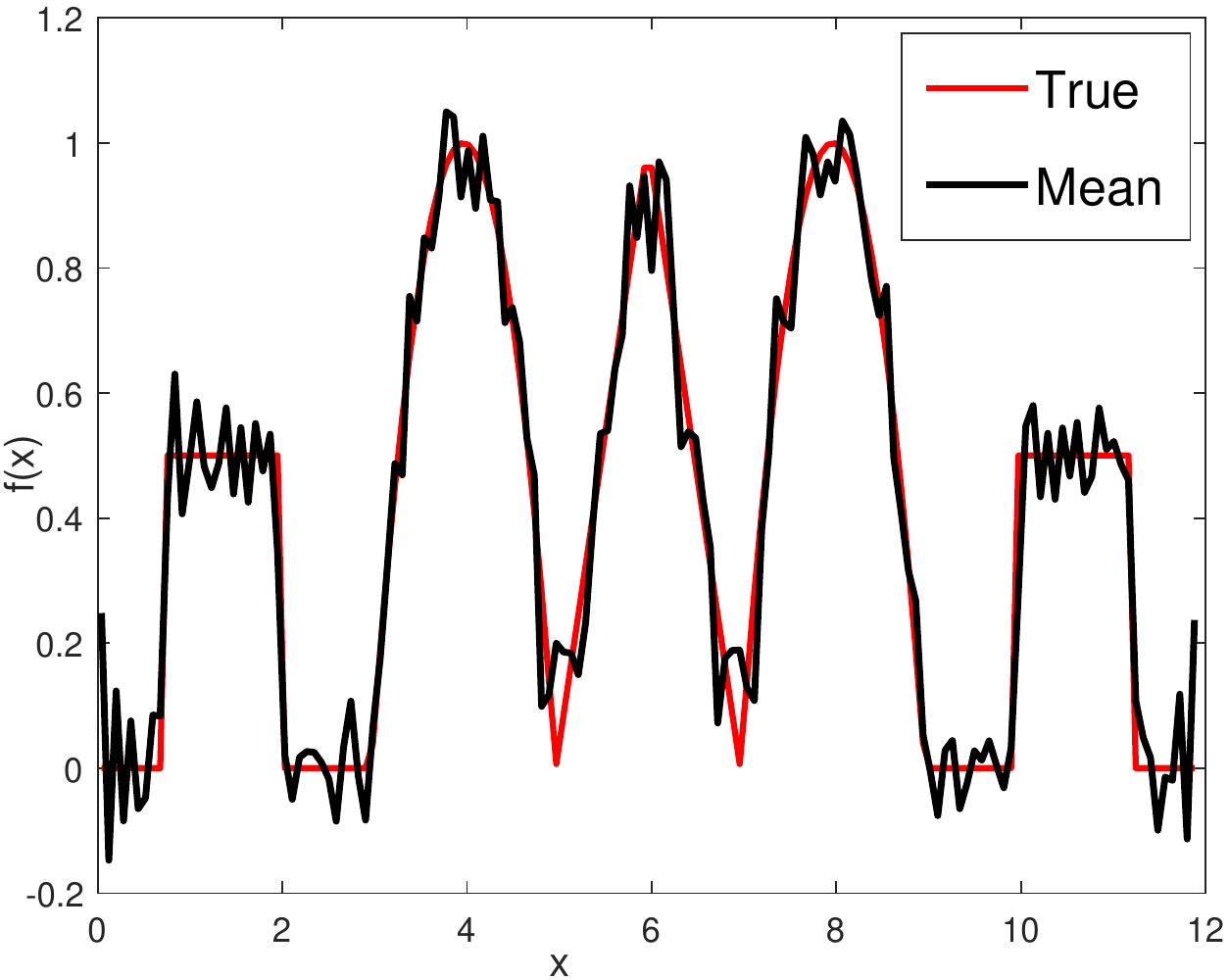}&
\includegraphics[width=0.45\textwidth, height=0.25\textwidth]{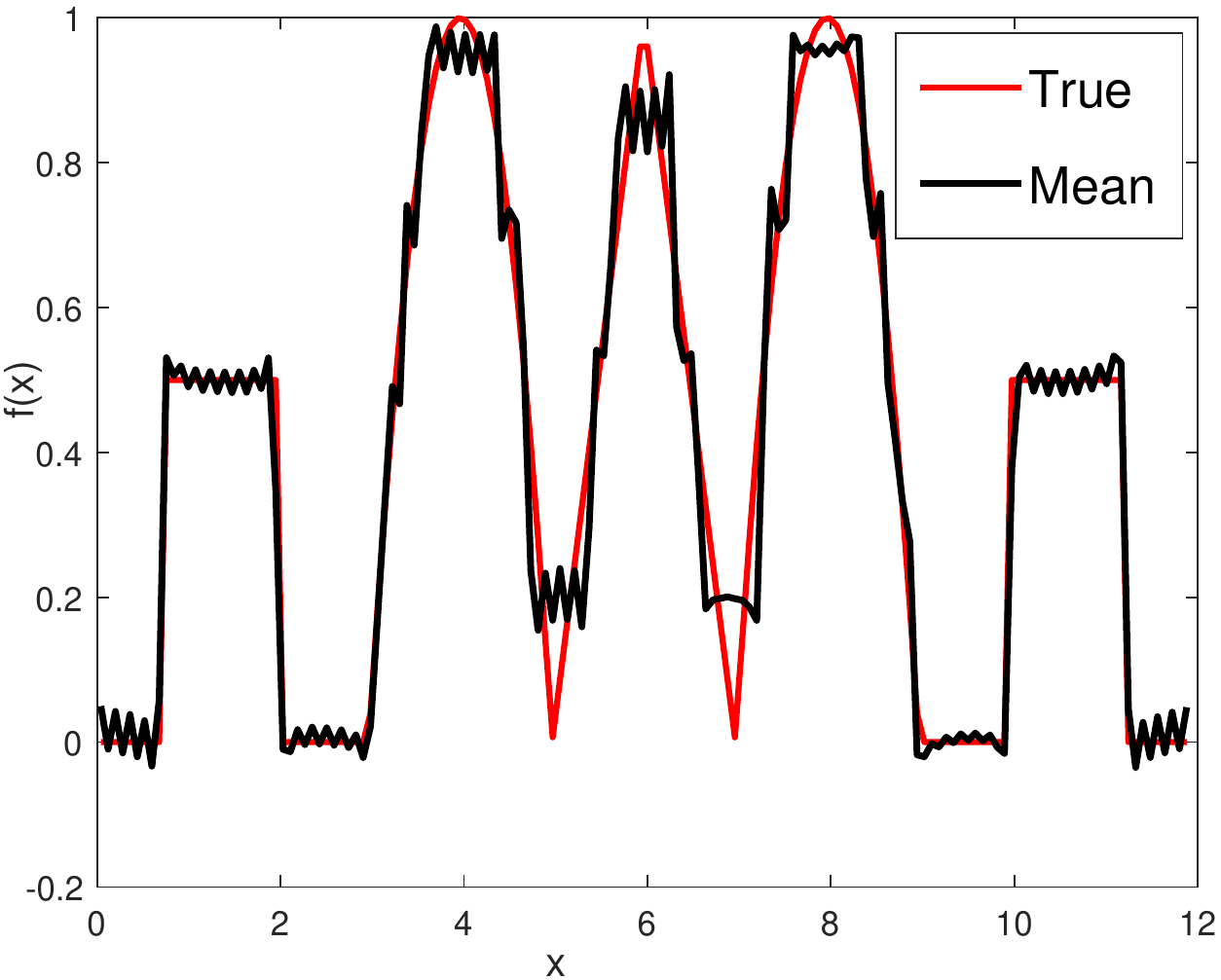} \\
\includegraphics[width=0.45\textwidth, height=0.12\textwidth]{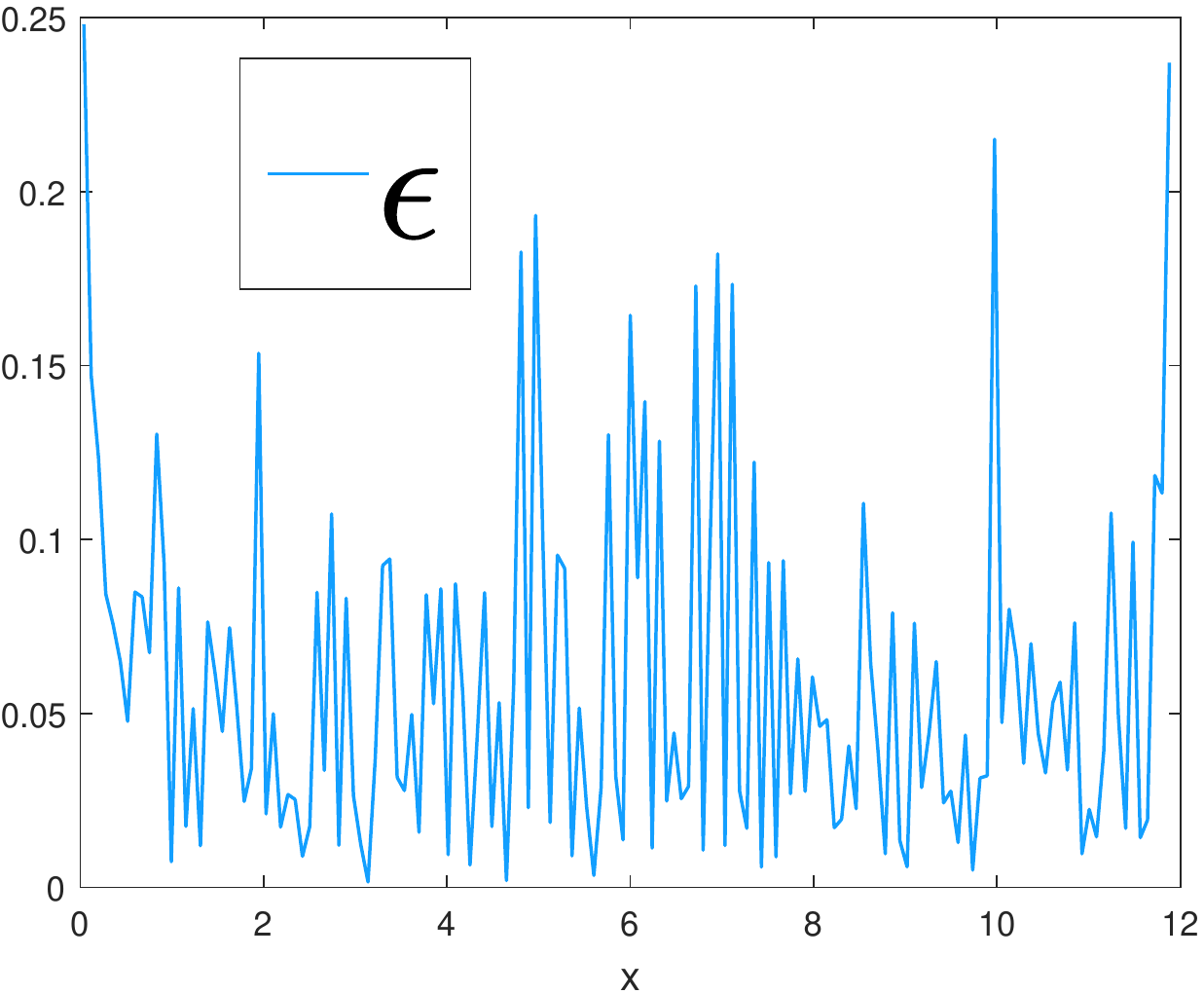}&
\includegraphics[width=0.45\textwidth, height=0.12\textwidth]{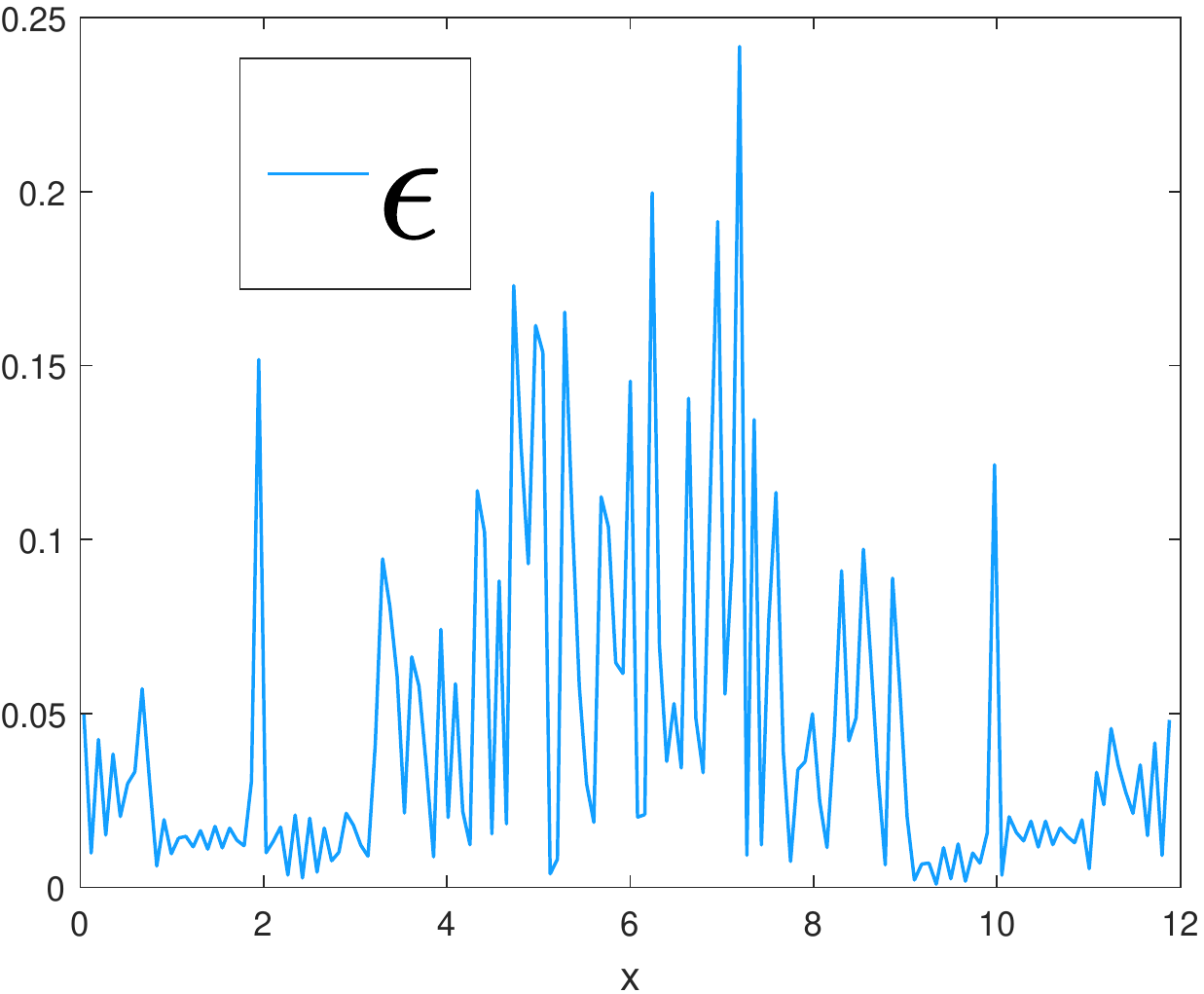}  \vspace{-5pt}\\
FTG with $\alpha=0.5$&FTG with $\alpha=0.9$\\  \hline
\includegraphics[width=0.45\textwidth, height=0.25\textwidth]{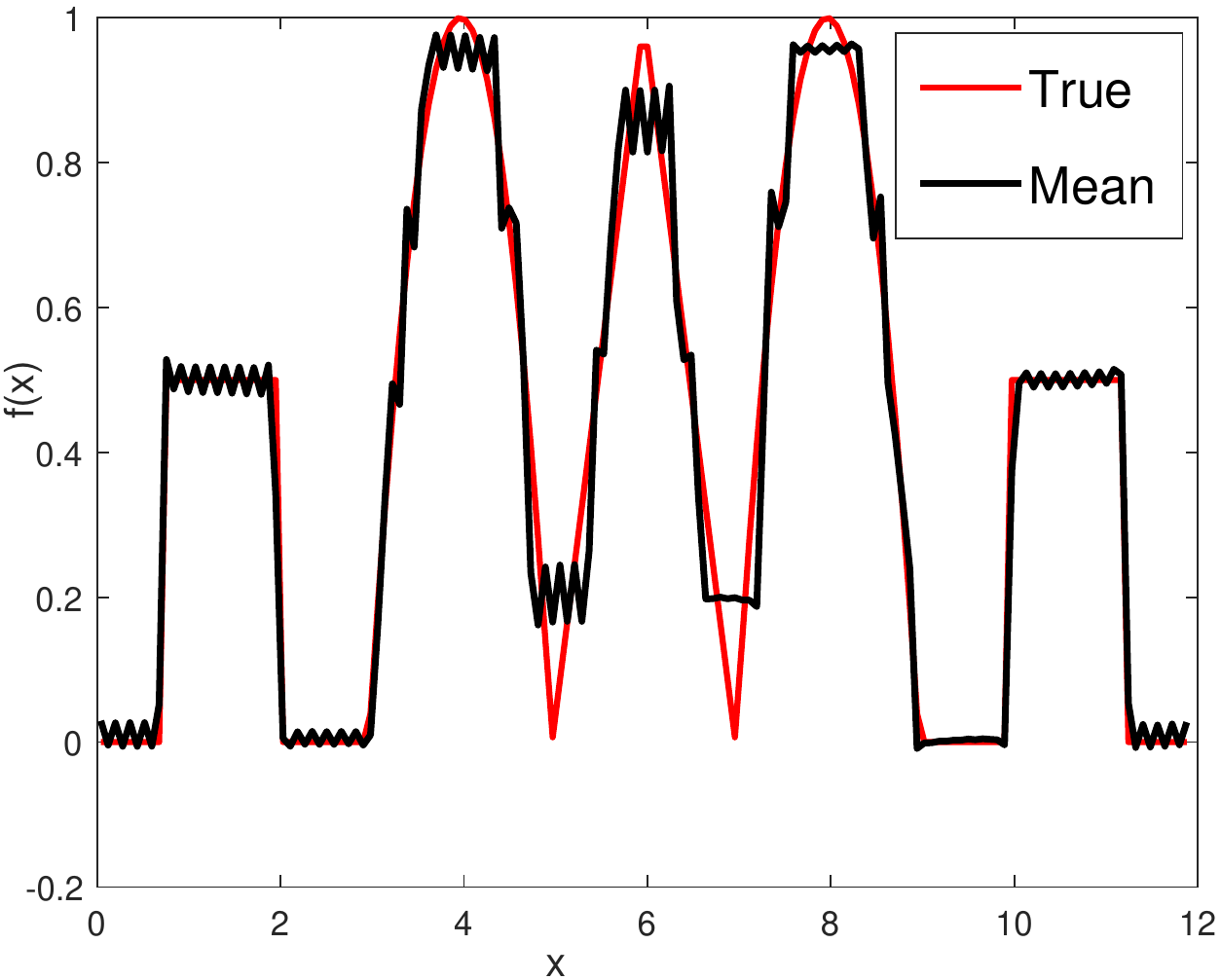} \hspace{2pt}&
\includegraphics[width=0.45\textwidth, height=0.25\textwidth]{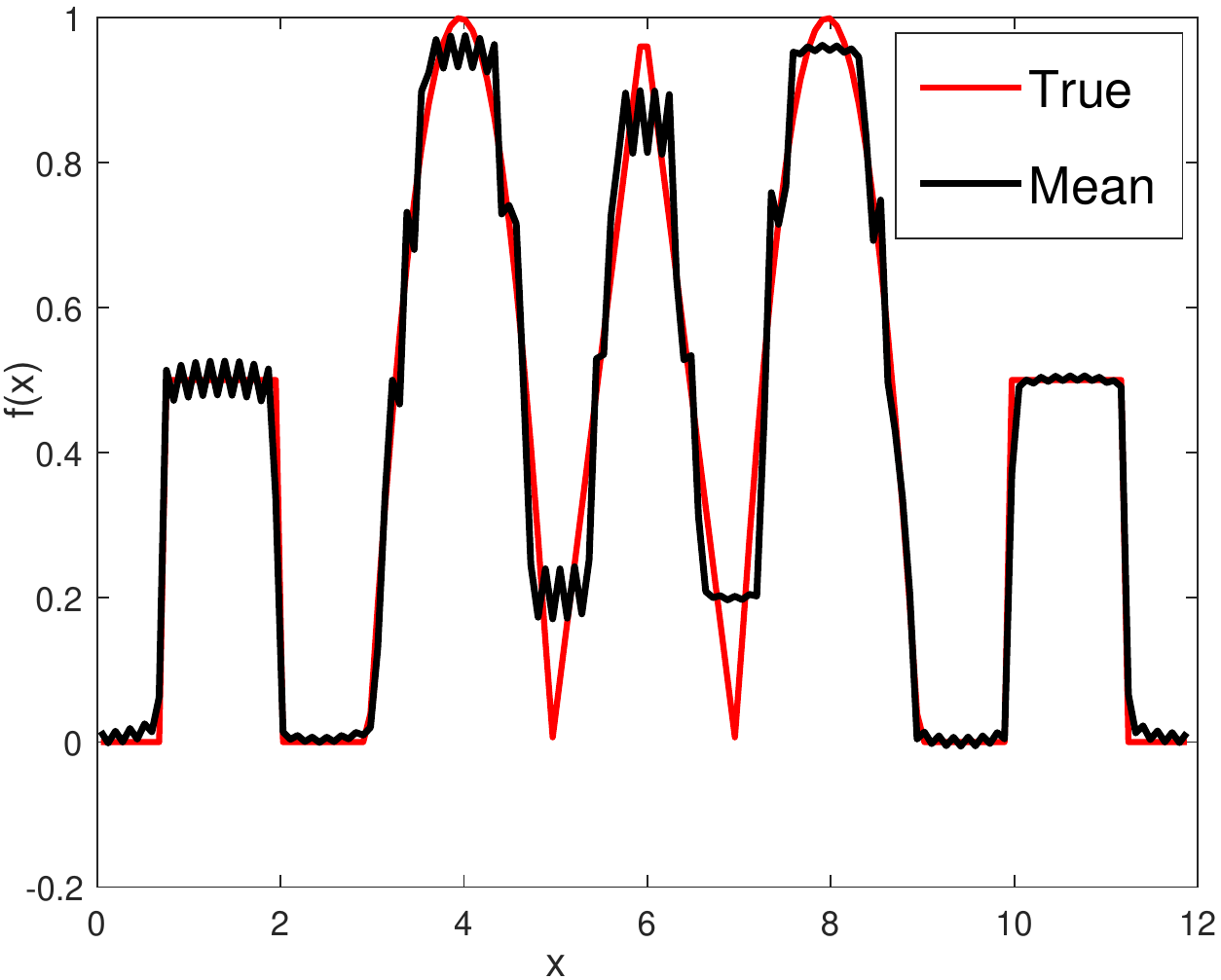}\\
\includegraphics[width=0.45\textwidth, height=0.12\textwidth]{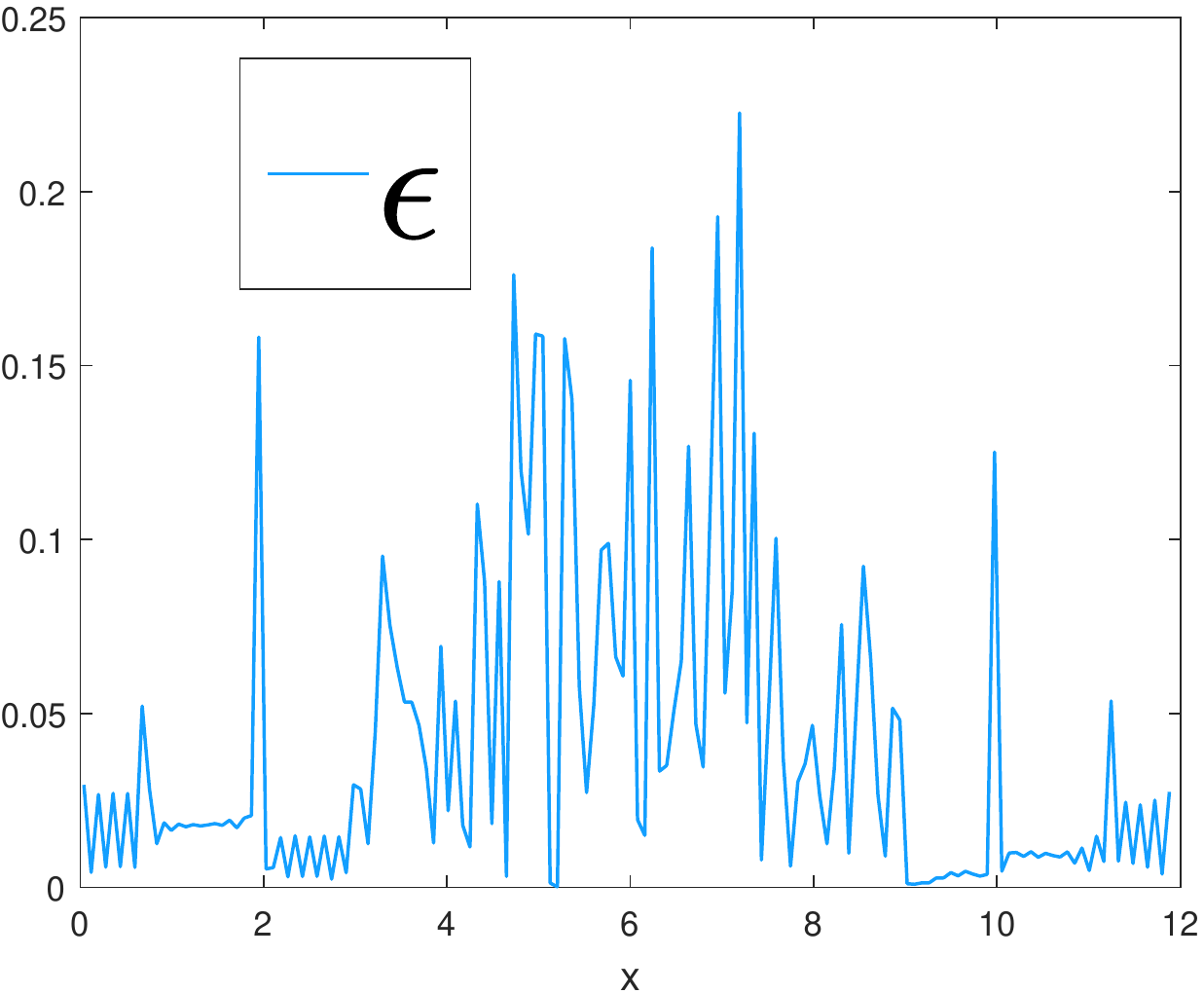}&
\includegraphics[width=0.45\textwidth, height=0.12\textwidth]{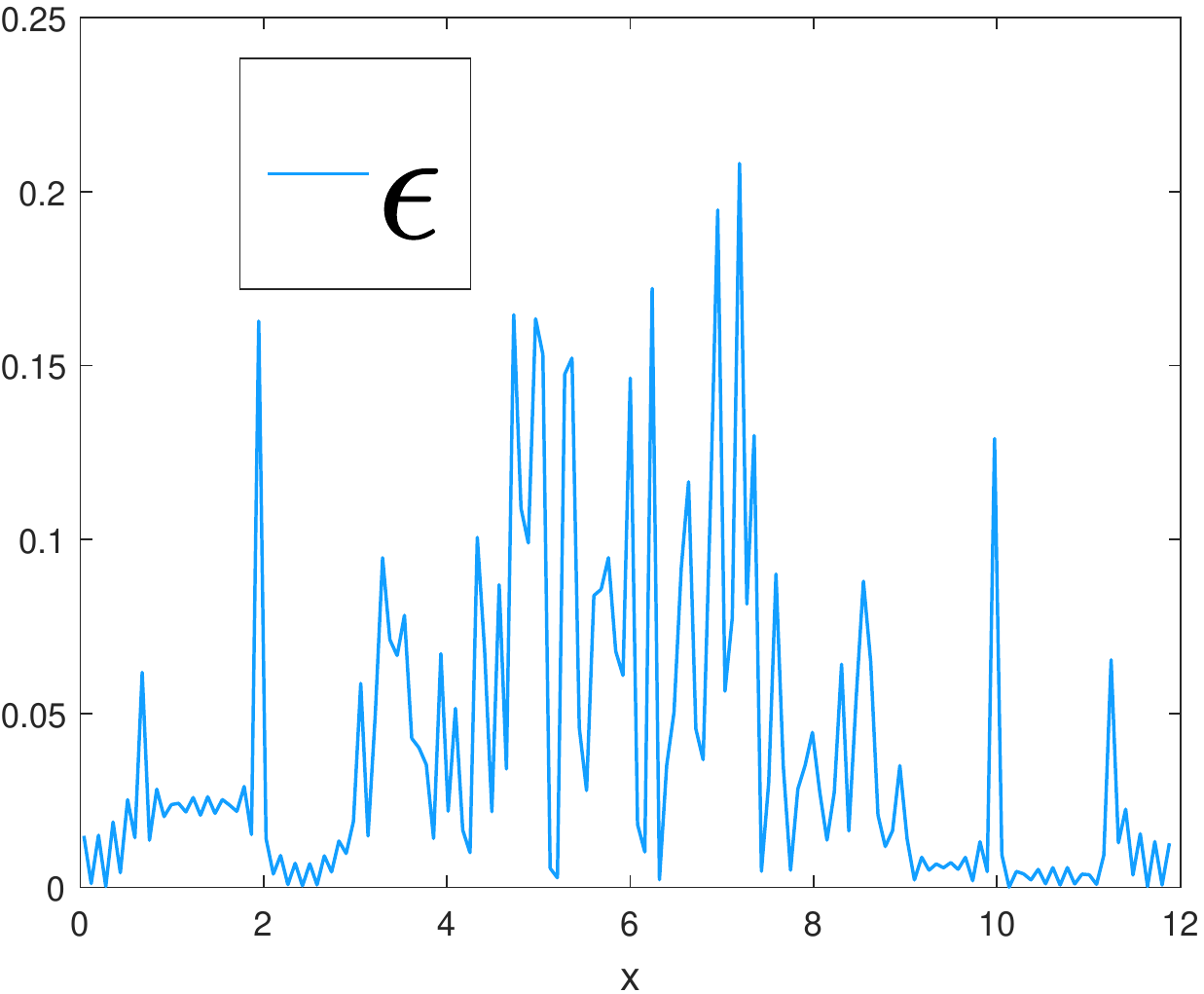}  \vspace{-5pt}\\
 FTG with $\alpha=0.95$ &FTG with $\alpha=0.99$ \\ \hline
 \end{tabular}
\caption{Reconstruction results for the heat source $f(\mathbf{x})$.
The posterior mean and absolute error $\epsilon$ using the linear diagonal map-based independence sampler Algorithm \ref{alg:diagonal_map} for FTG with $\alpha=0.5,\;0.9,\;0.95,\;0.99$ and TG prior.
}
\label{ex2_result_01}
\end{figure}

Another situation is for the FTG prior with $0<\alpha\le1$.
The posterior mean and absolute error $\epsilon$ using FTG prior with $\alpha=0.5,\;0.9,\;0.95,\;0.99$ are plotted in Figure \ref{ex2_result_01} and the relative error $RelErr$ are listed in Table \ref{table_heat01}.
From this table and figure, we conclude that the reconstruction results using FTG prior gradually converges to that using TG prior with the $\alpha \to 1^-$.
In Table \ref{table_heat01}, we can find that the FTG with $\alpha=0.2$ has lower $RelErr$ value than that with $\alpha=0.5$.
It seems that the selected Gaussian prior $\mu_0$ is appropriate for this example in the process of solving the numerical optimization, which are consistent with the above the example.
The Gaussian prior $\mu_0$ will play main role in the numerical results for the a small fractional order FTG prior.
And from this table, the $RelErr$ values generated by the FTG prior are approaching to that of TG prior as the fractional order $\alpha$ tends to $1^-$.
However, from absolute error curve in the Figure \ref{ex2_result_01}, the reconstructed results with FTG prior for $\alpha=0.99$ outperforms that with TG prior in some details such as the corner points.
\begin{table}[htbp]
\centering
\caption{The $RelErr$ values of the inverse source results using FTG for $0<\alpha\le1$ and TG prior.}
\begin{tabular}{cccccccccc}
  \toprule
          & TG & $\alpha=0.2$ & $\alpha=0.5$ & $\alpha=0.8$ & $\alpha=0.9$ & $\alpha=0.95$ & $\alpha=0.99$ \\
  \midrule
  $RelErr$ & 0.1054&0.1328& 0.1445& 0.1220 & 0.1202 & 0.1152 &0.1133\\
    \bottomrule
  \end{tabular}
  \label{table_heat01}
\end{table}

We now compare the efficiency performance of the linear diagonal map-based independence sampler in Algorithm \ref{alg:diagonal_map} and standard pCN in Algorithm \ref{alg:pCN}.
In this example, the stepsize $\beta$ is chosen so that the resulting acceptance probability is in the range $20\%-30\%$ and all other parameters of pCN take the same values as our diagonal map-based independence sampler method.
And we draw $10^6$ samples from the posterior with the first $0.5\times 10^6$ samples discarded as burn-in period for the standard pCN.

\begin{figure}[htbp]
 \centering
 \begin{tabular}{@{}c@{}c@{}c@{}}
\hspace{30pt} indep. sampler with a diagonal map&\hspace{30pt}standard pCN \vspace{10pt}\\
\rotatebox{90}{\hspace{5pt}FTG with $\alpha=0.9$}\hspace{15pt}\includegraphics[width=0.45\textwidth, height=0.2\textwidth]{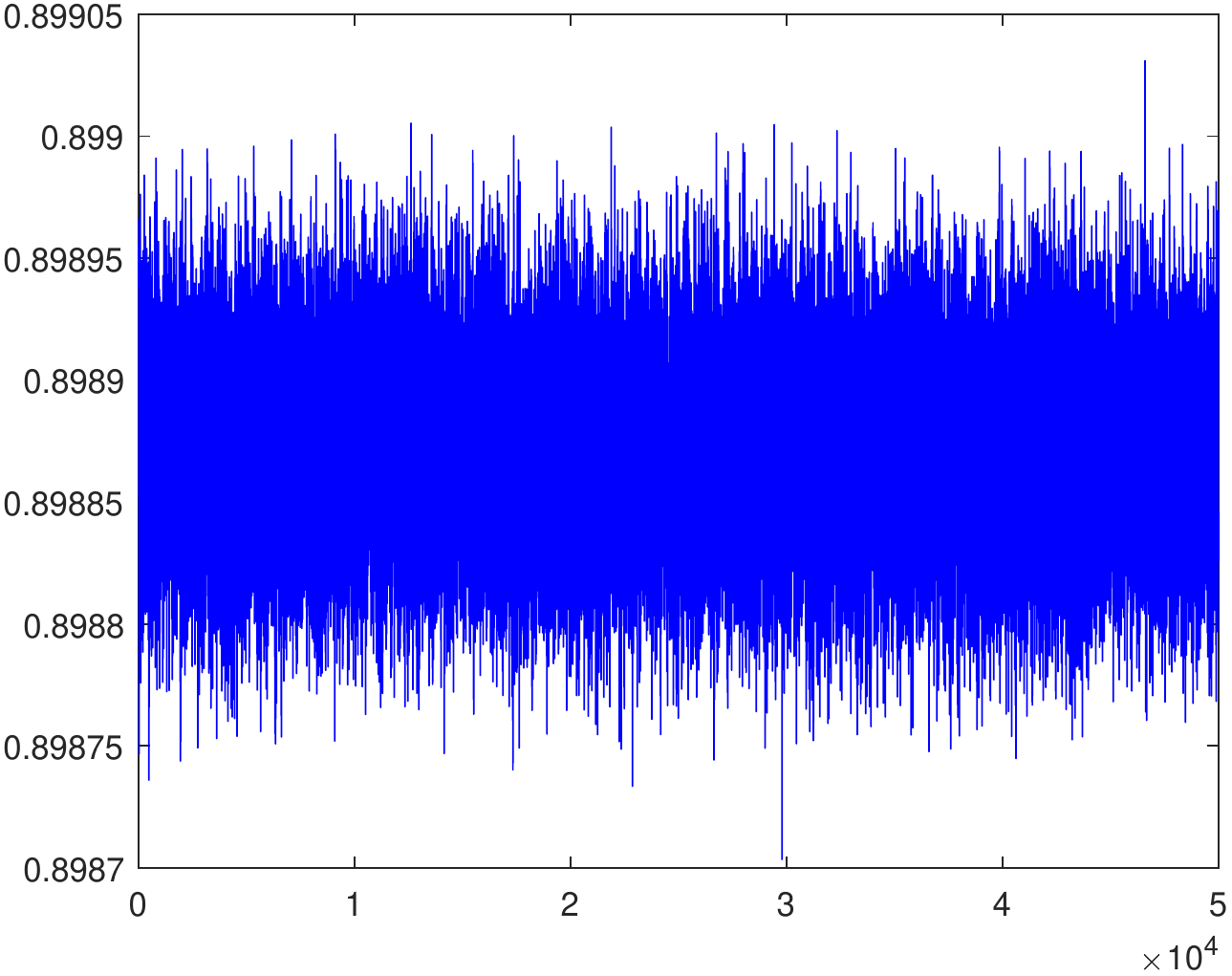}& \hspace{10pt}
\includegraphics[width=0.45\textwidth, height=0.2\textwidth]{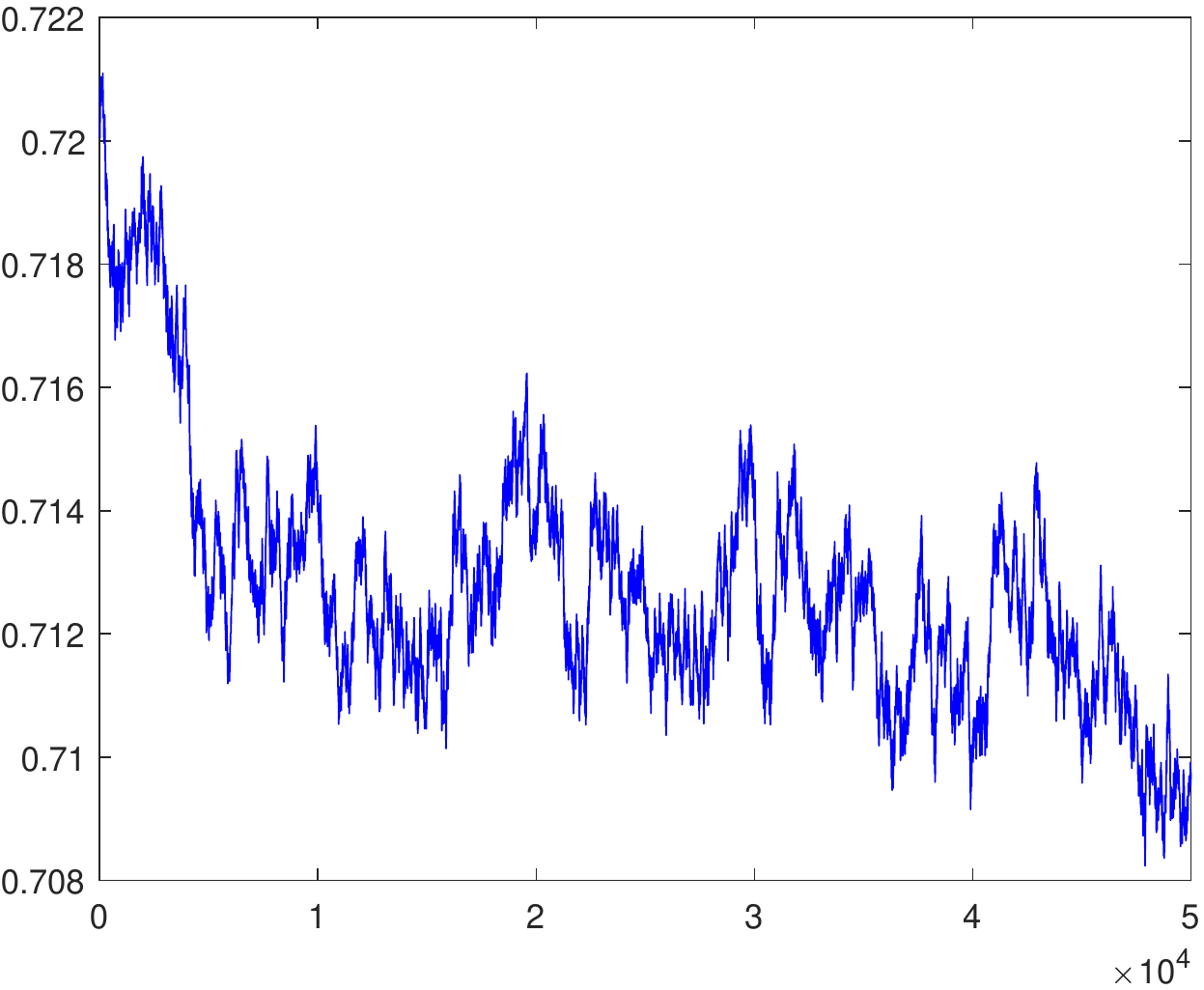}\\
\rotatebox{90}{\hspace{40pt}TG}\hspace{20pt}\includegraphics[width=0.45\textwidth, height=0.2\textwidth]{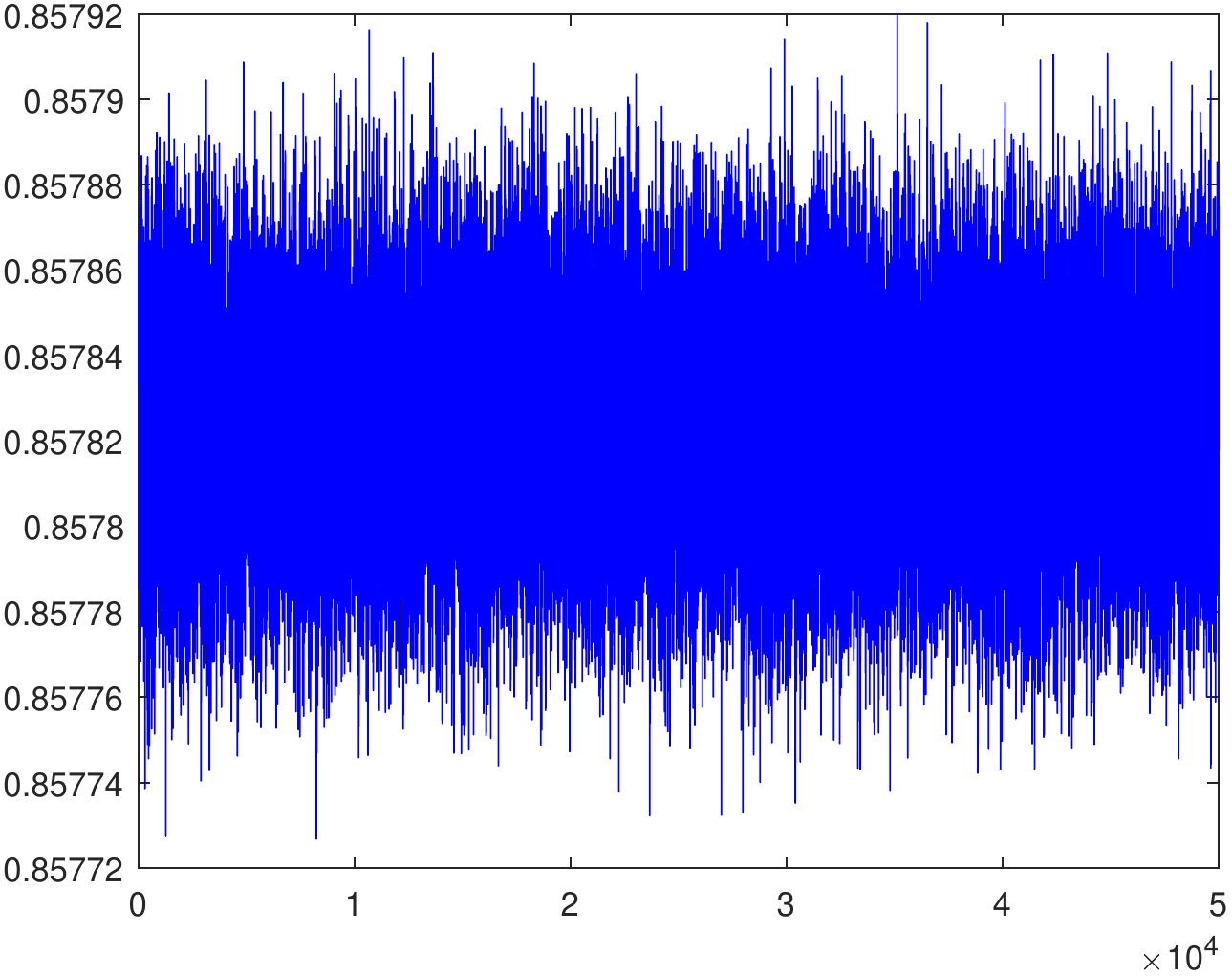}& \hspace{10pt}
\includegraphics[width=0.45\textwidth, height=0.2\textwidth]{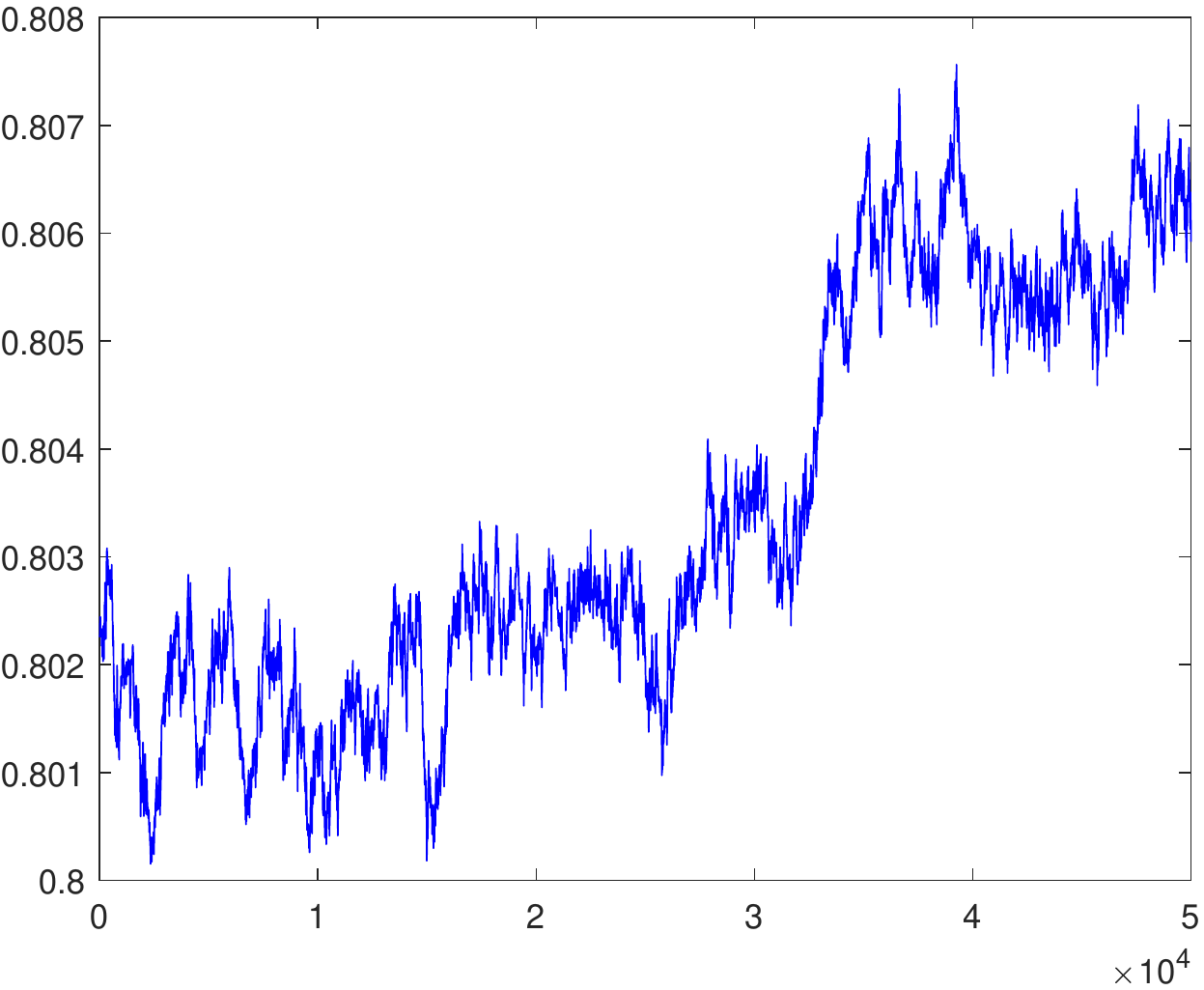}\\
\rotatebox{90}{\hspace{5pt}FTG with $\alpha=1.1$}\hspace{15pt}\includegraphics[width=0.45\textwidth, height=0.2\textwidth]{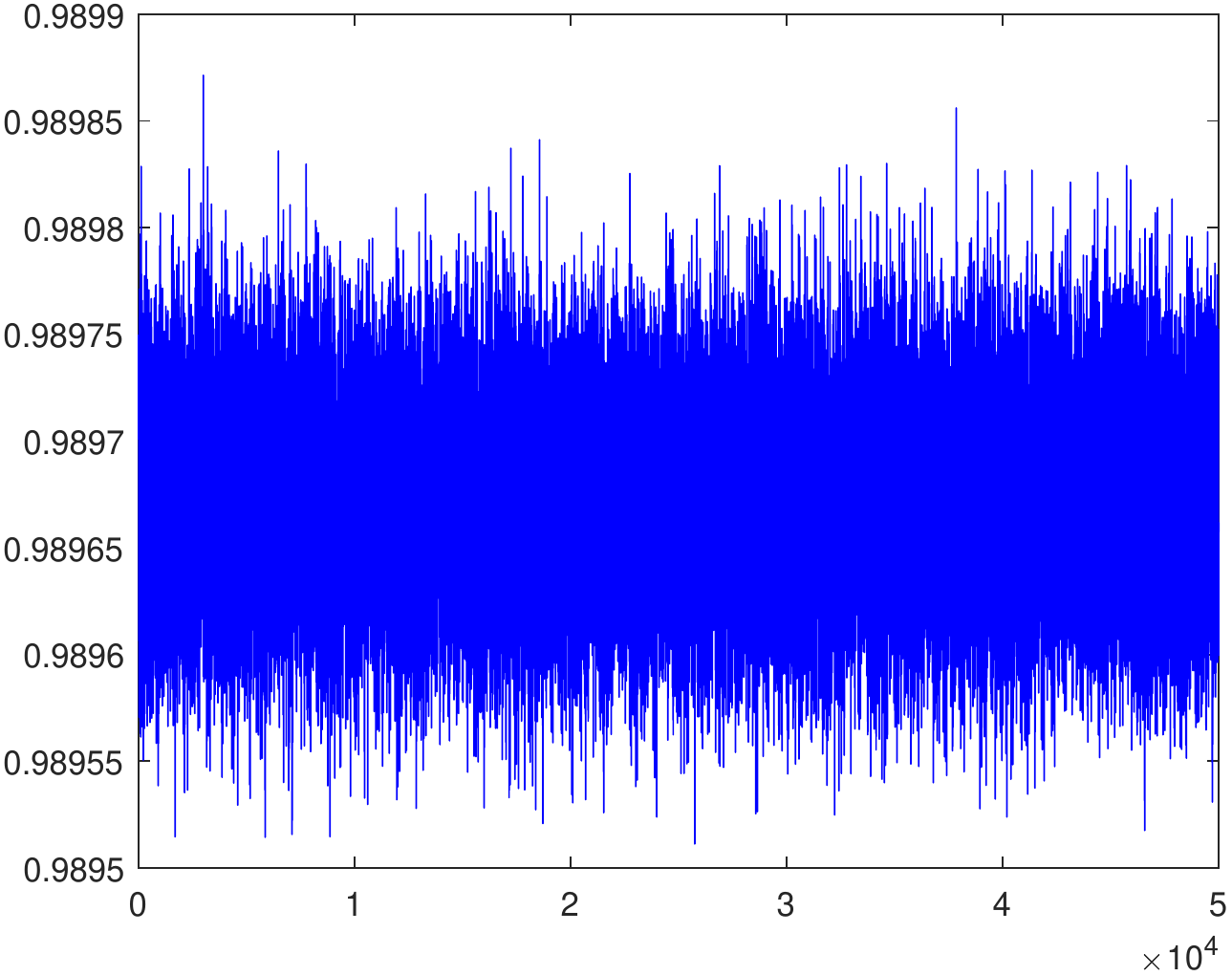}& \hspace{10pt}
\includegraphics[width=0.45\textwidth, height=0.2\textwidth]{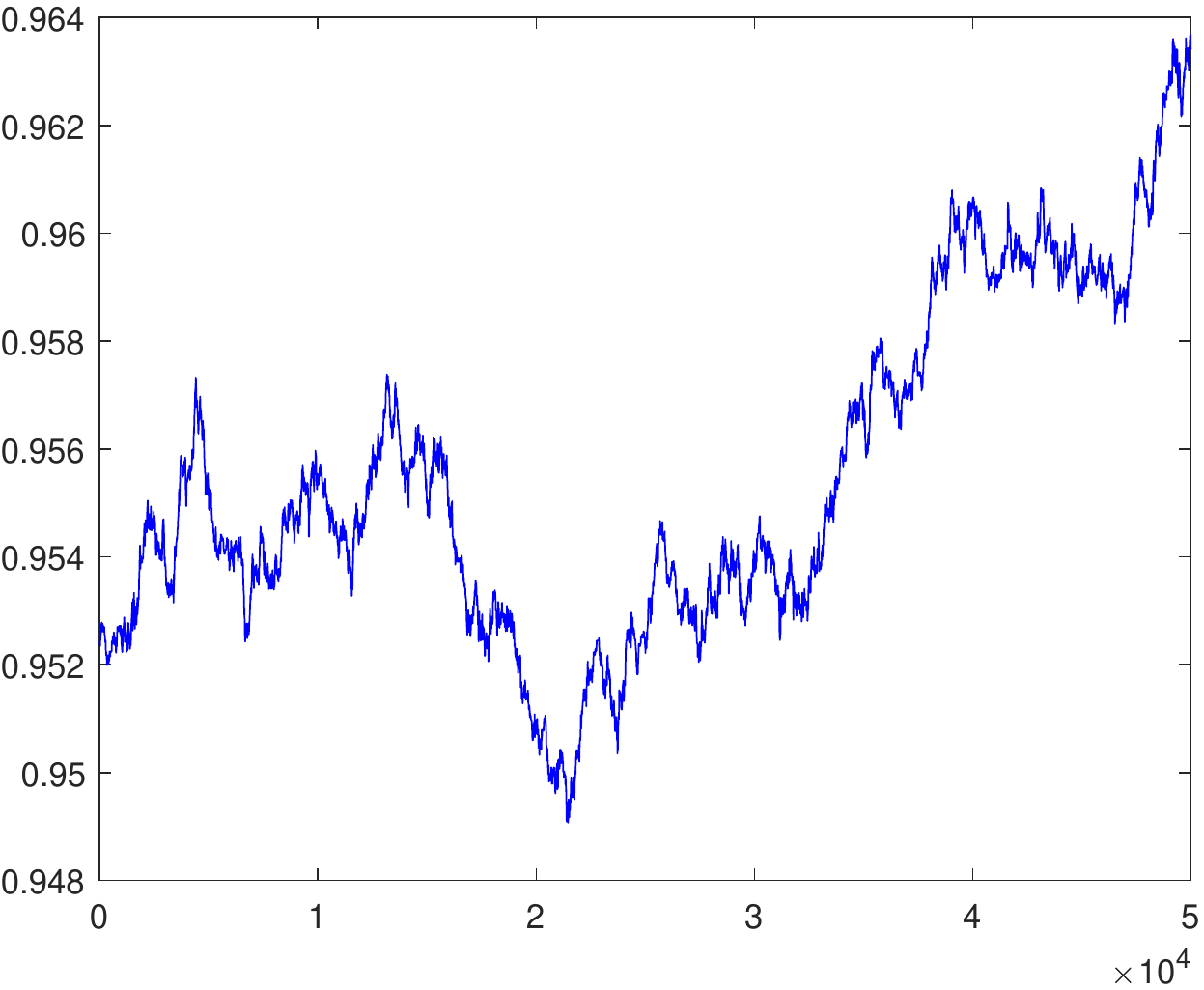}\\
 \end{tabular}
\caption{Inverse source identification problems. Trace plots of our independence sampler based on a linear diagonal map (the first column) and standard pCN (the second column) for the 75th component in chains (i.e. at $x=5.9603$) using the FTG with $\alpha=0.9,\; 1.1$ and TG prior.
}
\label{ex2_trace}
\end{figure}

\begin{figure}[htbp]
 \centering
 \begin{tabular}{@{}c@{}c@{}c@{}}
\hspace{30pt} indep. sampler with a diagonal map&\hspace{30pt}standard pCN \vspace{10pt}\\
\rotatebox{90}{\hspace{35pt}FTG with $\alpha=0.9$}\hspace{15pt}\includegraphics[scale=0.5]{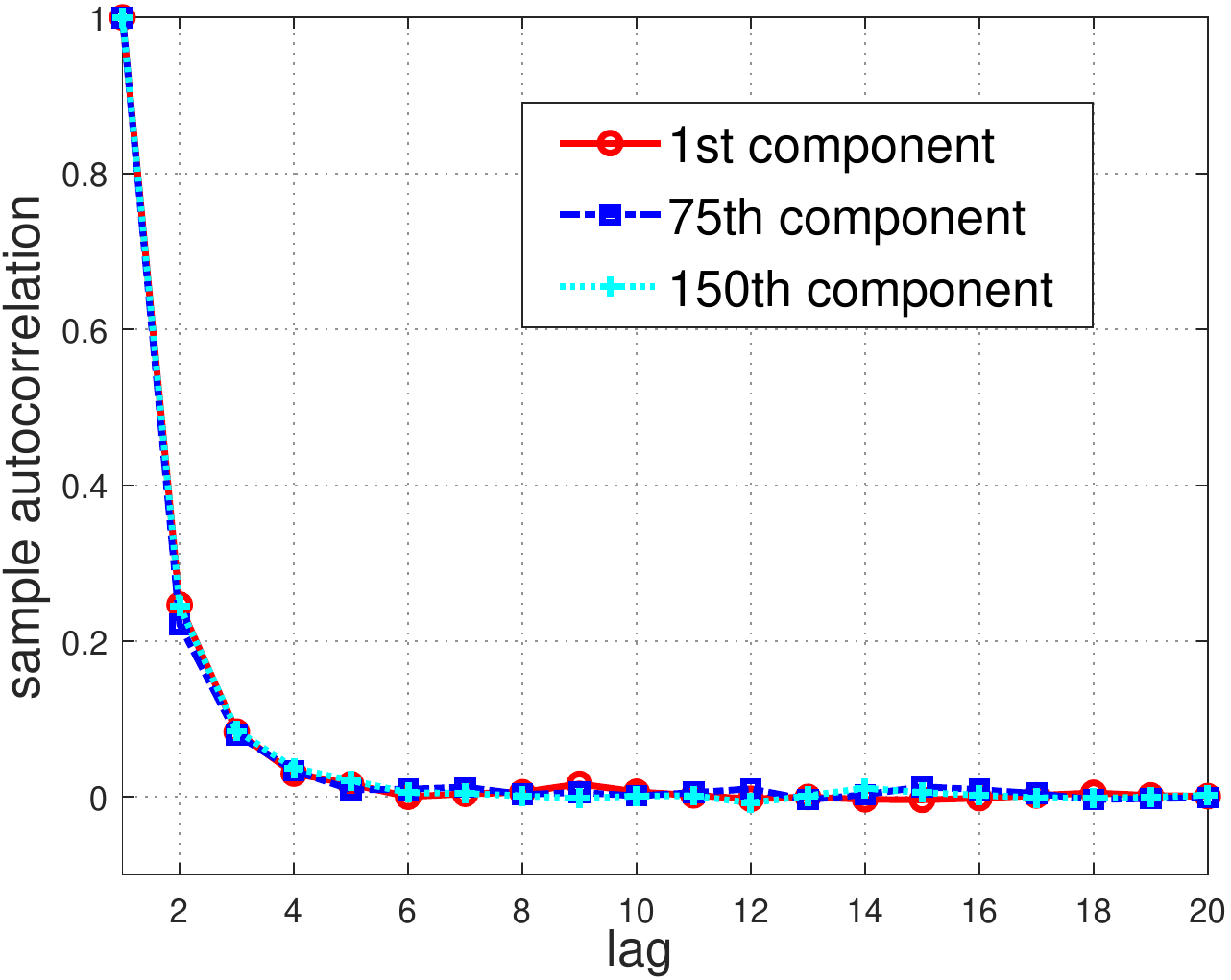}& \hspace{20pt}
\includegraphics[scale=0.5]{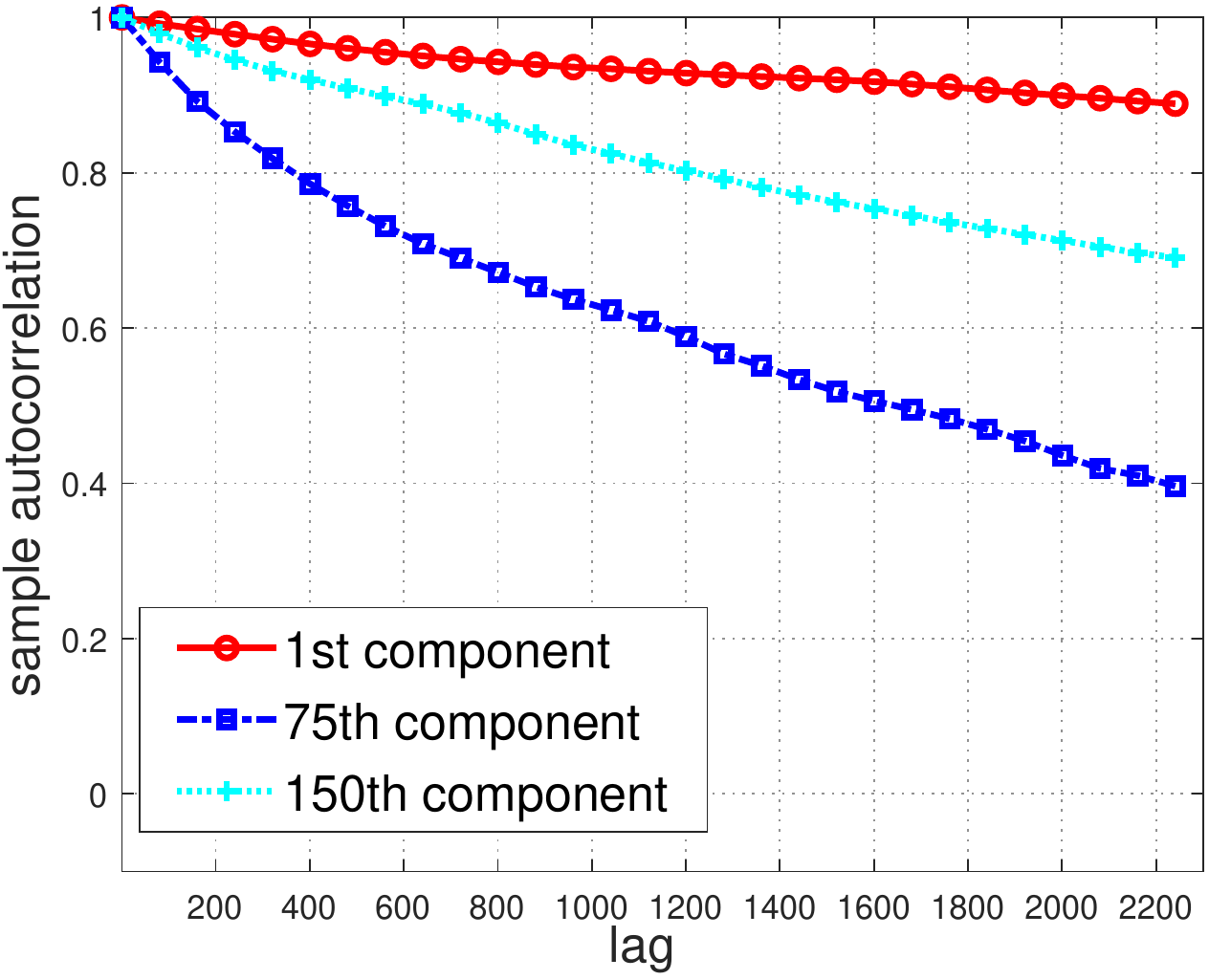}\\
\rotatebox{90}{\hspace{70pt}TG}\hspace{20pt}\includegraphics[scale=0.5]{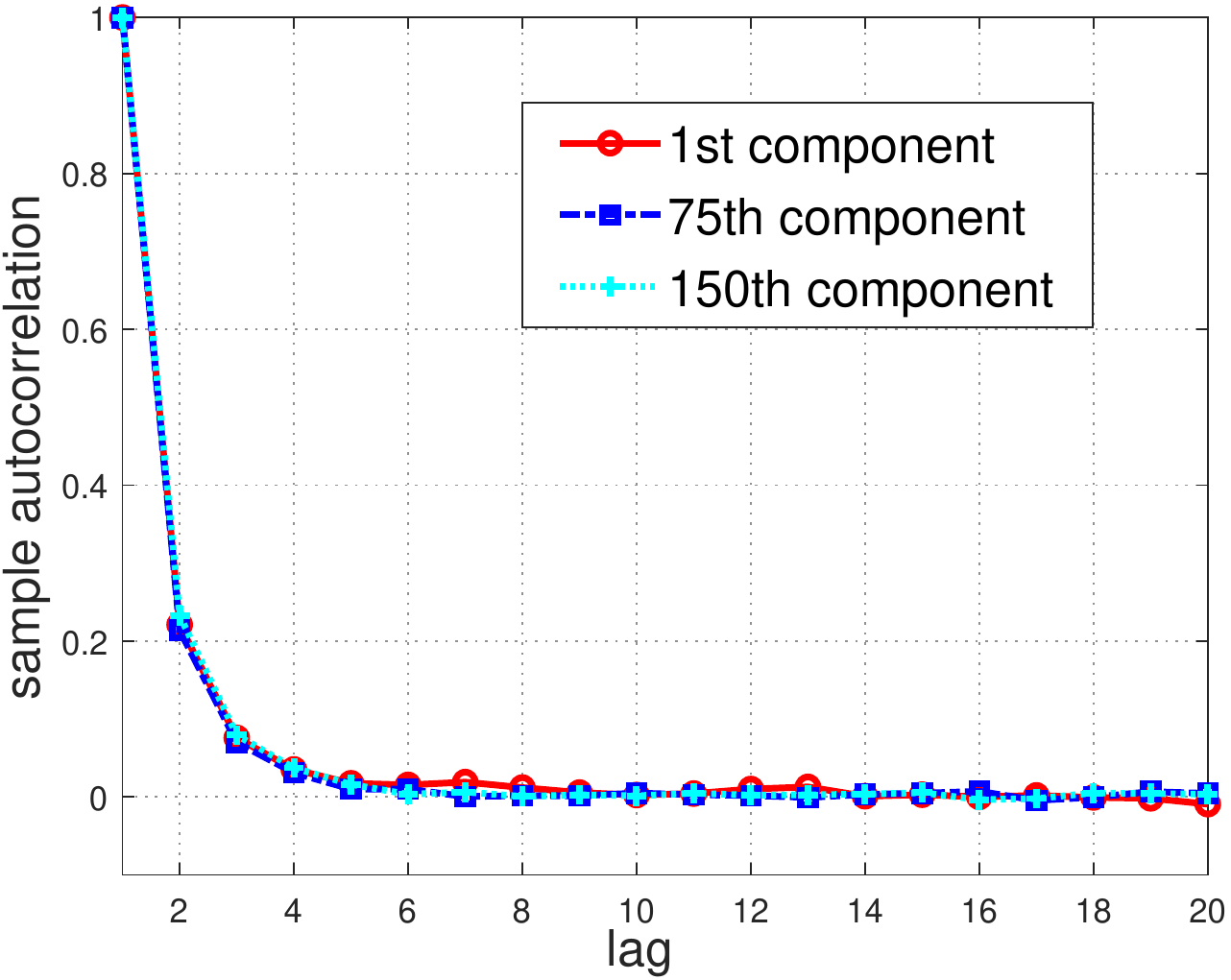}& \hspace{20pt}
\includegraphics[scale=0.5]{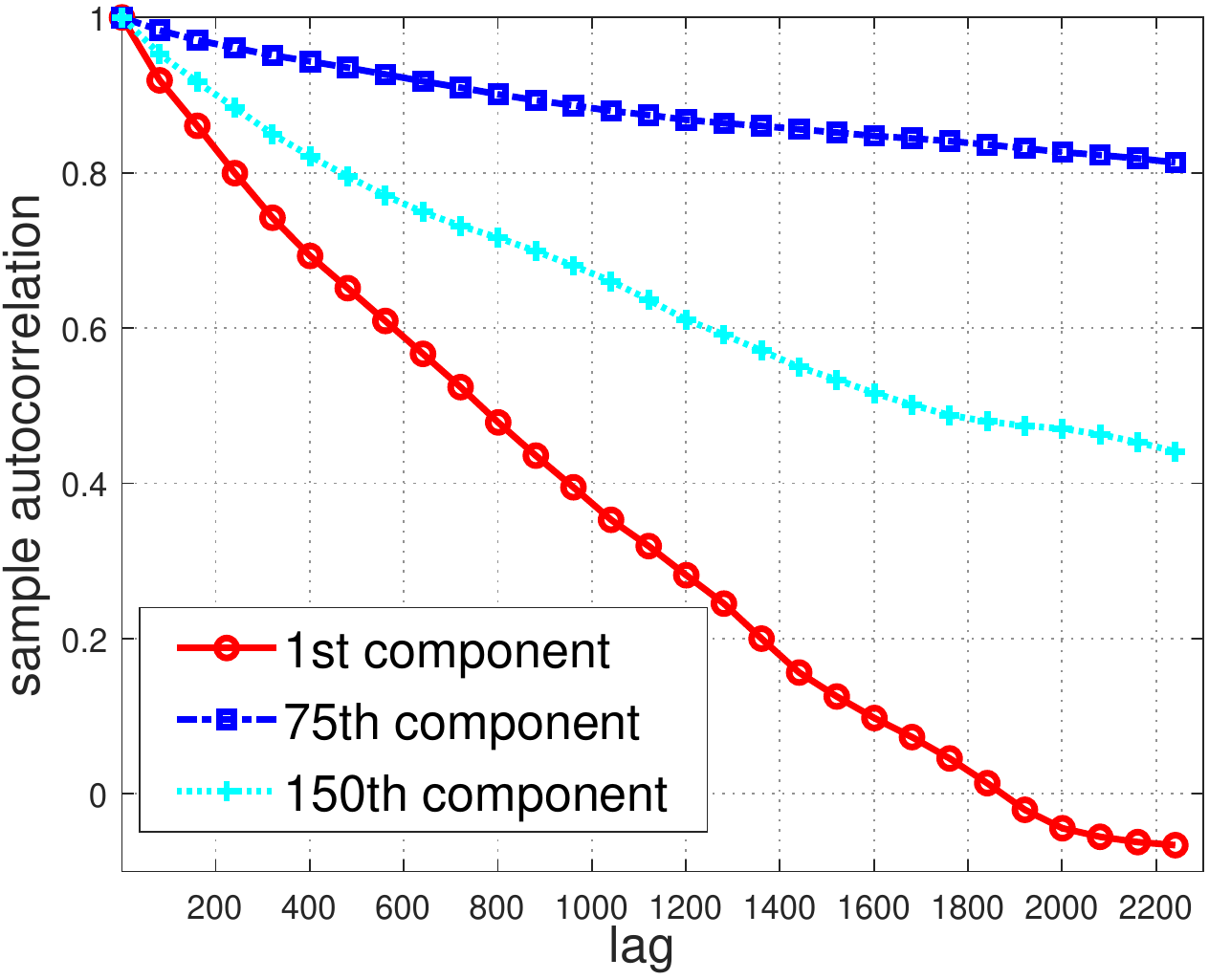}\\
\rotatebox{90}{\hspace{35pt}FTG with $\alpha=1.1$}\hspace{15pt}\includegraphics[scale=0.5]{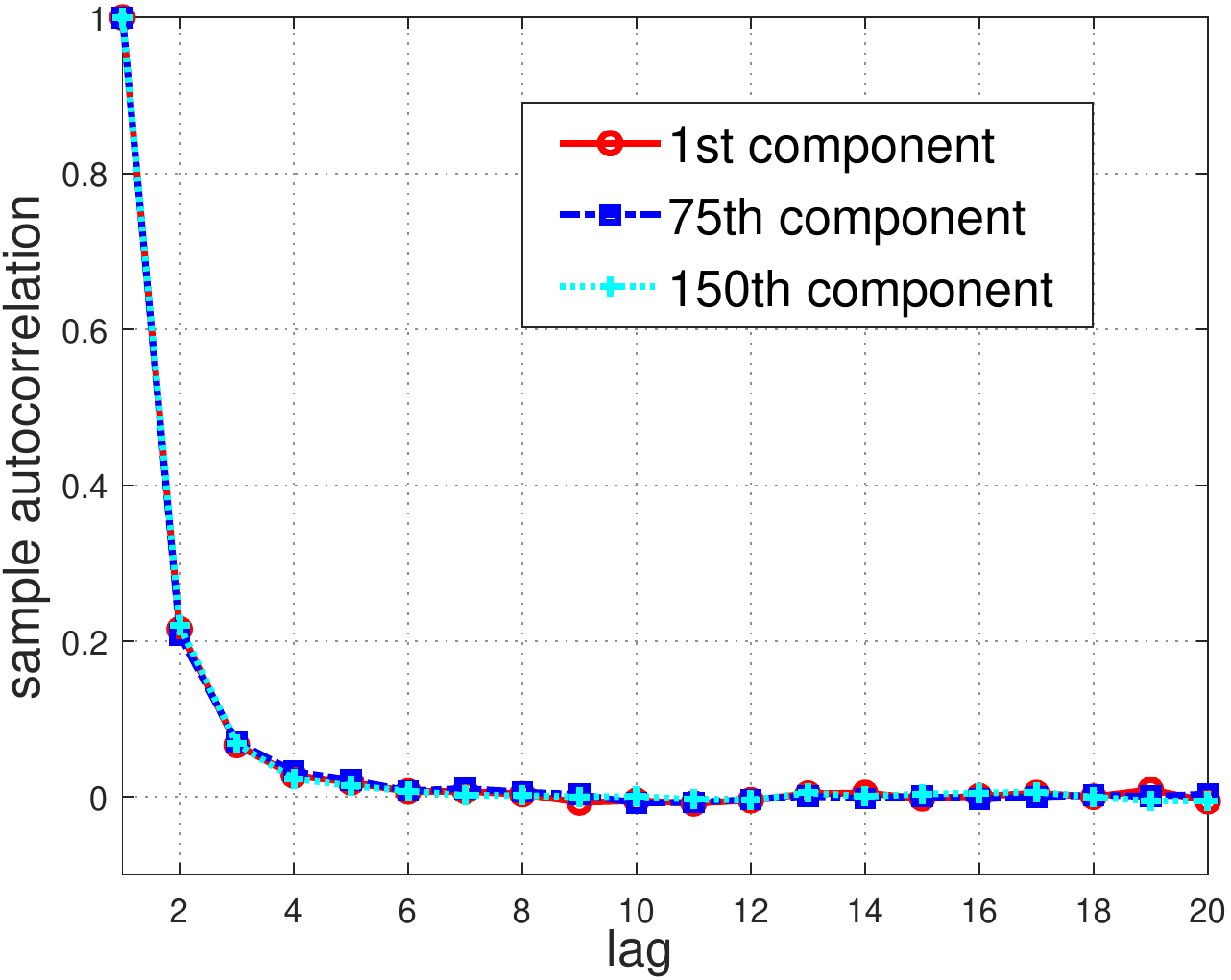}& \hspace{20pt}
\includegraphics[scale=0.5]{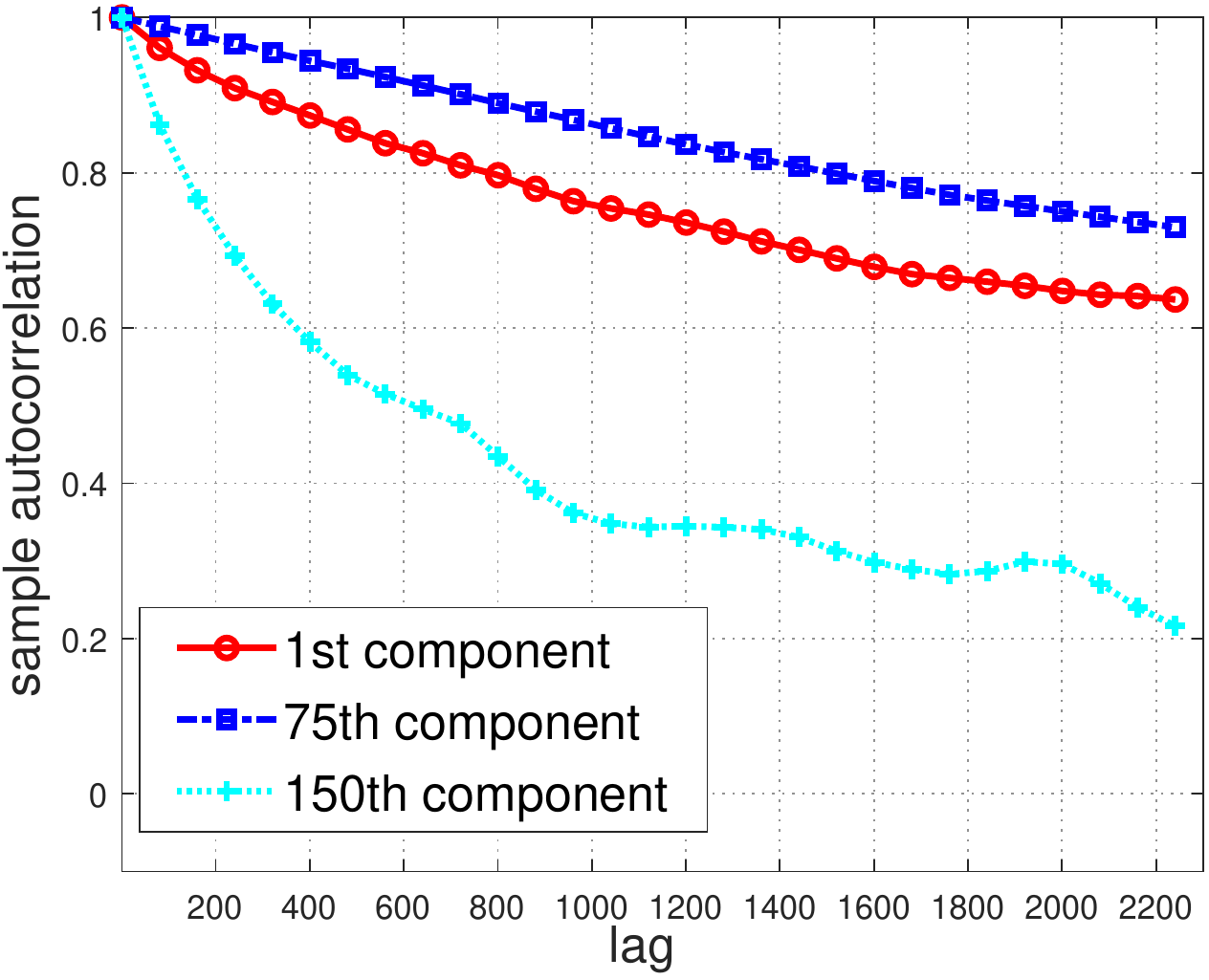}\\
 \end{tabular}
\caption{Inverse source identification problems. The autocorrelation functions corresponding to our independence sampler (the first column) and standard pCN (the second column) for the 1st, 75th, and 150th component in chains (i.e. at $x=0.0795,\; 5.9603,\; 11.9205$) using the FTG with $\alpha=0.9,\; 1.1$ and TG prior.
}
\label{ex2_acf}
\end{figure}

First, in Figure \ref{ex2_trace} we show the trace plots of the two methods for the unknown function $f(x)$ at $x=5.9603$ using the FTG with $\alpha=0.9,\; 1,1$ and TG prior.
High mixing rate in the plots indicate that successive iterations are highly independent and that the series of iterations have converged.
The trace plots, in this respect, indicate that our diagonal map-based independence sampler method achieves a much faster mixing rate than the standard pCN.
Next we compute the autocorrelation functions (ACF) of various quantities with the samples drawn by the two methods to further compare the efficiency.
We consider the points at $x=0.0795,\; 5.9603,\; 11.9205$ and plot the ACF results in Figure \ref{ex2_acf} using the FTG with $\alpha=0.9,\; 1,1$ and TG prior.
From this figure, the ACF for all three points of the chains generated by the diagonal-map based independence sampler decay faster than that corresponding to the standard pCN, which also suggests that the our method achieves a significantly better performance.
The average acceptance rate of the independence sampler is about $77\%$ and that of the standard pCN is about $31\%$.

\begin{figure}[htbp]
\centering
\subfigure[FTG with $\alpha=0.9$]{
\includegraphics[scale=0.5]{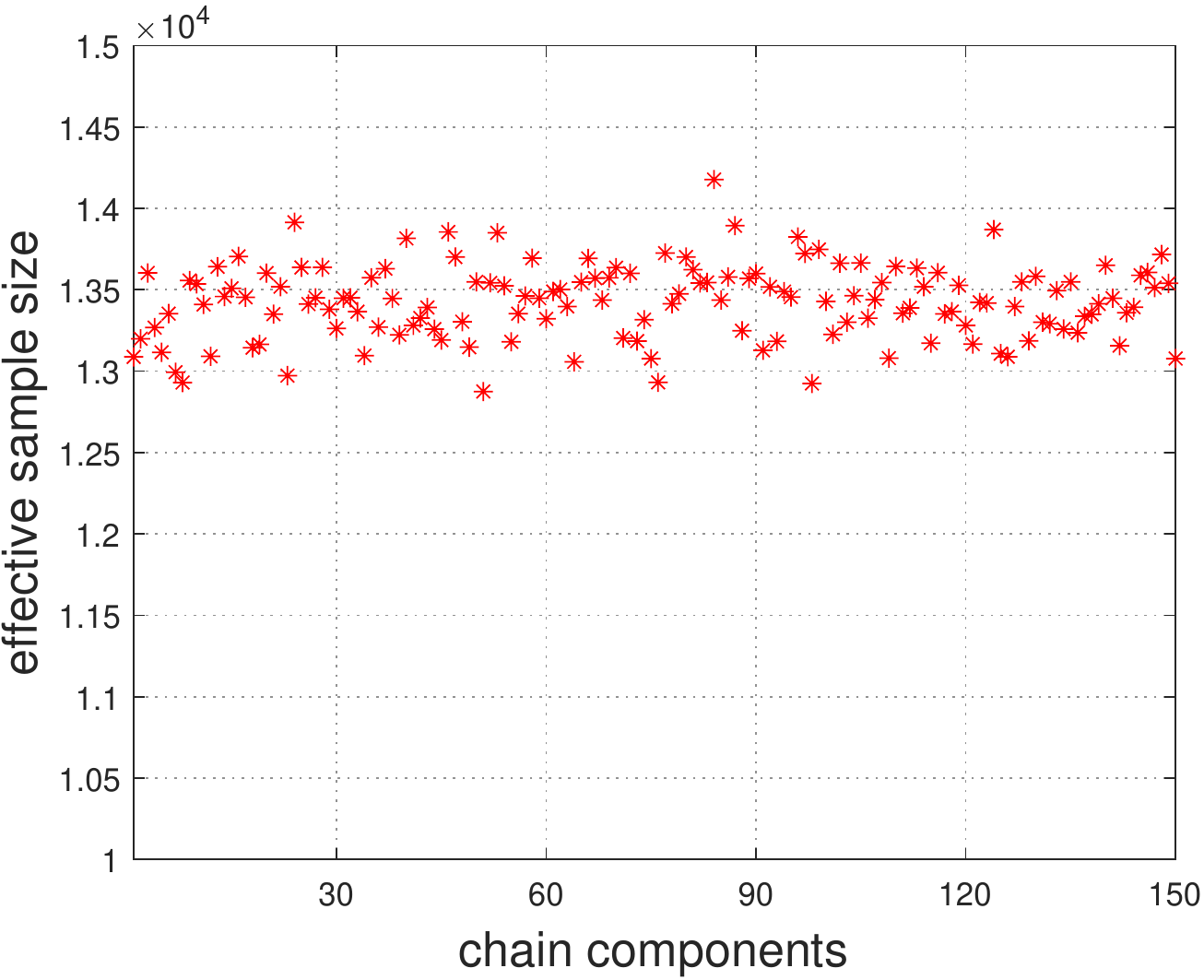}
}
\subfigure[TG]{
\includegraphics[scale=0.5]{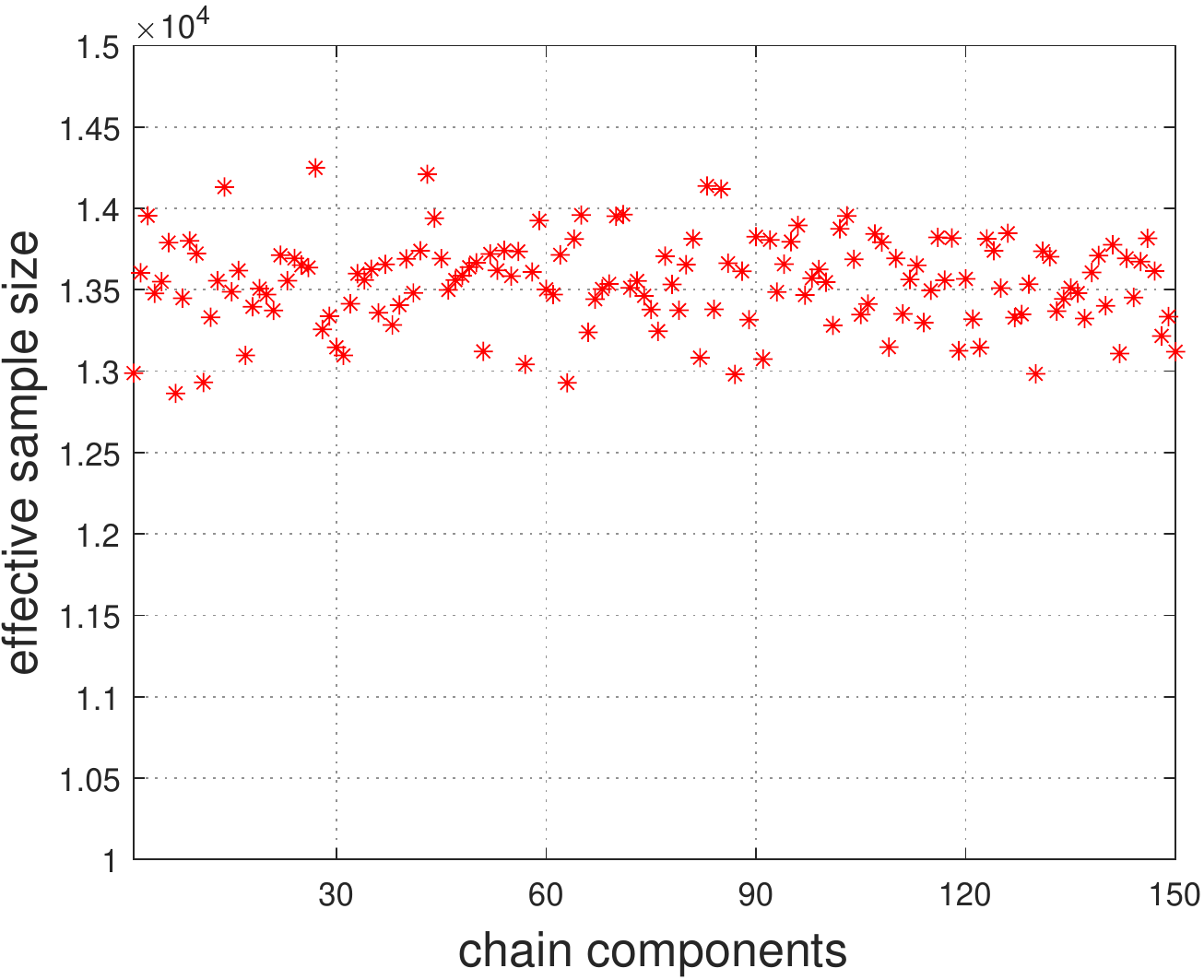}
}

\subfigure[FTG with $\alpha=1.1$]{
\includegraphics[scale=0.5]{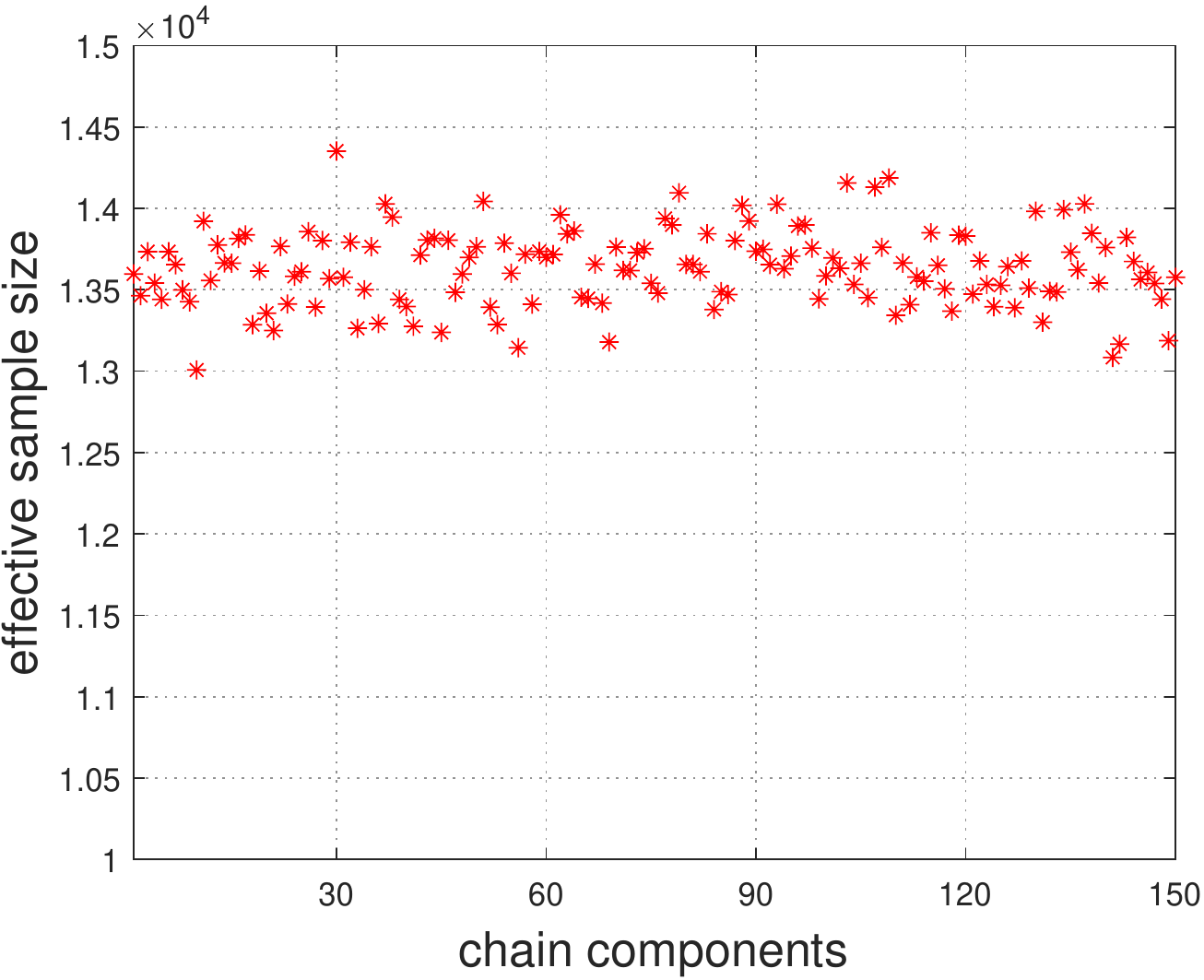}
}
\caption{Inverse source identification problems. Effective sample size of our independence sampler based on a linear diagonal map after $5\times 10^4$ MCMC iterations using the FTG with $\alpha=0.9,\; 1.1$ and TG prior.
}
\label{ex2_ess}
\end{figure}

Another common measure of the sampling efficiency of MCMC is the effective sample size (ESS) \cite{1998_Kass}.
The ESS is computed by
\[
\mathrm{ESS}=\frac{K}{1+2\varrho},
\]
where $K$ is the total sample size and $\varrho$ is the integrated autocorrelation time,  and an estimate of the number of effectively independent draws in the chain can be obtained by the ESS.
We compute the ESS of the unknown function $f(x)$ at each grid point through our diagonal map-based independence sampler method and show the results in Figure \ref{ex2_ess}.
The results show that our independence sampler algorithm on average produces significantly more effective independent samples, which agrees with faster decay of the ACF of this algorithm in Figure \ref{ex2_acf}.

\subsection{Limited computed tomography reconstruction}\label{2d_ex1}
We consider a classical inverse problem of X-ray computed tomography (CT), where X-rays travel from sources to detectors passing an object of interest.
The intensities from multiple sources are measured at the detectors, the goal is to reconstruct the density of the object.

\begin{figure}[htbp]
\centering
\includegraphics[width=0.45\textwidth, height=0.4\textwidth]{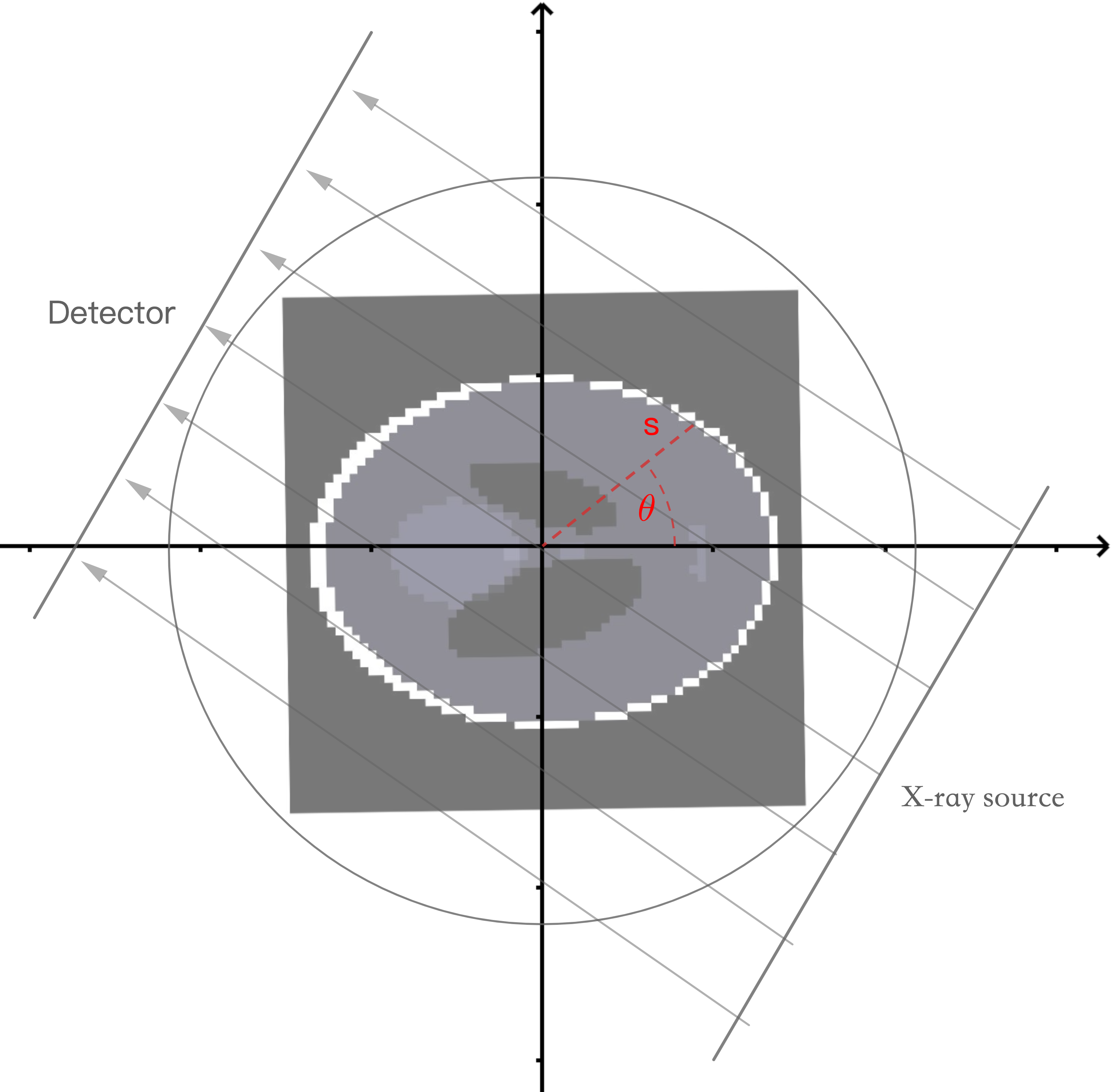}
\caption{X-ray tomography problem.}
\label{Pbg}
\end{figure}

\subsubsection{Problem setup}
In this case, the forward operator is characterized by the Radon transform \cite{radon_1986, Mueller_2012}, which can model the attenuation of the x-ray when traveling from the source to the detector through the target, Shepp-Logan phantom, (cf. Figure \ref{Pbg}).
Let the angle $\theta \in [0,\pi]$, and denote by $\boldsymbol{\zeta}:=[\cos\theta, \sin\theta]^T$ the unit vector with angle $\theta$ with respect to the horizontal axis.
The Radon transform, denoted by $\mathfrak{R}$, of the function $f: \Omega \subset \mathbb{R}^2 \to \mathbb{R}$ then calculates the integral, which depends on the $\theta$ and a linear parameter $s\in \mathbb{R}$, along the line $L:=\{\mathbf{x}\in\mathbb{R}^2: \mathbf{x}\cdot\boldsymbol{\zeta}=s\}$

\begin{align}
\mathfrak{R}f(s,\theta)=\int_Lf(\mathbf{x})\mathrm{d}\mathbf{x}^{\perp},
\end{align}
where the $\mathrm{d}\mathbf{x}^{\perp}$ is the Lebesgue measure along the line $L$ and the $f$ is called an attenuation coefficient function \cite{Mueller_2012}.
Then the result is the logarithm of the ratio of the intensity between the source, $I_0$, and the detector, $I_d$, through the Beer's law
\begin{align}
\mathfrak{R}f(s,\theta)=-\log(\frac{I_d}{I_0})=\mathbf{y},
\end{align}
where the $\mathbf{y}$ is called sinogram, which is also called the measurement \eqref{Model} and can take the form
\begin{align}
\label{CTeq1}
\mathbf{y}=\mathfrak{R}f(s,\theta)+\boldsymbol{\eta}=
\begin{bmatrix}
\int_{L_1}f(\mathbf{x})\mathrm{d}s_1 \\[8pt]
\int_{L_2}f(\mathbf{x})\mathrm{d}s_2 \\
\vdots \\
\int_{L_n}f(\mathbf{x})\mathrm{d}s_n
\end{bmatrix}
+\boldsymbol{\eta},
\end{align}
where $\boldsymbol{\eta} \sim \mathcal{N}(0,\sigma^2I)$ is assumed to be the measurement errors and $n$ is sum of the number of the x-ray in all discrete angles.
If the domain of $f$ is discretized into a grid and the density (i.e. the $f$) is assumed to be constant within each grid, then the line integrals are approximated as
\begin{align}
\label{CTeq2}
\int_{L_i}f(\mathbf{x})\mathrm{d}s\approx \sum^{d}_{j=1}a_{ij}\mathbf{f}_j,
\end{align}
where $d$ is the total number of grid cell and the $a_{ij}$ is the length of the intersection between line $L_i$ and the $j$th grid cell.
Then \eqref{CTeq1} and \eqref{CTeq2} yield a linear model $\mathbf{y}=A\mathbf{f}+\boldsymbol{\eta}$, where the matrix $A$ is defined by $A=(a_{ij})_{n\times d}$.

\subsubsection{Set up of inverse problems}
In this case, the target square Shepp-Logan phantom $f(\mathbf{x})$, defined on $\Omega=[-1,1]^2$, is divided uniformly into $64\times 64$ grid cell (or pixels), seeing the top left in Figure \ref{2d_result}, i.e. $d=4096$.
We consider the measurement data $\mathbf{y}=A\mathbf{f}+\boldsymbol{\eta}$ and $\boldsymbol{\eta}$ adds Gaussian noise with zero mean and standard deviation $0.0115$, which corresponds to $0.07\%$ noise with respect to the maximum norm of the output $A\mathbf{f}$.
Note that the measurement data $\mathbf{y}$ is computed from a twice finer ($128\times 128$) grid Shepp-Logan phantom (at same measurement angles but finer arrangement of X-rays) and interpolated to lower resolution, while we use a $64\times 64$ grid for inference.
The Gaussian prior $\mu_0$ is taken by zero mean and covariance matrix $C_0=10^{-5}I$ and the FTG prior with the fractional order $\alpha=0.5,\; 0.8,\; 0.9,\;1.1,\; 1.2,\; 1.5$.
The shape and rate parameter of Gamma distribution are taken by $k=2.55\times10^6,\;\vartheta=1$, respectively.

We then construct a transport map from the Gaussian distribution $\mu_0$ to the posterior \eqref{H_postprior}.
The numerical optimization problem \eqref{d_opt_problem} is performed with MATLAB's fmincon optimizer, where the step tolerance (StepTolerance) is set to $10^{-3}$; the SpecifyConstraintGradient and SpecifyOdjectiveGradient are set to true and we use $M=4096$ samples of the Gaussian distribution $\mu_0$ to approximate the expected value through the SAA in the objective function (see Section \ref{numercal opt}).
The transport map and the regularization parameter via the \eqref{de_reg} are then used to precondition for our independent sampler as shown in Algorithm \ref{alg:diagonal_map}.
And the proposal $\mu_{ref}$ is the Gaussian distribution with zero mean and the standard deviation $10^{-5}$ in the Algorithm \ref{alg:diagonal_map}.
We compare the posterior mean using our independent sampler with FTG to TG prior and the result of filtered back-projection (FBP) \cite{Mueller_2012} inversion technique.
The structural similarity (SSIM) \cite{Wang_2004}, which is good at measuring quality of the reconstructed images in terms of image structure, is also adopted here.

\subsubsection{Result}
\begin{figure}[htbp]
 \centering
 \begin{tabular}{|@{}c@{}|@{}c@{}|@{}c@{}|}
 \hline
\includegraphics[scale=0.5]{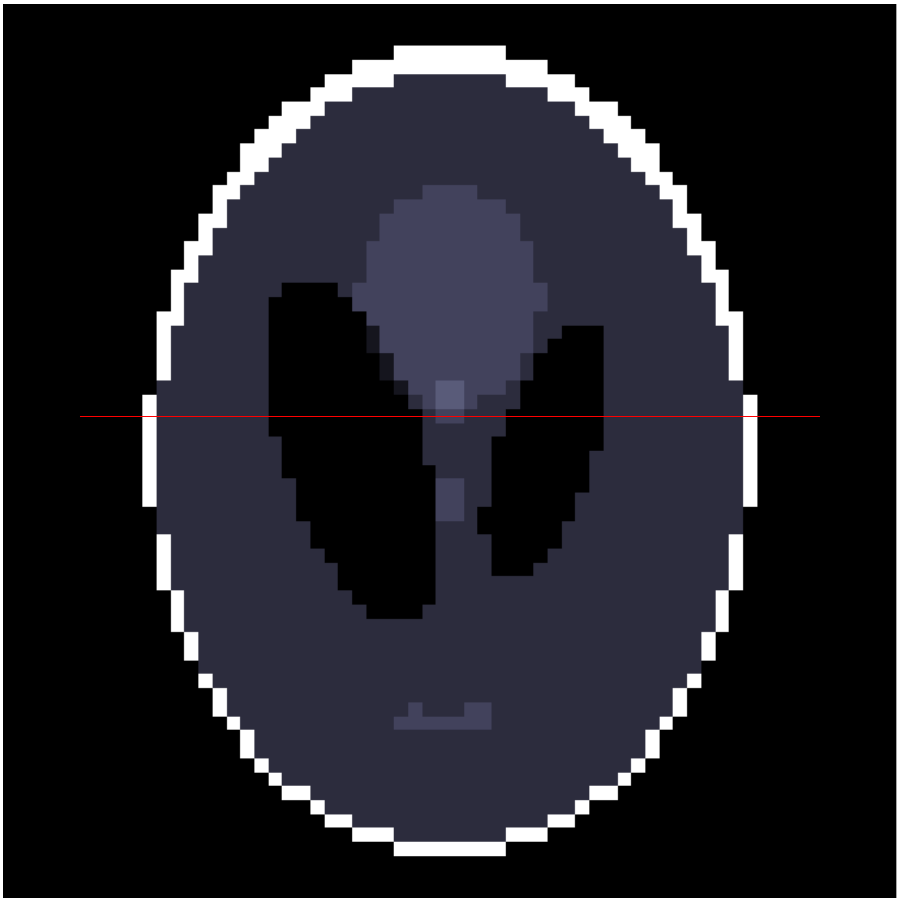}&
\includegraphics[scale=0.5]{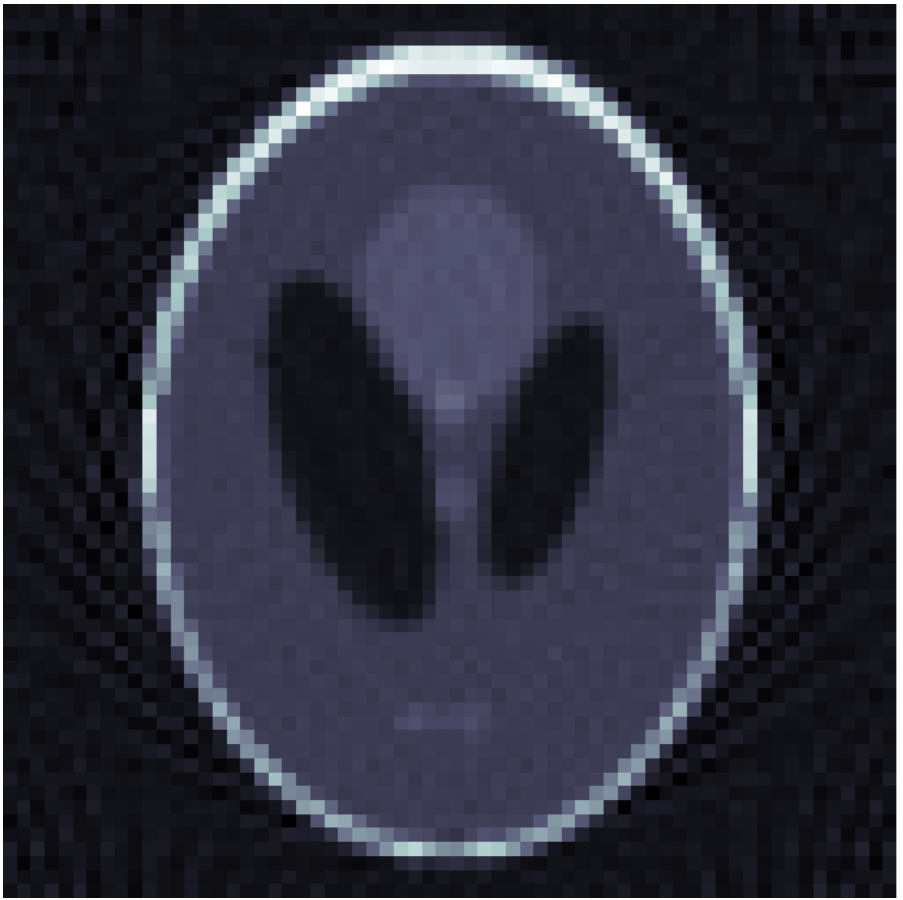}&
\includegraphics[scale=0.5]{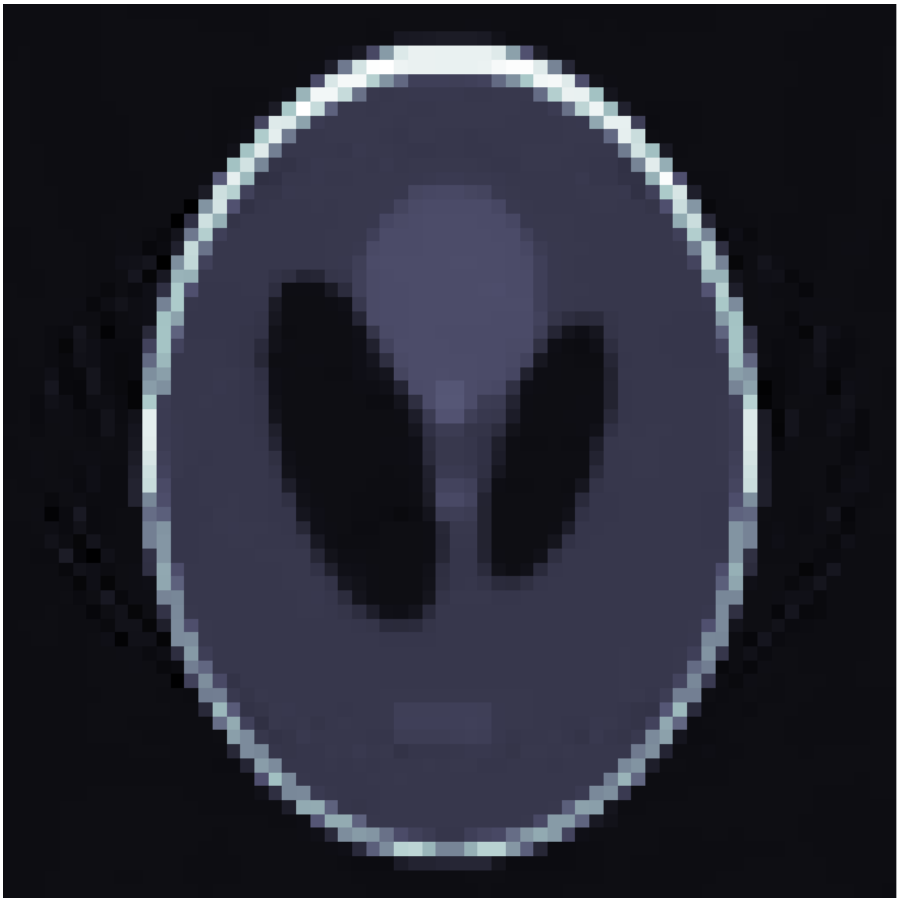} \\
\hspace{-15pt}Target&\hspace{-15pt}FBP&\hspace{-15pt}TG \\ \hline
\includegraphics[scale=0.5]{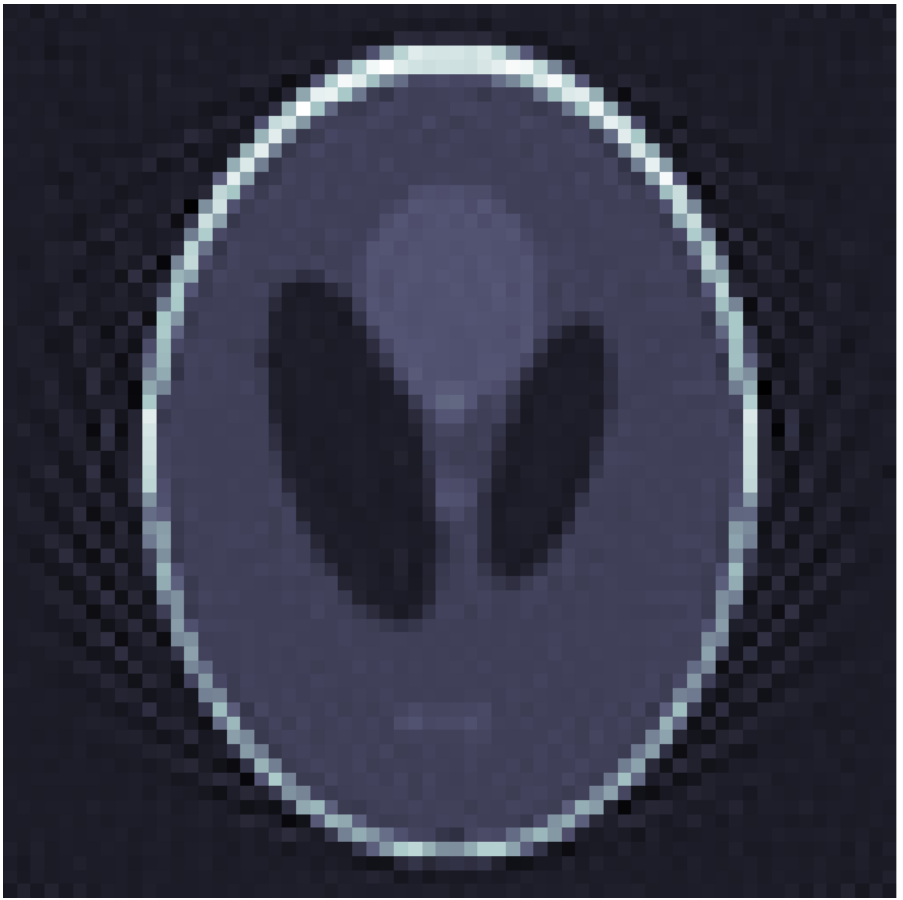}&
\includegraphics[scale=0.5]{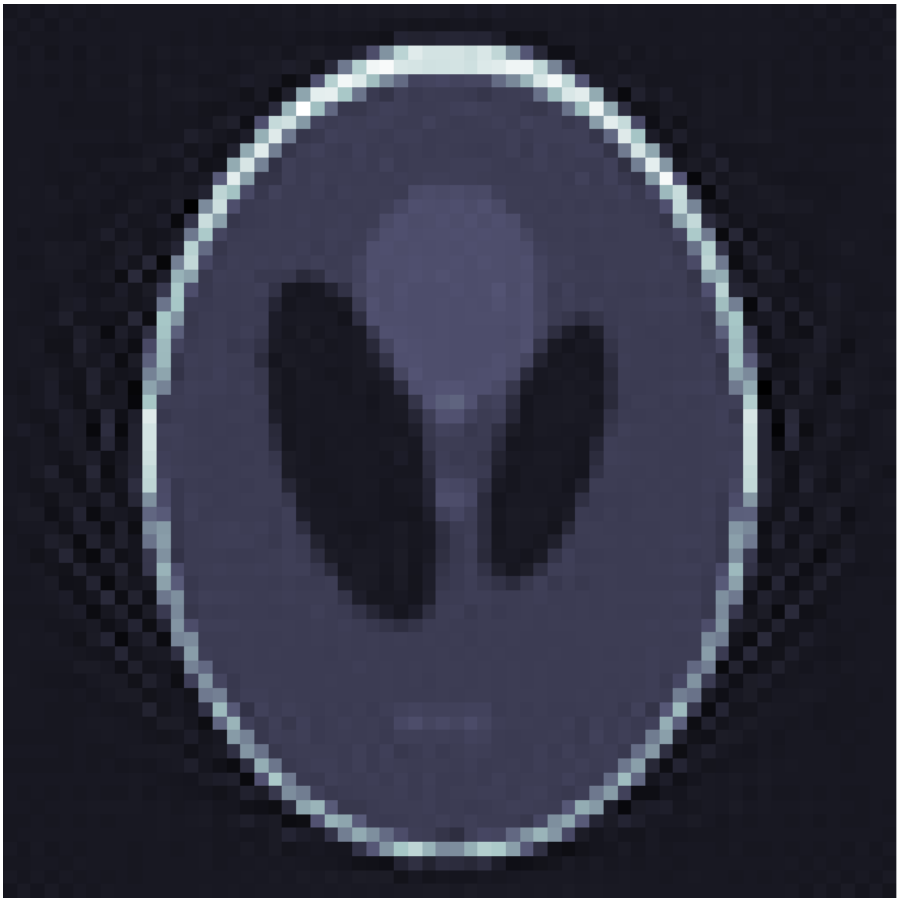}&
\includegraphics[scale=0.5]{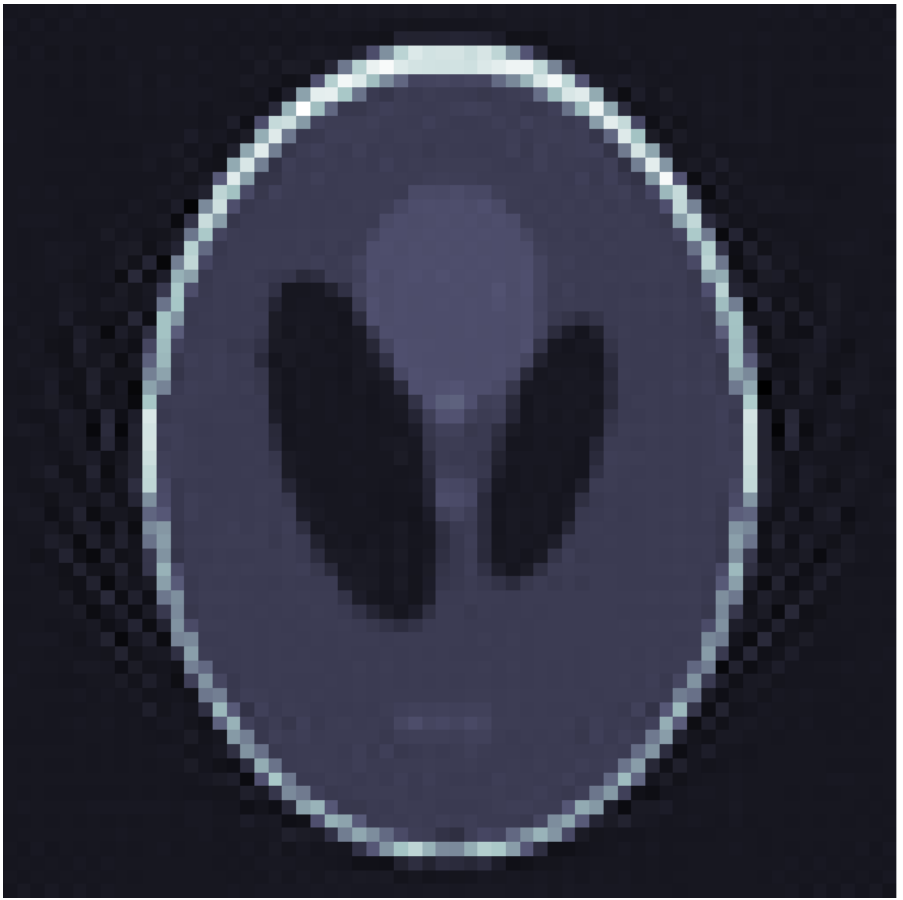} \\
FTG with $\alpha=0.5$&FTG with $\alpha=0.8$ &FTG with $\alpha=0.9$ \\ \hline
\includegraphics[scale=0.5]{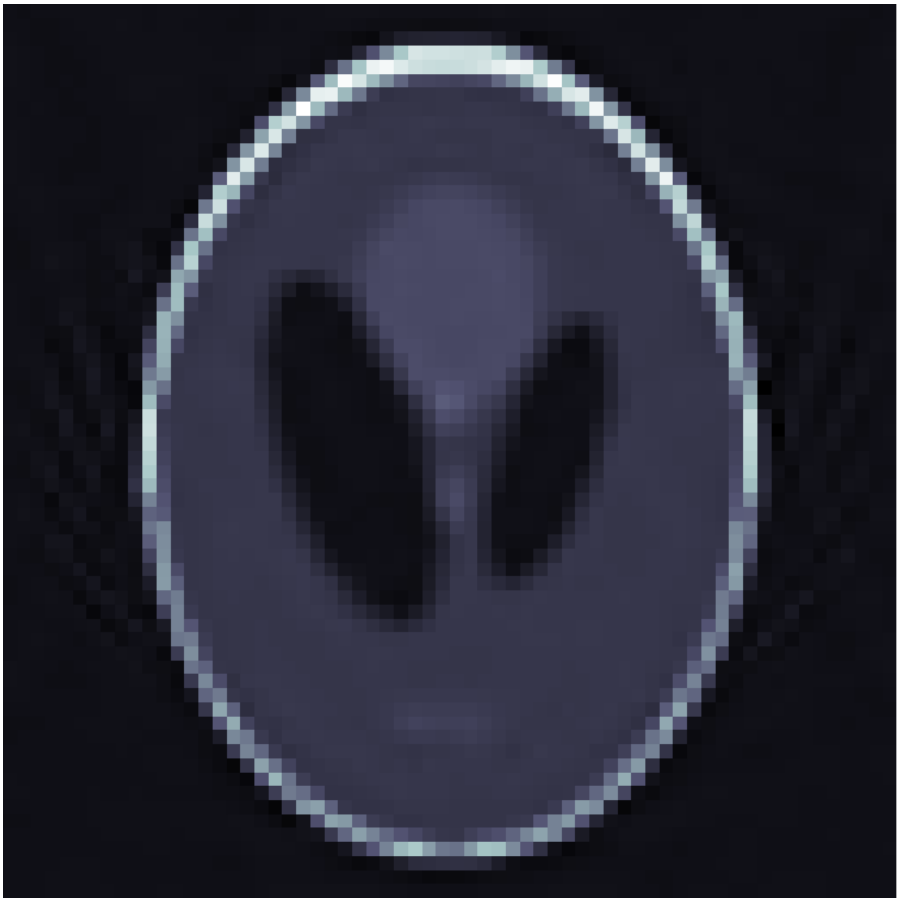}&
\includegraphics[scale=0.5]{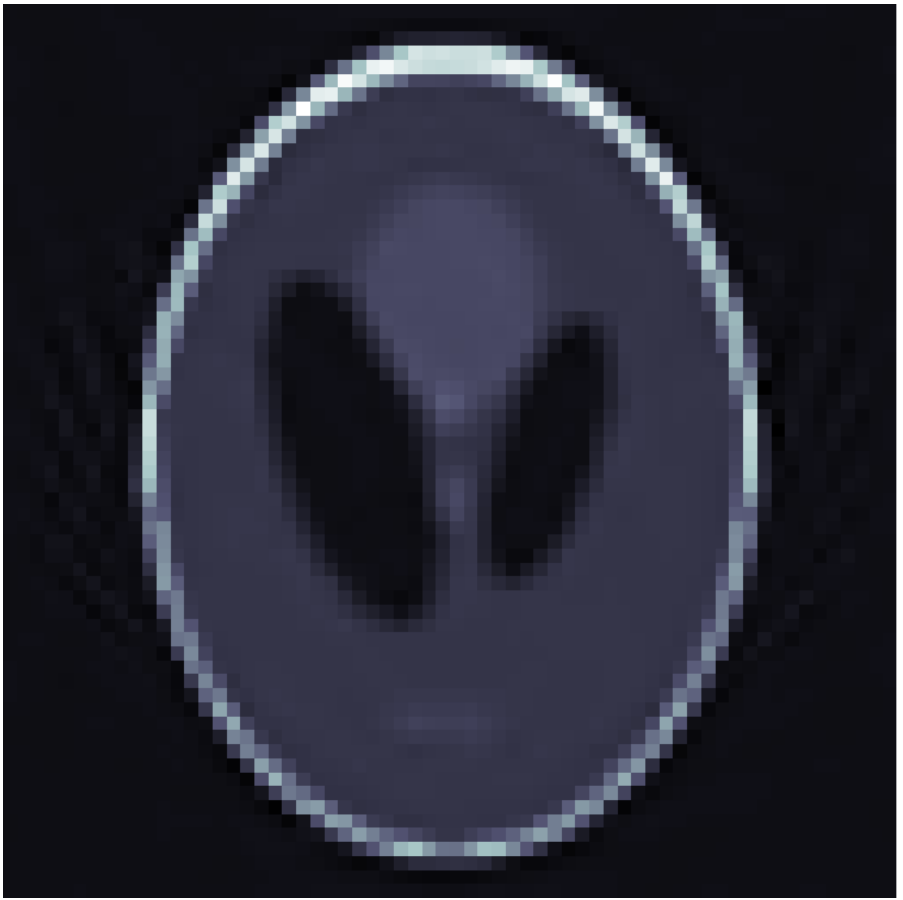}&
\includegraphics[scale=0.5]{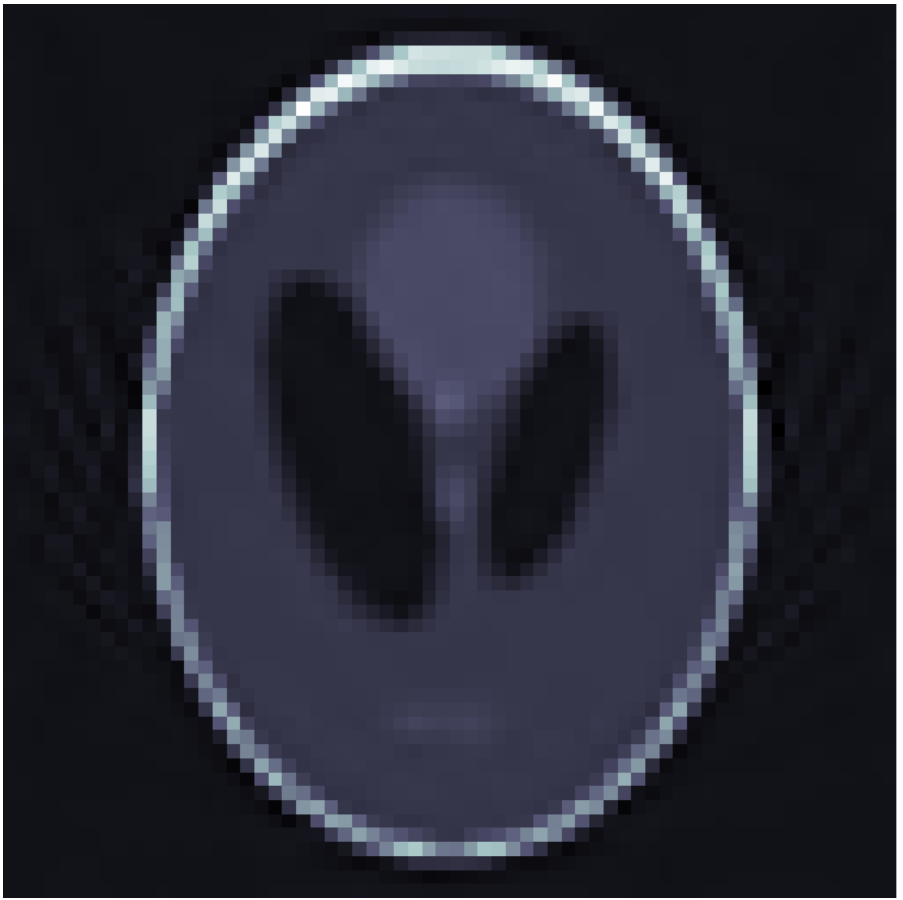} \\
FTG with $\alpha=1.1$ &FTG with $\alpha=1.2$ &FTG with $\alpha=1.5$ \\  \hline
\end{tabular}
\caption{Reconstruction results for the phantom of Shepp-Logan.
Target (top left), the Shepp-Logan phantom, $64\times 64$ pixel, and the red lines in the target image indicate the areas used in the line plots.
The corresponding FBP reconstruction (top middle), the reconstruction results with TG prior (top right) and the FTG prior (middle and bottom row) for the fractional order $\alpha=0.5,\; 0.8,\; 0.9,\;1.1,\; 1.2,\; 1.5$.
}
\label{2d_result}
\end{figure}

\begin{figure}[htbp]
 \centering
 \begin{tabular}{@{}c@{}c@{}c@{}}
 \multicolumn{3}{@{}c@{}}{\includegraphics[scale=0.4]{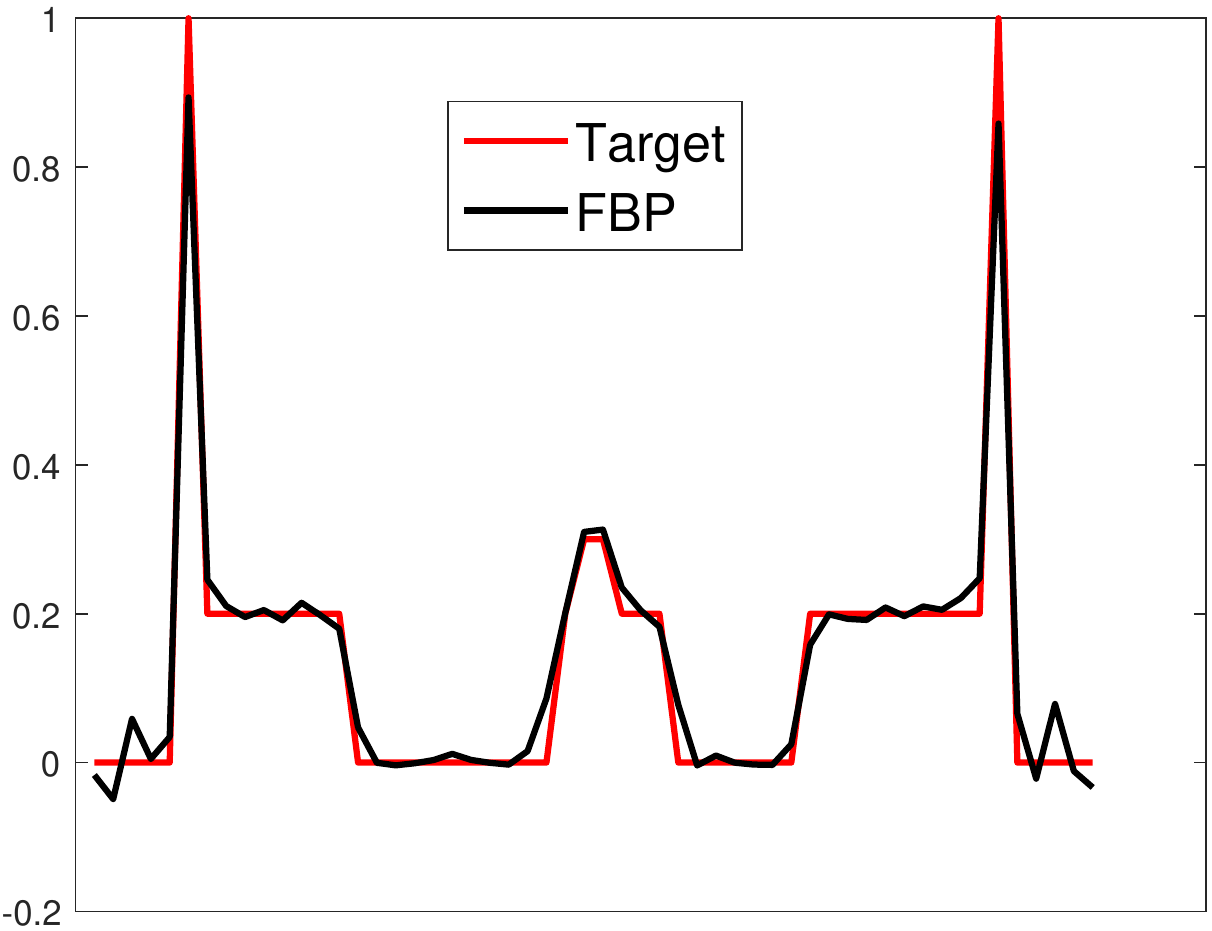}
 \includegraphics[scale=0.4]{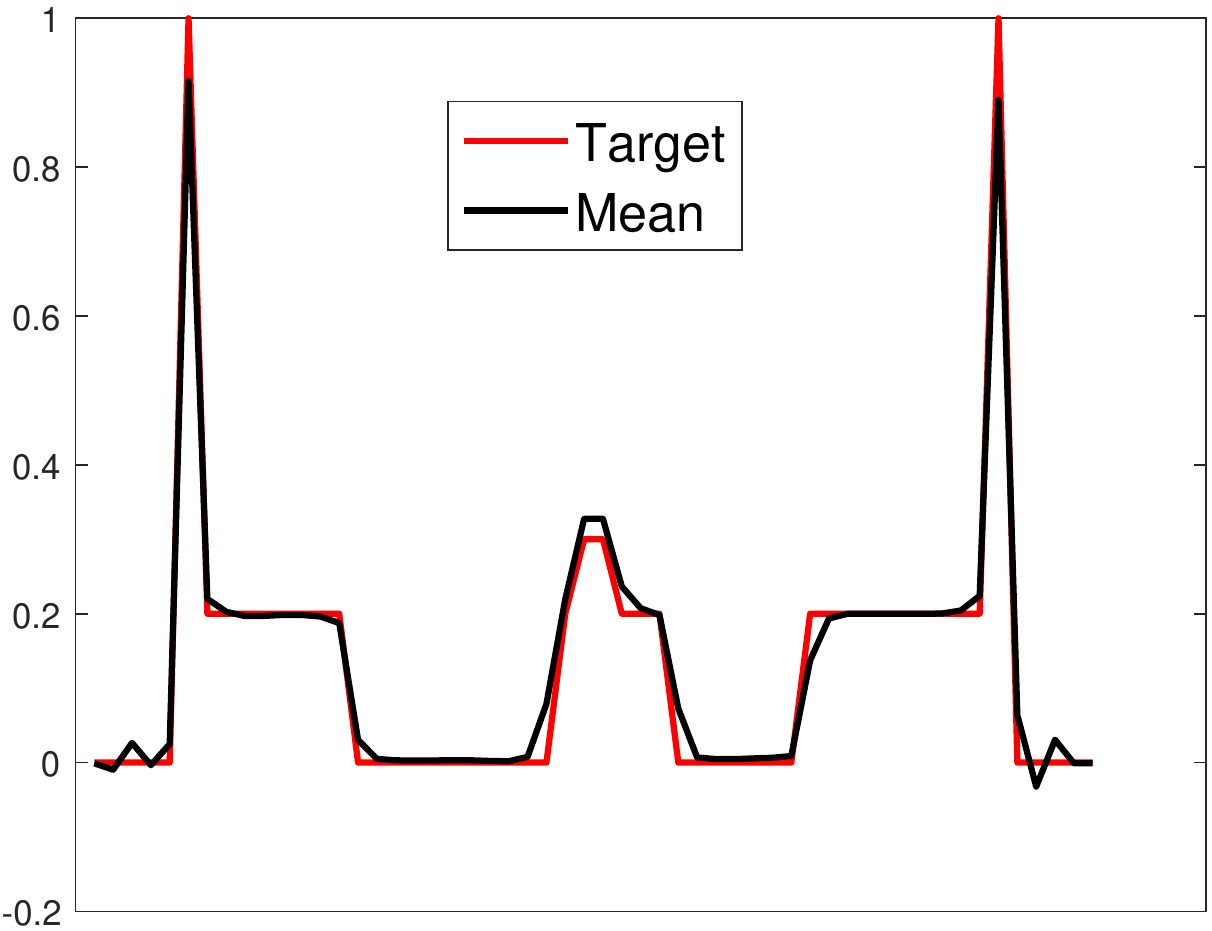}}\\
 \multicolumn{3}{@{}c@{}}{FBP \hspace{110pt} TG} \vspace{10pt}\\ \hline
\includegraphics[scale=0.4]{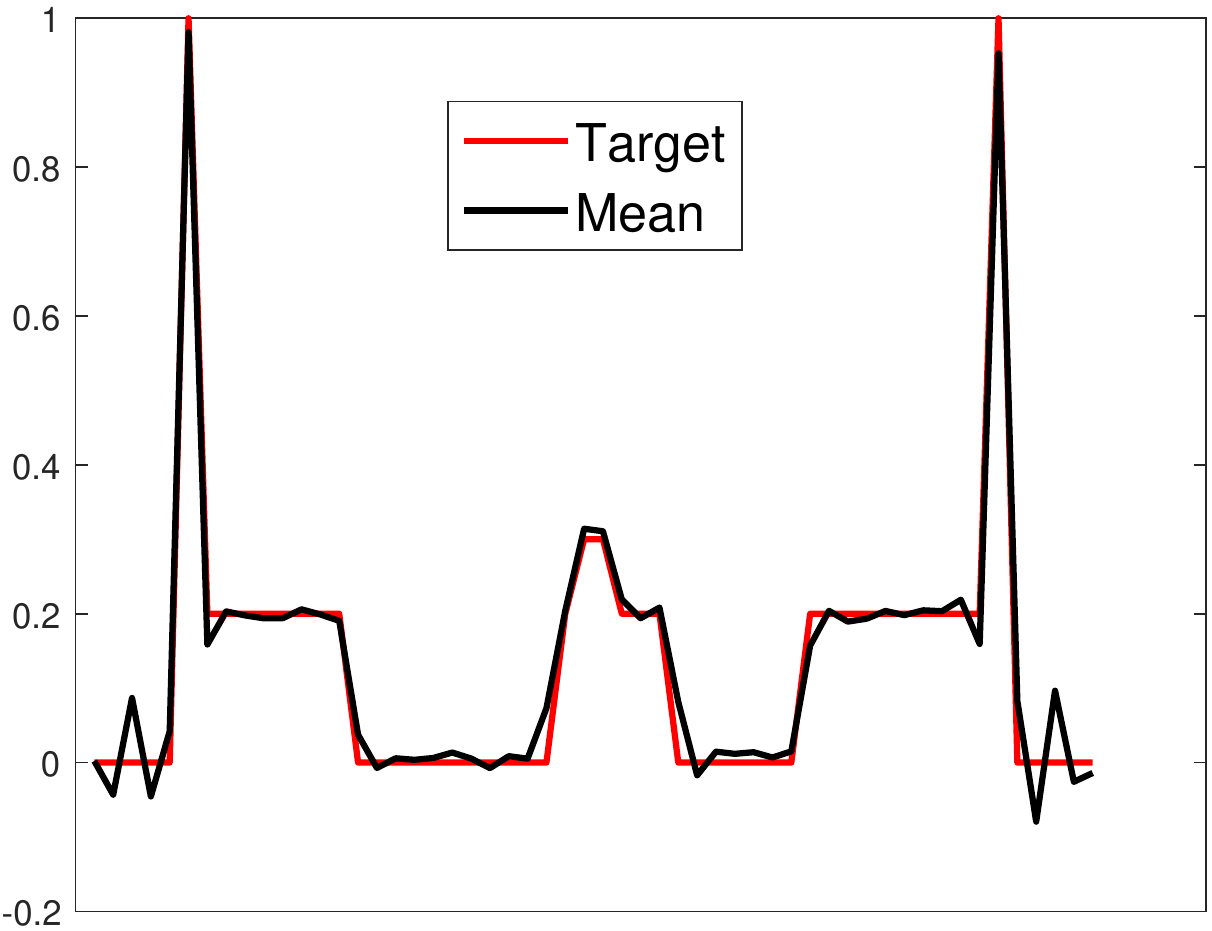}&
\includegraphics[scale=0.4]{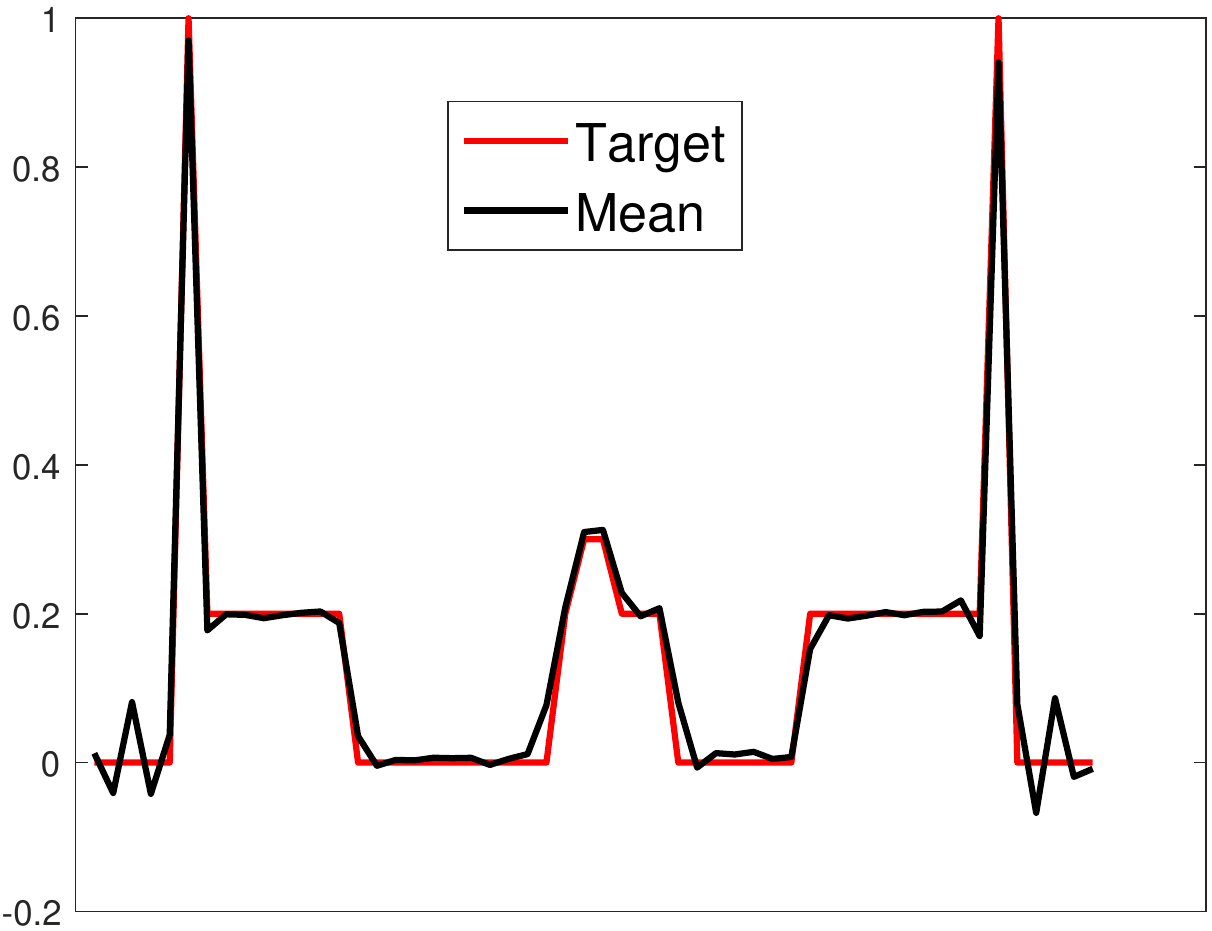}&
\includegraphics[scale=0.4]{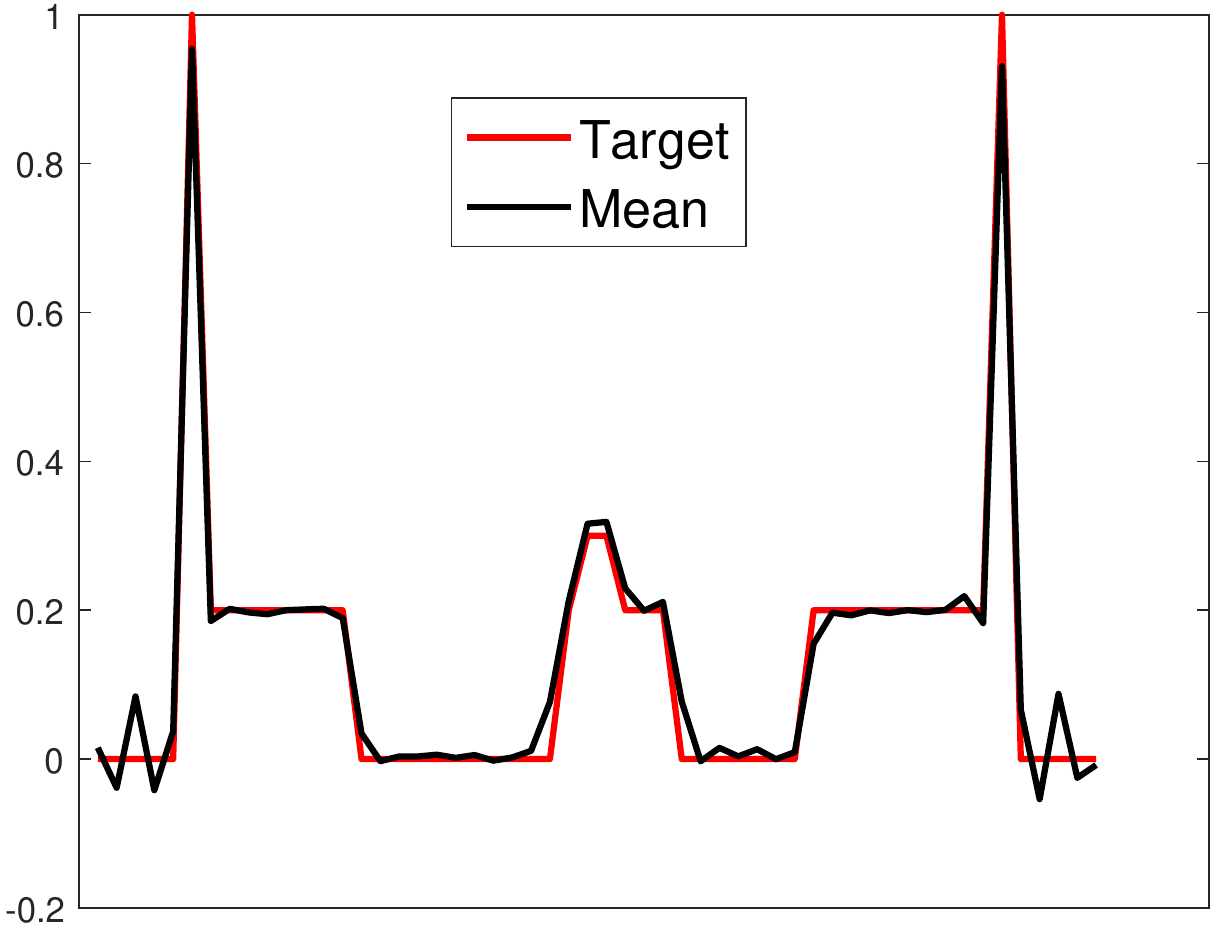} \\
FTG with $\alpha=0.5$&FTG with $\alpha=0.8$ &FTG with $\alpha=0.9$  \\ \hline
\includegraphics[scale=0.4]{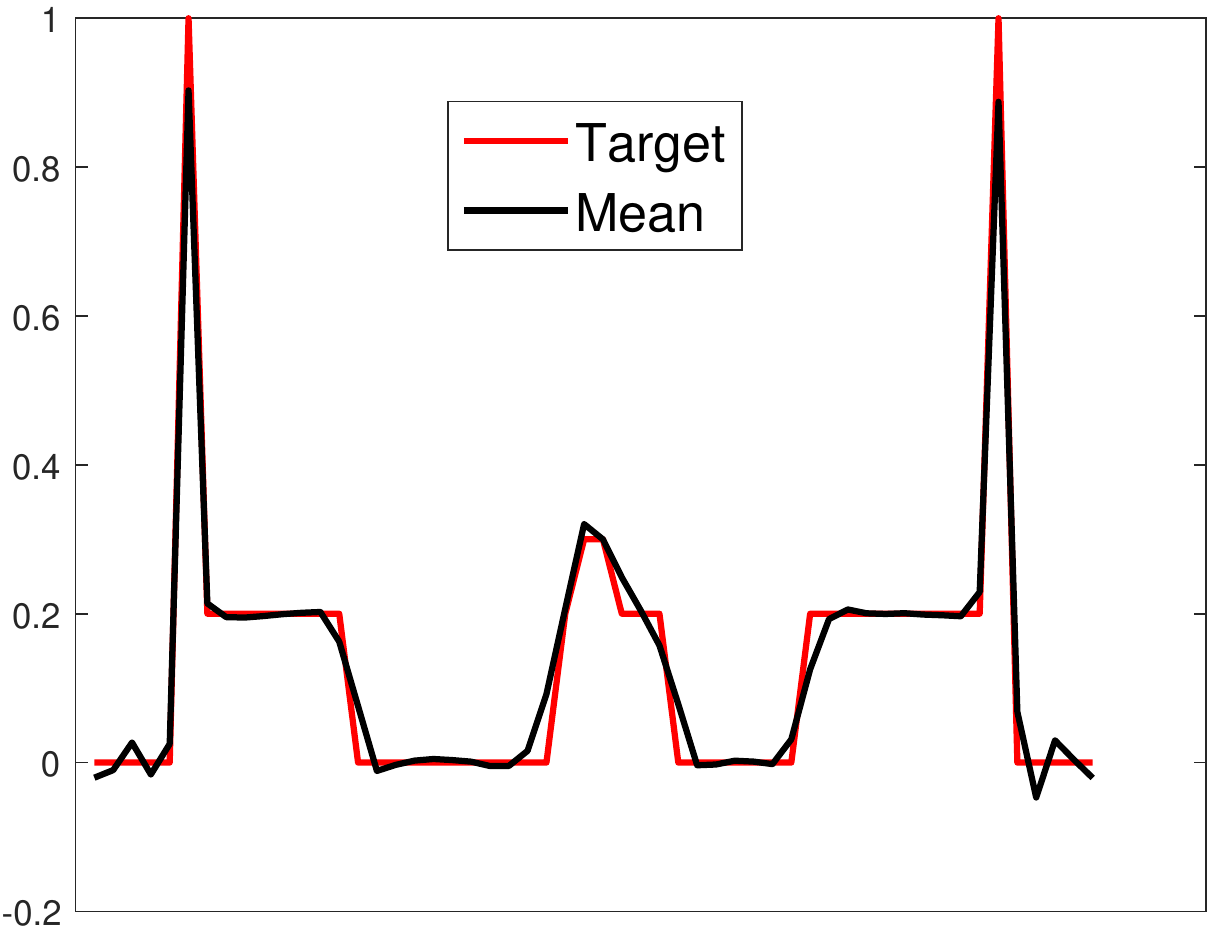}&
\includegraphics[scale=0.4]{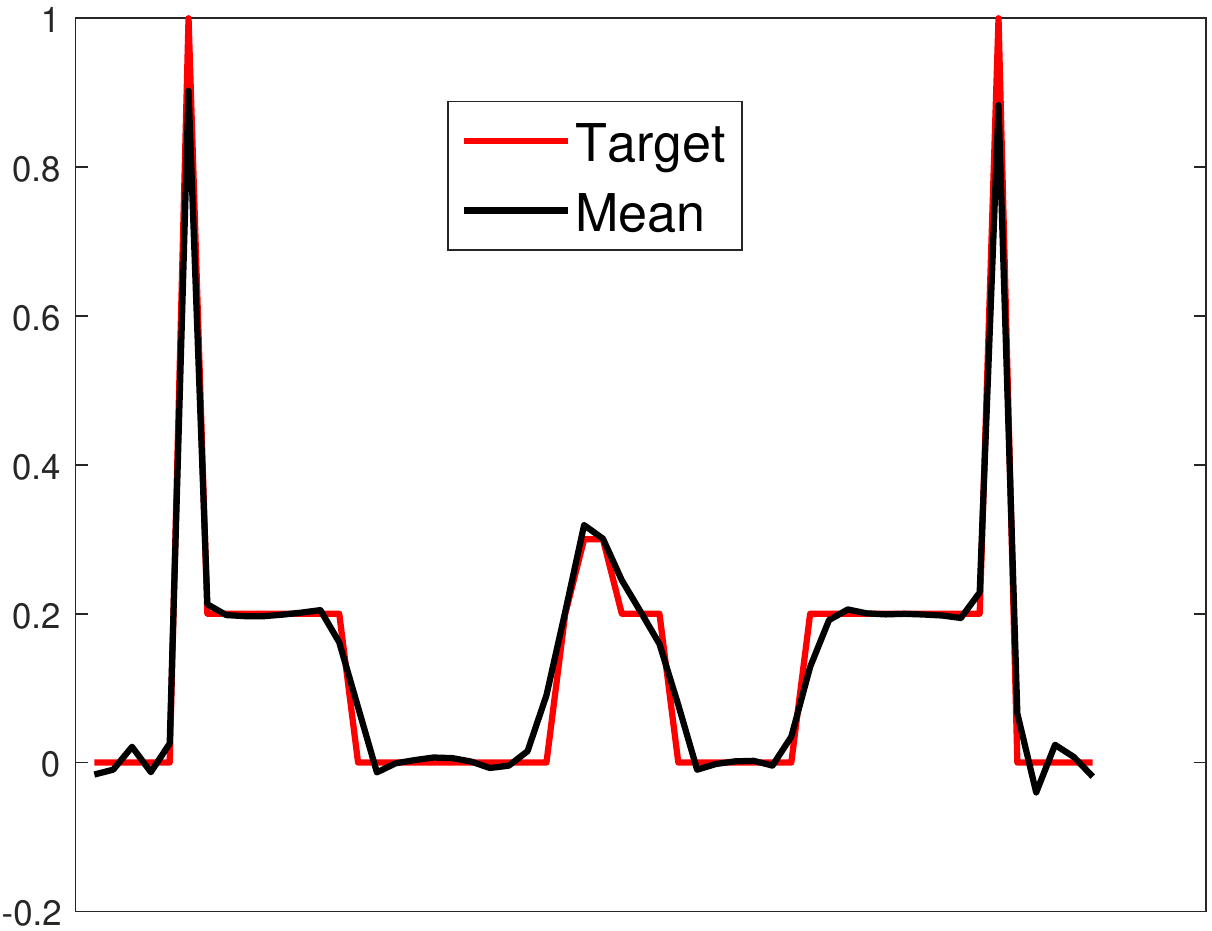}&
\includegraphics[scale=0.4]{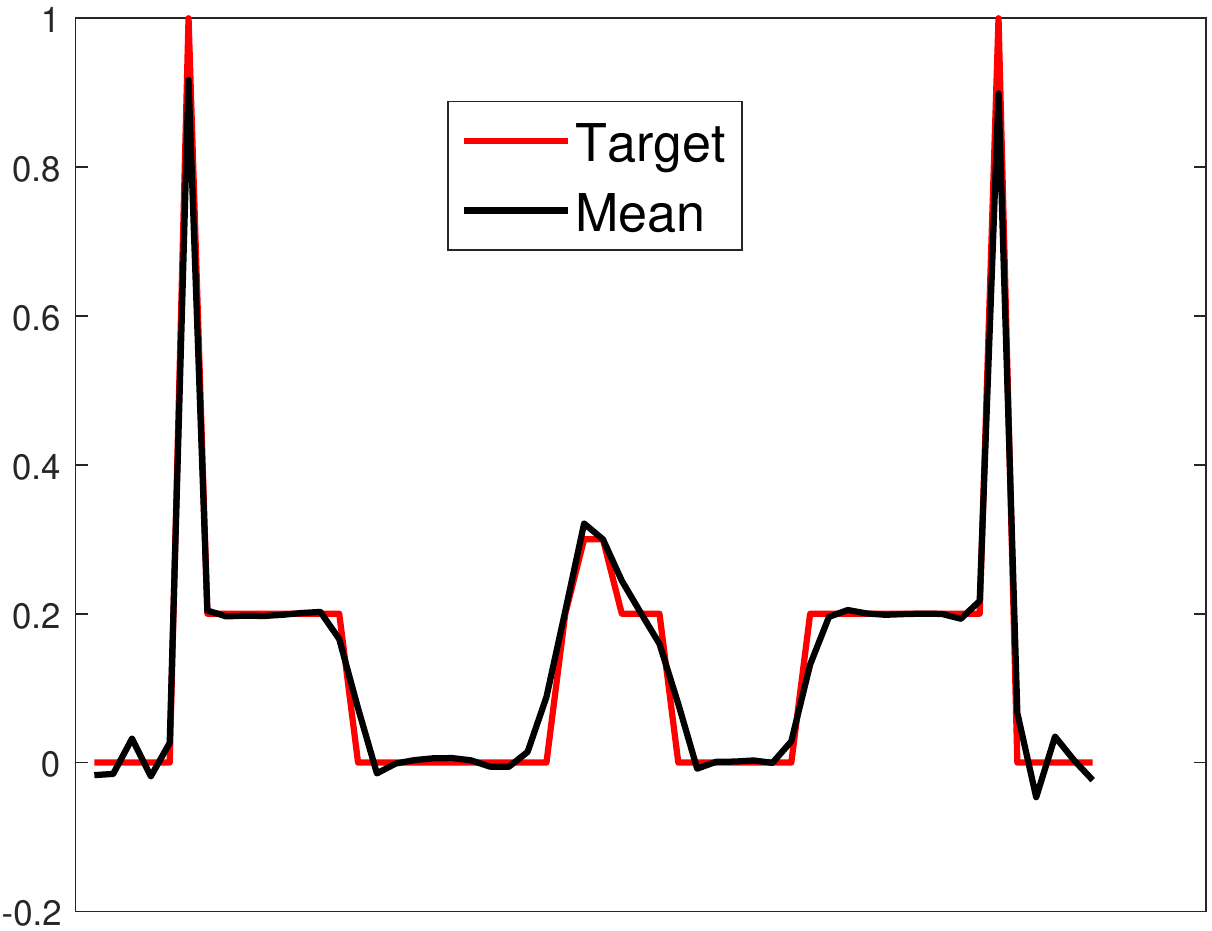} \\
FTG with $\alpha=1.1$ &FTG with $\alpha=1.2$ &FTG with $\alpha=1.5$ \\ \hline
 \end{tabular}
\caption{Limited computed tomography reconstruction. Reconstruction results for the red lines areas in the target image (see, top left in Figure\ref{2d_result}).
The corresponding FBP reconstruction (top left), the reconstruction results with TG prior (top right) and the FTG prior (middle and bottom row) for $\alpha=0.5,\; 0.8,\; 0.9,\;1.1,\; 1.2,\; 1.5$.
}
\label{1d_result}
\end{figure}

Top left in Figure \ref{2d_result} shows the target Shepp-Logan phantom $64\times 64$ pixel and the area on the image where the red line passes are used to make a line plot in Figure \ref{1d_result}.
The  top middle image in Figure \ref{2d_result} is the FBP reconstructed result using the \textsl{iradon.m} in the MATLAB Image Processing Toolbox, and the top right and middle row as well bottom in Figure \ref{2d_result} are our reconstructed images with the TG and FTG prior, respectively.
Figure \ref{1d_result} shows the reconstruction images for the lines plot in term of the red line areas in target image (see top left in Figure \ref{2d_result}).

It can be seen that the reconstructed images with FTG prior for $\alpha =1.1,\; 1.2,\; 1.5$, in Figure \ref{2d_result} and the $RelErr$ and SSIM values in Table \ref{table_ct}, are consistent with the TG prior and outperforms results through the FBP reconstruction method.
And the FTG prior with $\alpha =1.1,\; 1.2,\; 1.5$ can eliminate the blocky effect and be able to well remove under-sampling artifacts while preserving high resolution information, cf. the bottom row in Figure \ref{1d_result} and Figure \ref{2d_result}.
From the Table \ref{table_ct}, the FTG prior with $\alpha =0.5,\; 0.8,\; 0.9$ can yield the lower $RelErr$ value compared with TG and the FBP.
Because the phantom of Shepp-Logan has much details information, and for the FTG prior with $\alpha \to1$ the reconstruction results can better preserve this small details such as textural information, cf. the middle row in Figure \ref{1d_result} and Figure \ref{2d_result}.

Generally speaking, the Shepp-Logan phantom image is of piecewise constant gray level, which is the case where the TG prior is most applicable.
However, our FTG prior with $\alpha \to1$ still has a good performance for the reconstruction result.
\begin{table}[htbp]
\centering
\caption{Limited computed tomography reconstruction. Error measures of reconstruction results using FBP and FTG with $\alpha=0.5,\; 0.8,\; 0.9,\; 1.1,\; 1.2,\; 1.5$ and TG prior.}
\begin{tabular}{ccccccccc}
  \toprule
          &FBP & TG & $\alpha=0.5$ & $\alpha=0.8$ & $\alpha=0.9$ & $\alpha=1.1$ & $\alpha=1.2$ & $\alpha=1.5$ \\
  \midrule
  $RelErr$ &0.4031& 0.3827 & \textbf{0.3805} & \textbf{0.3759} & \textbf{0.3744}&0.3883&0.3900 & 0.3858 \\
  SSIM    &0.7480& 0.8670 & 0.7976& 0.8433 &0.8531& 0.8484& 0.8488 & 0.8501 \\
  \bottomrule
  \end{tabular}
  \label{table_ct}
\end{table}
\subsection{Image denoising}\label{2d_ex2}
Image denoising is a simple yet heavily addressed problem in image processing.
The noise in digital images can be caused by the failure or poor performance of image sensors, or failure of the data transmission process.
The denoising process deals with the removal of noise from the noisy images considered.

\subsubsection{Problem setup}
The denoising problem is often modeled by a degradation model, which, in this paper, is given by
\begin{align}
\label{de_model}
\mathbf{y}=f(\mathbf{x})+\boldsymbol{\eta},\quad \mathbf{x}\in\Omega\subset\mathbb{R}^2,
\end{align}
where the $f$ represents the original image, $\mathbf{y}$ stands for a degraded/noisy image (given data) and the  random variable $\boldsymbol{\eta}$ is assumed a Gaussian white noise with the known standard deviation.
For a given degraded image $\mathbf{y}$, the image denoising is to reconstruct the original noise-free image $f$.
In fact, the degradation model \eqref{de_model} is equivalent to the model \eqref{Model} for the $\mathcal{A}$ as an identity operator.

\subsubsection{Set up the inverse problrm}
In this example, the original image$-$cameraman $f(\mathbf{x})$, defined on $\Omega=[-1,1]^2$, is choosing $128\times128$ pixel, i.e. the dimension $d=1.6384\times 10^4$.
The noisy image $\mathbf{y}$ is taken as the original image corrupted by Gaussian noise with zero mean and standard deviation $0.03$, which corresponds to $3\%$ noise with respect to the maximum norm of the original image, seeing the top right of Figure \ref{2d_result_de}.
The Gaussian prior $\mu_0$ is taken to be the standard Gaussian distribution and the FTG prior with the fractional order $\alpha=0.5,\; 0.9,\; 1.1,\;1.5,\; 1.8$.
The shape and rate parameter about Gamma distribution are set to $k=3\times10^4,\;\vartheta=1$, respectively.

Then, using the alternating direction method we construct a transport map between the standard Gaussian distribution $\mu_0$ and the posterior \eqref{H_postprior}.
The numerical optimization problem \eqref{d_opt_problem}, analogy to above CT reconstruction, is performed with MATLAB's fmincon optimizer, where the the step tolerance (StepTolerance) is set to $10^{-3}$; the SpecifyConstraintGradient and SpecifyOdjectiveGradient are set to true and we use $M=2000$ samples of the standard Gaussian distribution $\mu_0$ to approximate the expected value via the SAA in the objective function (see Section \ref{numercal opt}).
The transport map $\widetilde{T}$ and the regularization parameter $\lambda$ via the formula \eqref{de_reg} are then used as precondition for the linear diagonal map-based independence sampler.
And the proposal $\mu_{ref}$ is taken to be a Gaussian distribution with zero mean and the standard deviation $0.01$ in the sampling process.
The \textsl{Peak signal-to-noise ratio} (PSNR), which does well in measuring quality of the reconstructed images in terms of image denoising, is also adopted in this case.

\subsubsection{Result}

\begin{figure}[htbp]
 \centering
 \begin{tabular}{@{}|@{}c@{}|@{}c@{}|@{}c@{}|}
 \multicolumn{3}{@{}c@{}}{\includegraphics[scale=0.55]{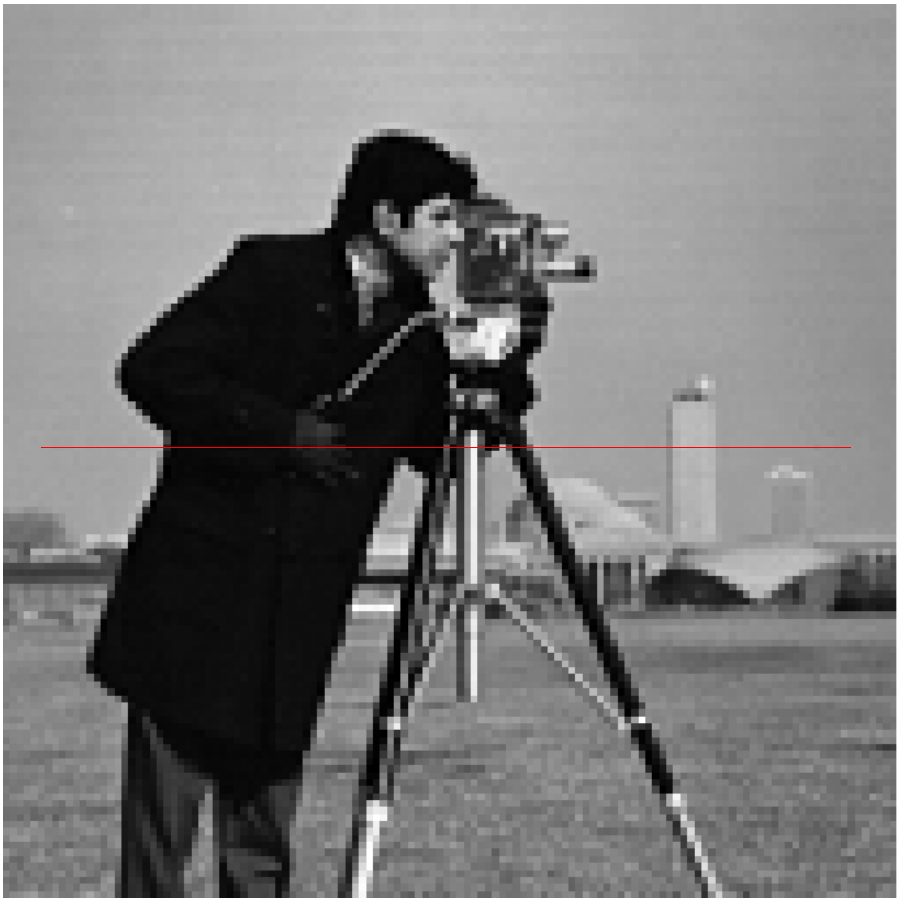}
 \includegraphics[scale=0.55]{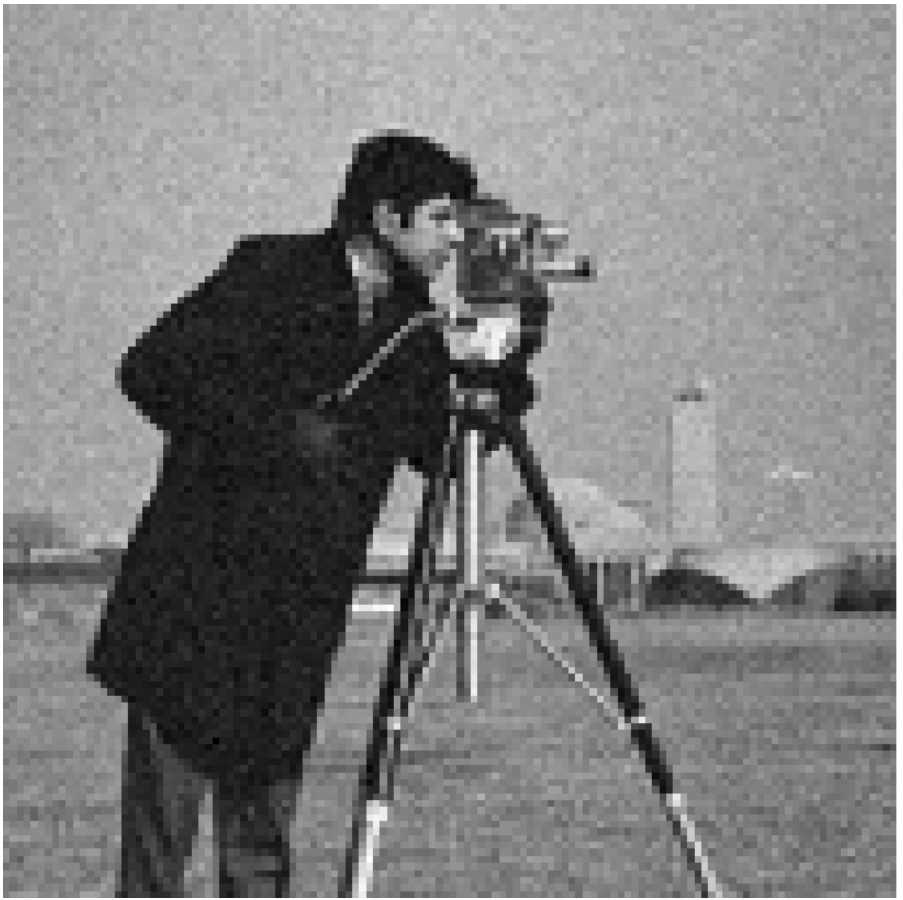}}\\
 \multicolumn{3}{@{}c@{}}{original image \hspace{70pt} noisy image} \vspace{20pt}\\
  \hline
\includegraphics[scale=0.55]{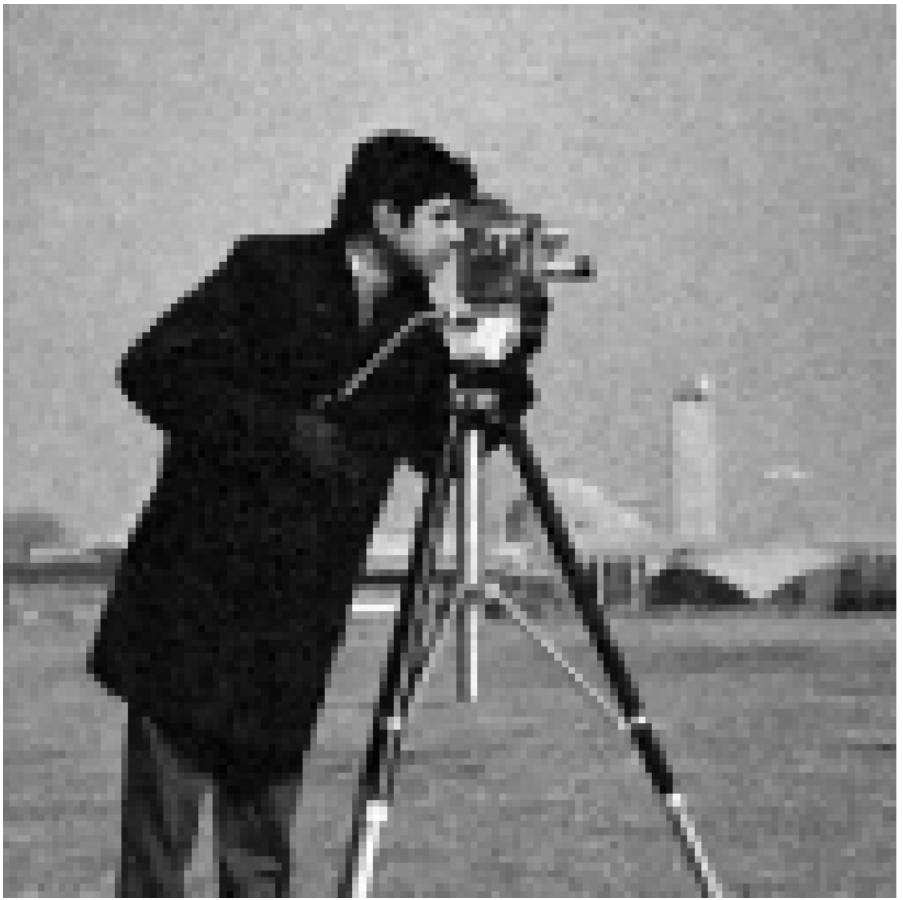}&
\includegraphics[scale=0.55]{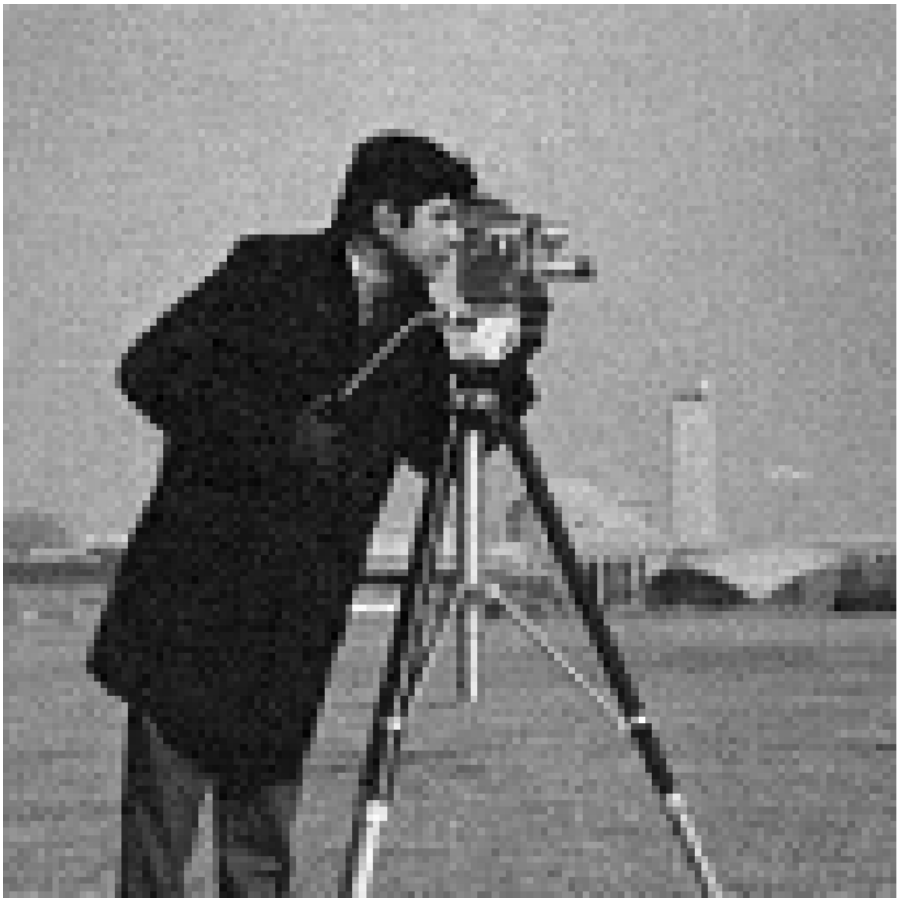}&
\includegraphics[scale=0.55]{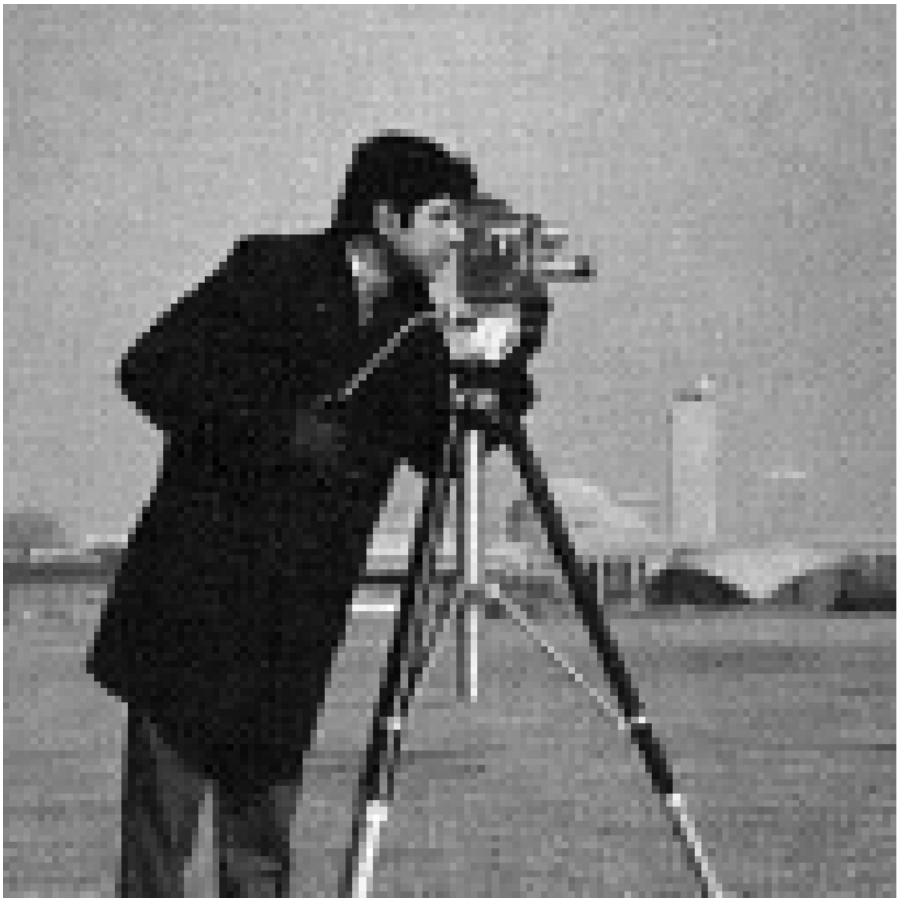} \\
TG&FTG with $\alpha=0.5$ &FTG with $\alpha=0.9$ \\
\includegraphics[scale=0.55]{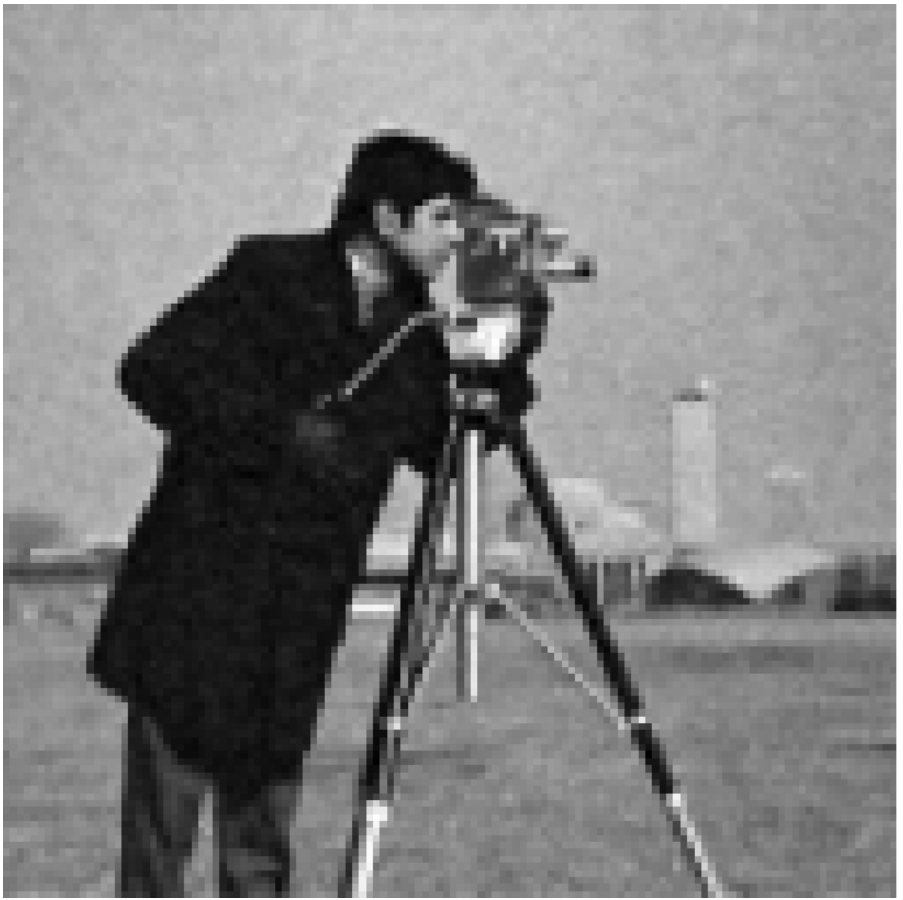}&
\includegraphics[scale=0.55]{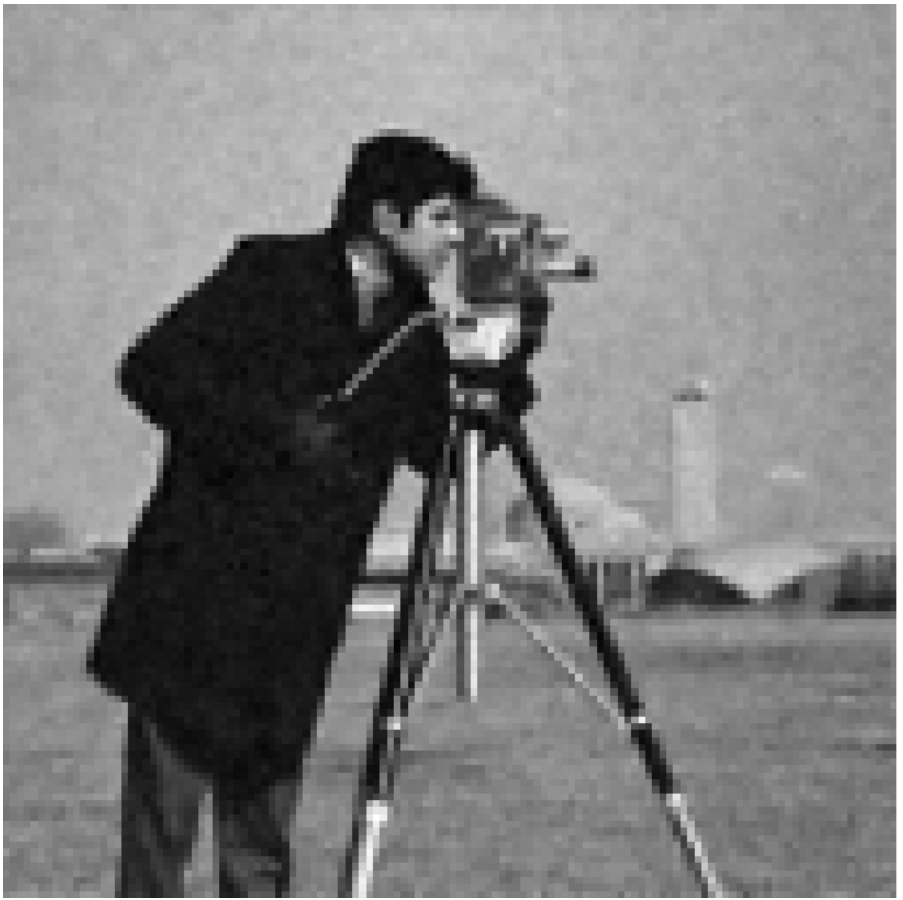}&
\includegraphics[scale=0.55]{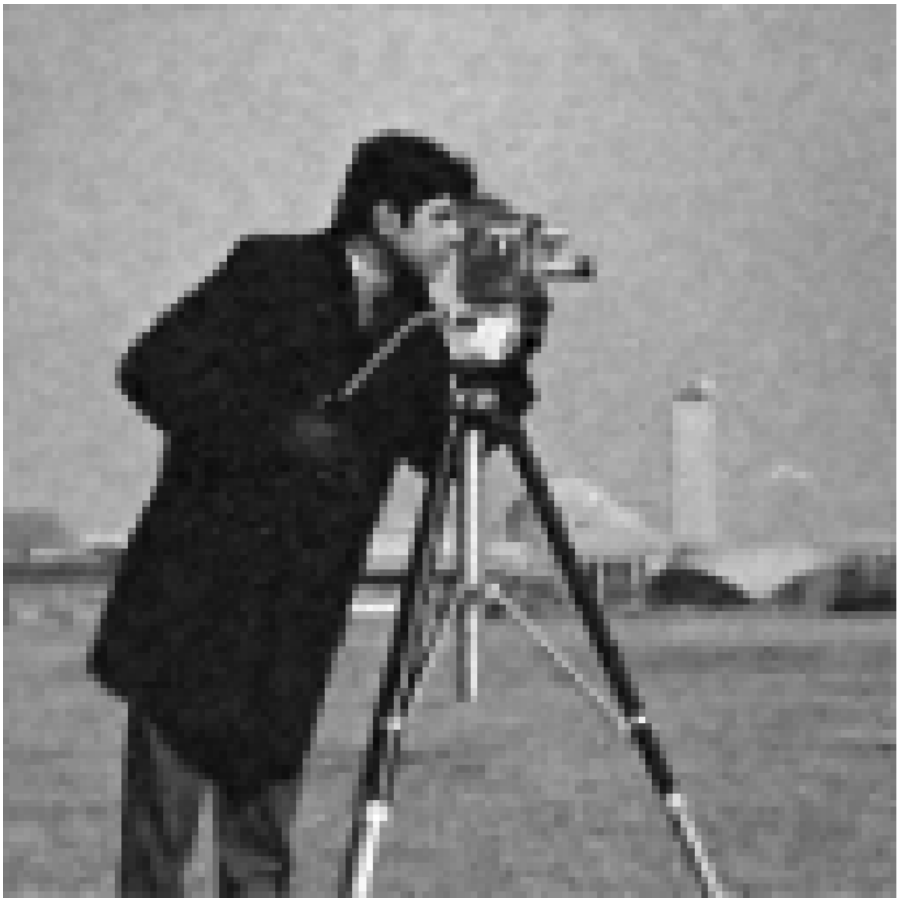} \\
FTG with $\alpha=1.1$ &FTG with $\alpha=1.5$ &FTG with $\alpha=1.8$ \\ \hline
 \end{tabular}
\caption{Reconstruction results for the image denoising.
Original image (top left), the Cameraman, $128\times 128$ pixel, and the red lines in the original image indicate the areas used in the line plots.
The noisy image (top right) is obtained by adding $3\%$ noise on the basis of the maximum norm of original image.
The reconstruction results using the FTG with the fractional order $\alpha=0.5,\; 0.9,\; 1.1,\;1.5,\; 1.8$ and TG prior (middle and bottom row).}
\label{2d_result_de}
\end{figure}

\begin{figure}[htbp]
 \centering
 \begin{tabular}{@{}c@{}c@{}}
\includegraphics[width=0.45\textwidth, height=0.2\textwidth]{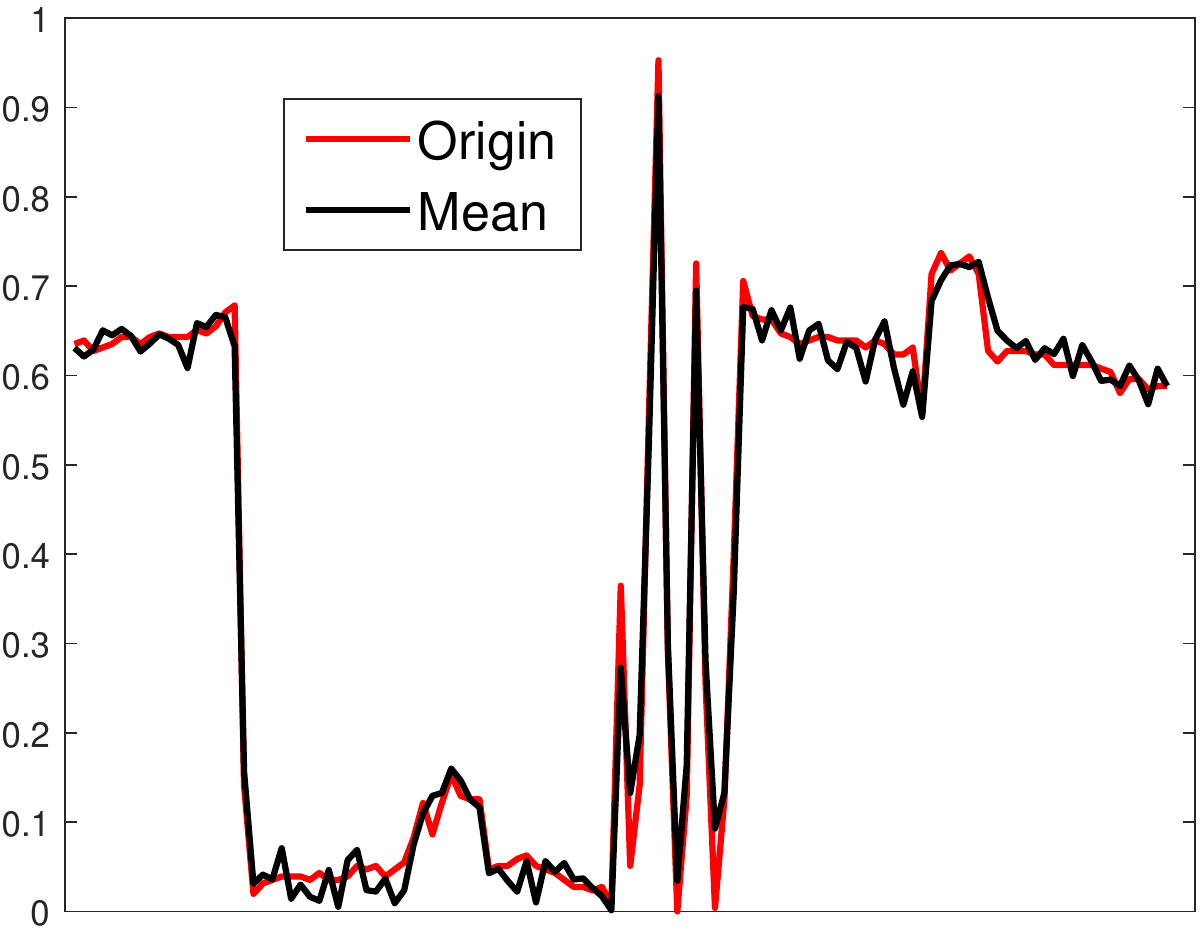}&
\includegraphics[width=0.45\textwidth, height=0.2\textwidth]{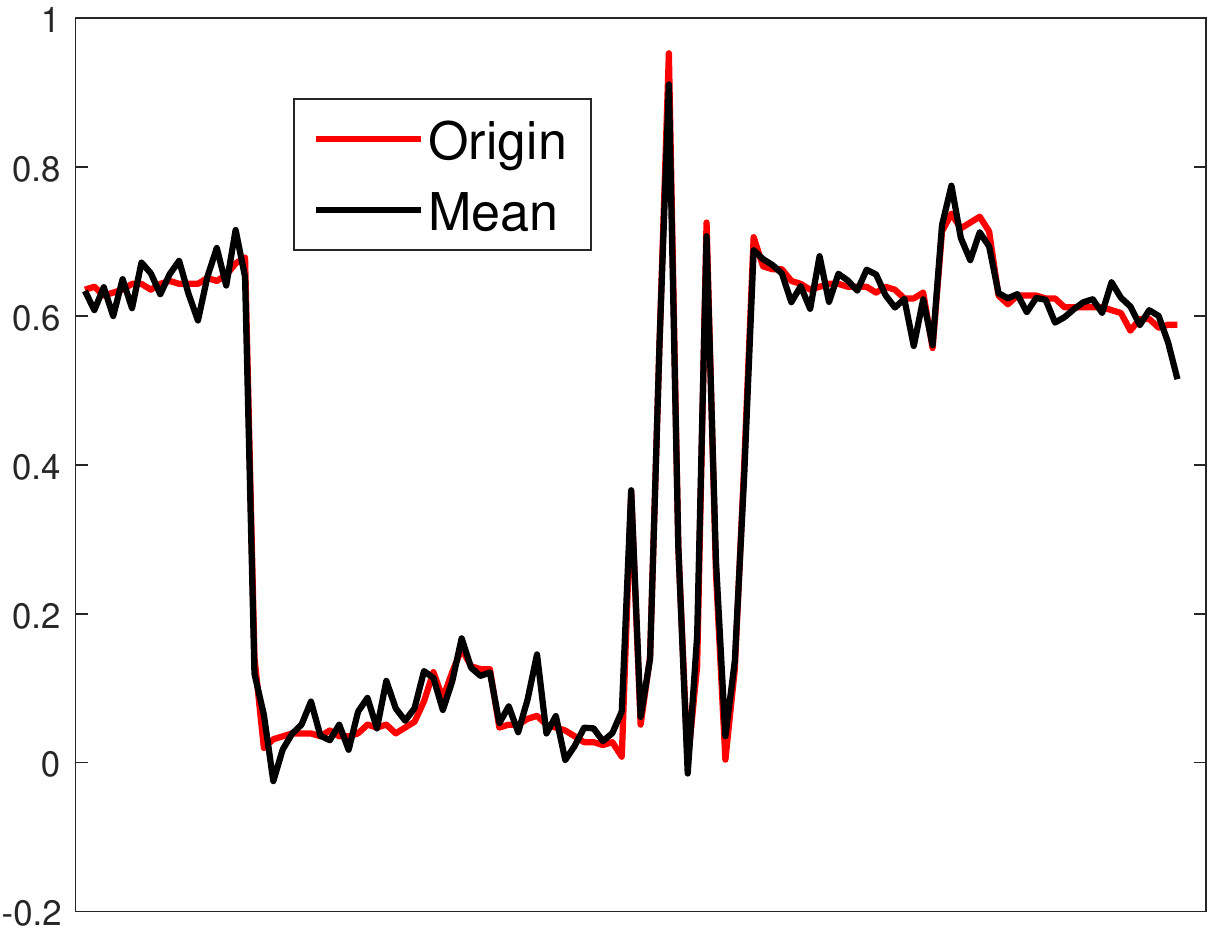}\\
TG&FTG with $\alpha=0.5$ \\  \hline
\includegraphics[width=0.45\textwidth, height=0.2\textwidth]{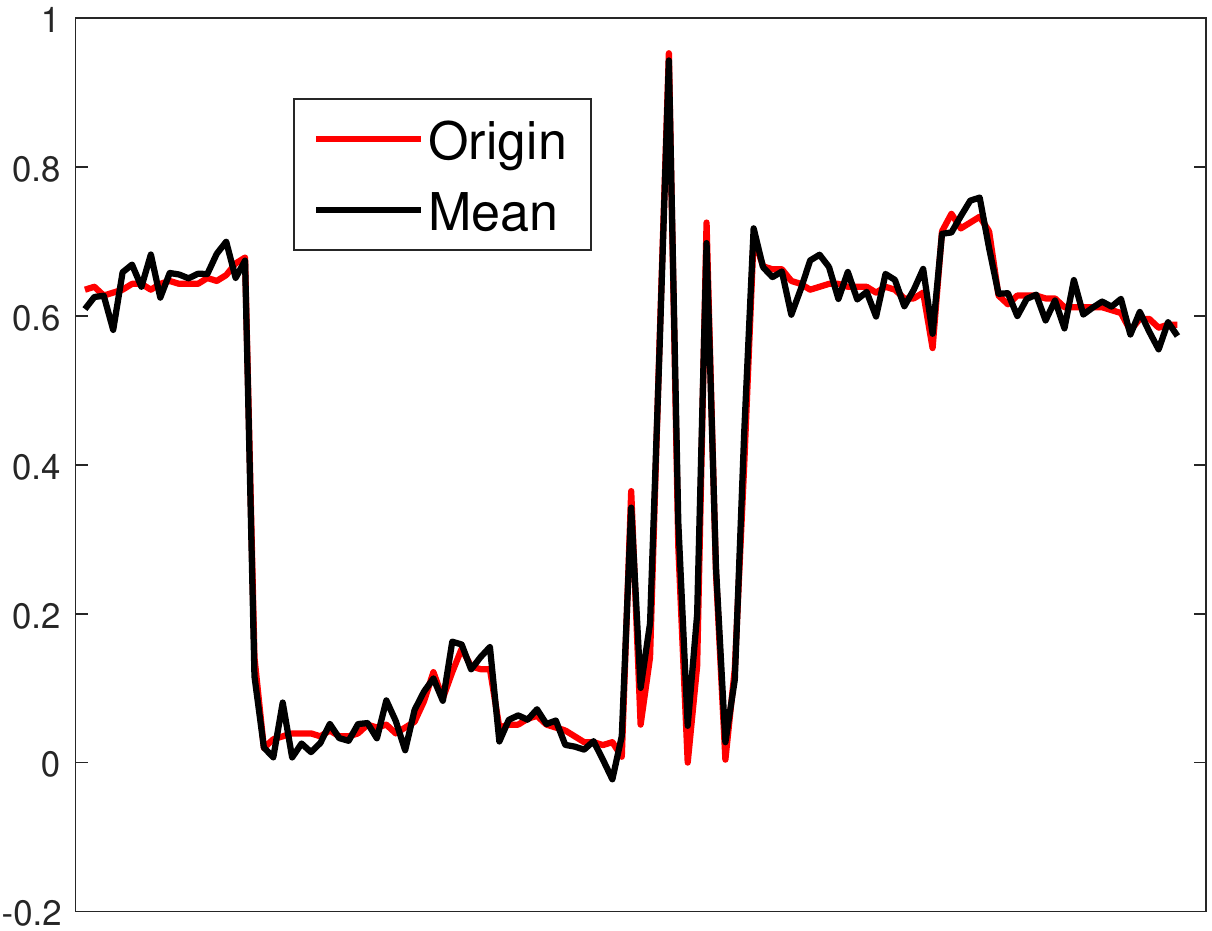} \hspace{2pt}&
\includegraphics[width=0.45\textwidth, height=0.2\textwidth]{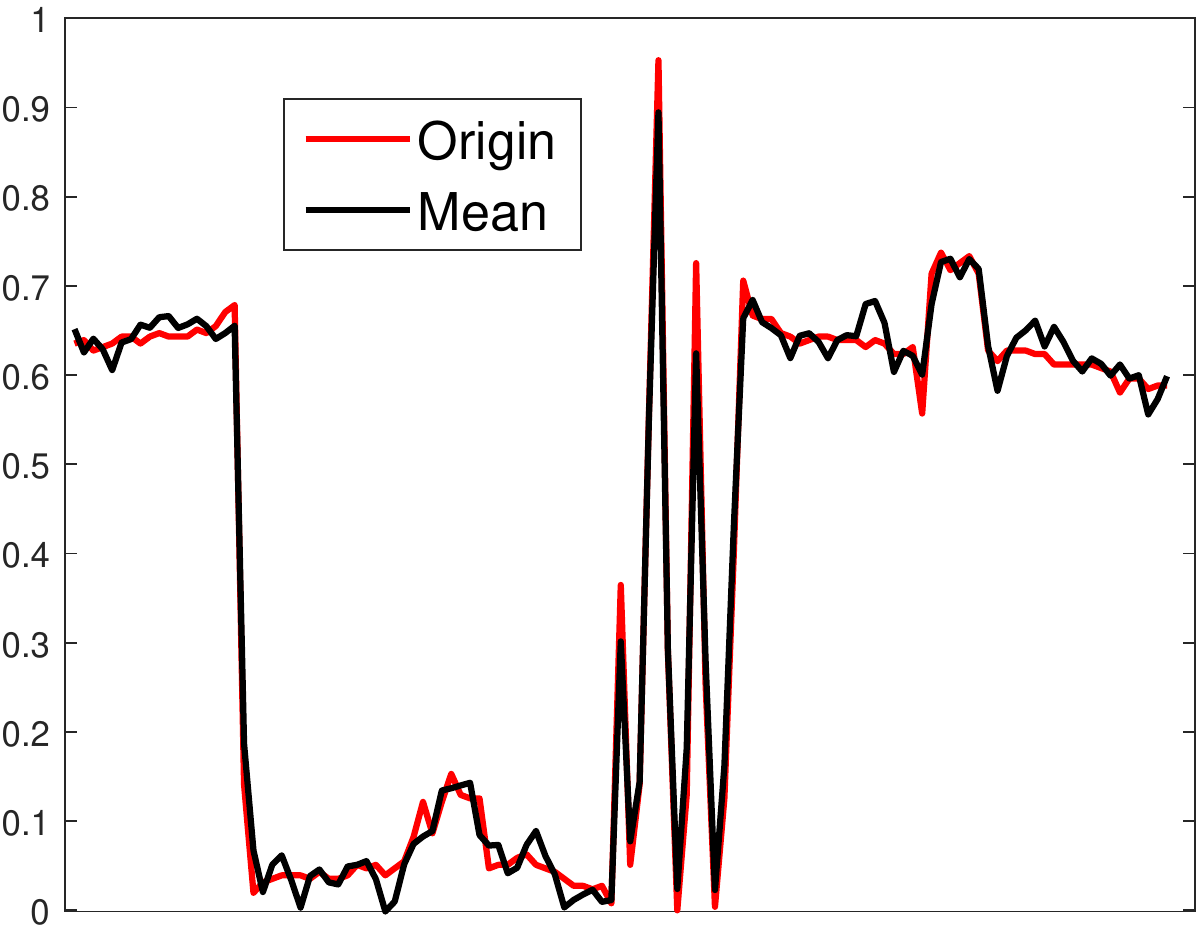}\\
 FTG with $\alpha=0.9$ &FTG with $\alpha=1.1$ \\ \hline
\includegraphics[width=0.45\textwidth, height=0.2\textwidth]{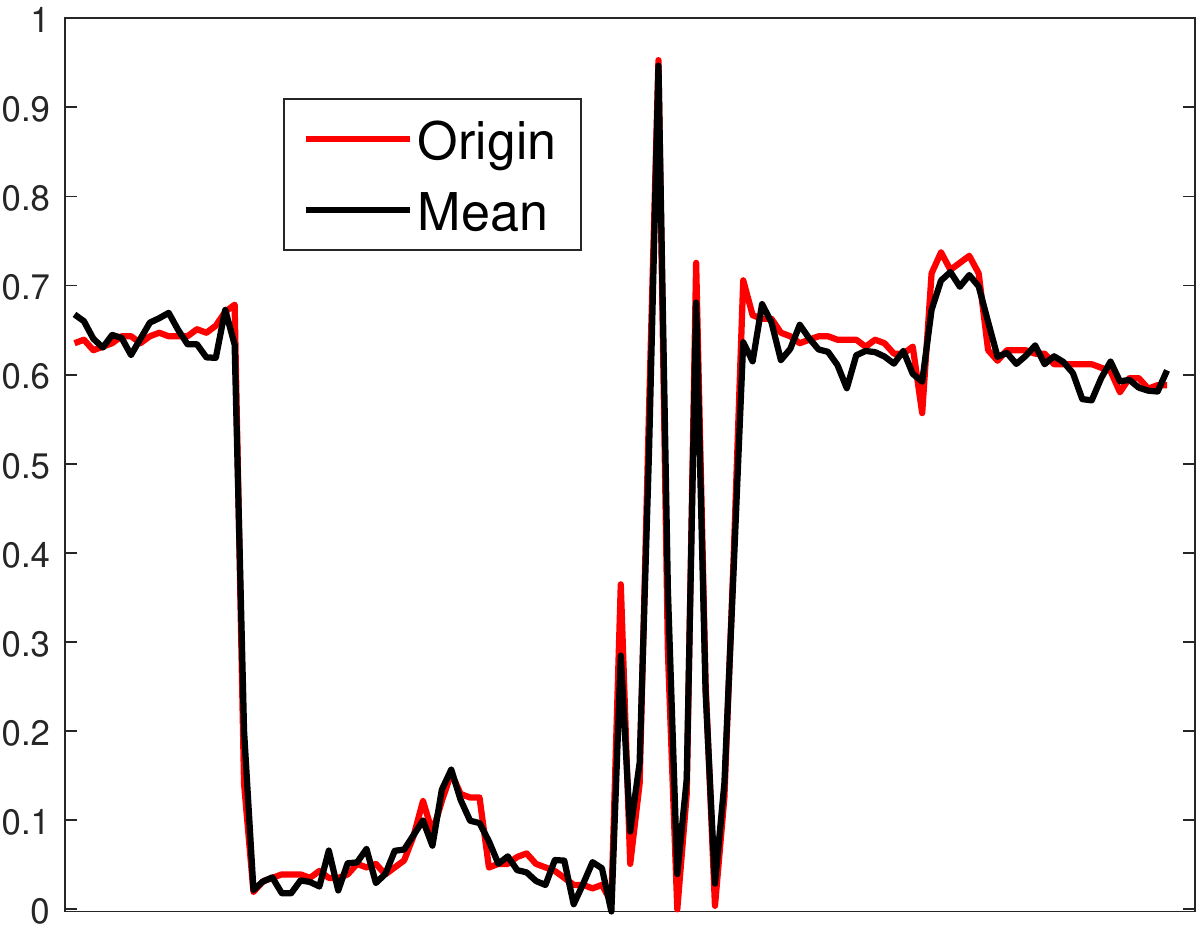}\hspace{2pt}&
\includegraphics[width=0.45\textwidth, height=0.2\textwidth]{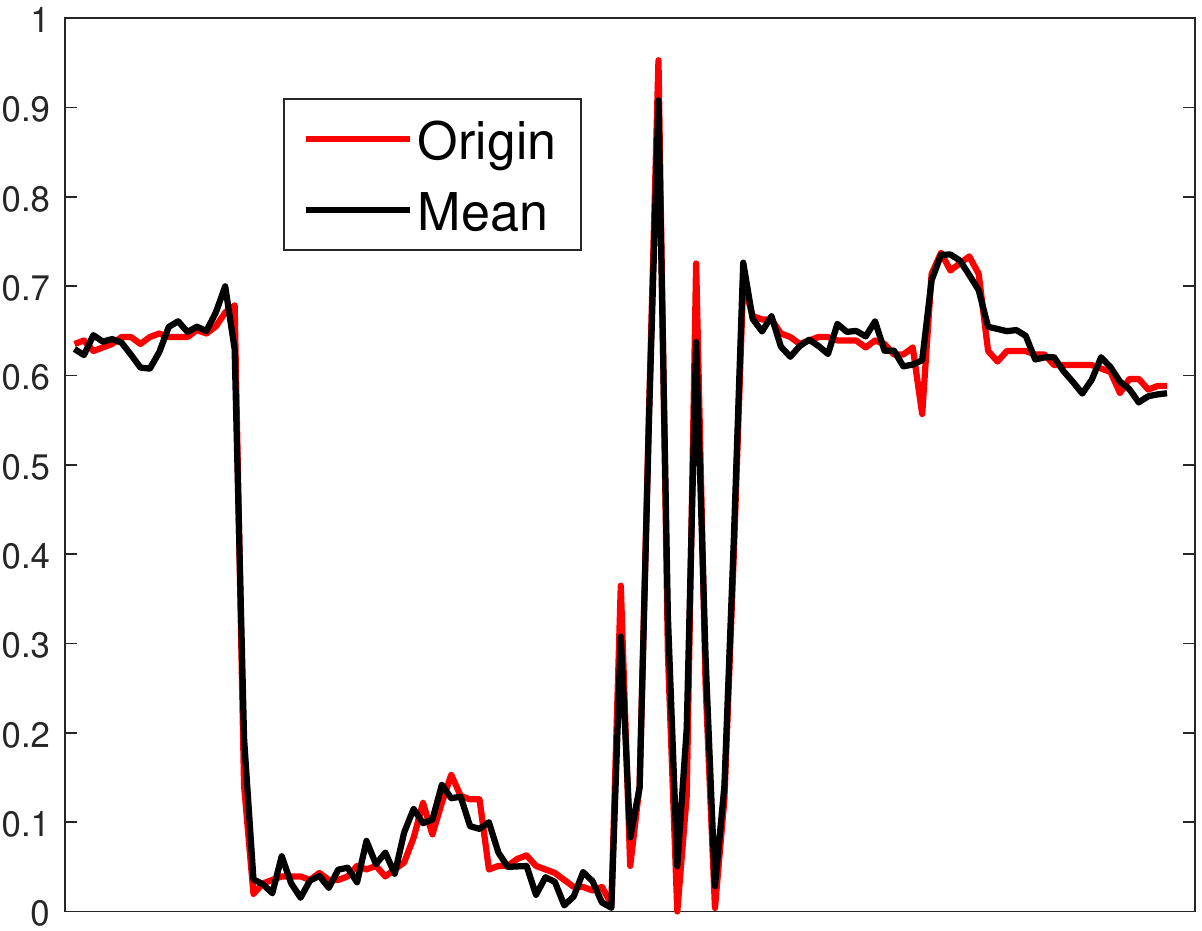}\\
FTG with $\alpha=1.5$ &FTG with $\alpha=1.8$ \\  \hline
 \end{tabular}
\caption{Image denoising. Reconstruction results for the red lines areas in the original image (see, top left in Figure \ref{2d_result_de}).
The corresponding reconstruction results using the FTG with the fractional order $\alpha=0.5,\; 0.9,\; 1.1,\;1.5,\; 1.8$ and TG prior.
}
\label{1d_result_de}
\end{figure}

Top left in Figure \ref{2d_result_de} shows the original image, the Cameraman,  $128\times 128$ pixel and the area on the image where the red line passes are used to make a line plot in Figure \ref{1d_result_de}.
The middle and bottom row in Figure \ref{2d_result} are our reconstructed images with the FTG prior for the fractional order $\alpha=0.5,\; 0.9,\; 1.1,\;1.5,\; 1.8$ and TG using $3\%$ noise with respect to the maximum norm of the original image.
Figure \ref{1d_result_de} shows the reconstruction images for the lines plot in term of the red line areas in original image (see top left in Figure \ref{2d_result_de}).

It can be seen that the reconstructed images with FTG prior in Figure \ref{2d_result} and the PSNR and SSIM values in Table \ref{table_de} are in agreement with the TG prior.
Of course, the TG prior yields a better result in term of the PSNR value from the Table \ref{table_de}.
Because the original image has some edges.
However, from Table \ref{table_de}, the FTG prior with $\alpha=1.1,\; 1.5,\; 1.8$ yield the higher SSIM value compared with the TG prior.
This is because the original image also has some texture generating the nonconstant gray level in the domain, which leads to significant blocky effect for the reconstructed result with TG prior, cf. the top left in Figure \ref{1d_result_de}, while the FTG prior can eliminate well the blocky effect caused by TG and better preserve some small details such as textural information and corner points.
Thus, for the target image with some small details, our FTG with some appropriate $\alpha$ in $(1,2]$ prior not only maintains a low relative error, which is consistent with that of TG, but also has a better reconstructed result in terms of the structure.

\begin{table}[htbp]
\centering
\caption{Image denoising. The SSIM and PSBR values of noisy image and the reconstruction images using the FTG with the fractional order $\alpha=0.5,\; 0.9,\; 1.1,\;1.5,\; 1.8$ and TG prior.}
\begin{tabular}{cccccccc}
  \toprule
          &noisy image& TG & $\alpha=0.5$ & $\alpha=0.9$ & $\alpha=1.1$ & $\alpha=1.5$ & $\alpha=1.8$ \\
  \midrule
  SSIM    &0.7526& 0.8670 & 0.7932& 0.8381 &\textbf{0.8674}& \textbf{0.8679}& \textbf{0.8684}  \\
  PSNR   &30.45&32.74& 31.15    & 32.01 &32.06 & 32.11& 32.01 \\
  \bottomrule
  \end{tabular}
  \label{table_de}
\end{table}

%
%


\section{Conclusions}\label{section5}

In this work, we have presented a FTG prior for infinite-dimensional Bayesian inverse problems.
We use the FTV term to improve the ability to capture the detail information and use the Gaussian reference measure to ensure that it results in a well defined posterior measure.
And the hierarchical Bayesian framework is also applied here, where the regularization parameter can be flexibly determined.
Moreover, we also propose an efficient diagonal map-based independence sampler for the linear inverse problems in infinite-dimensional setting.
This sampler has two stages: firstly we construct a diagonal map that can approximately pushforwards the reference measure to the posterior measure and secondly the posterior measure is explored by an independence sampler using a proposal distribution derived from the diagonal map.
Finally, we provide some numerical examples to demonstrate the performance of the FTG prior and the efficiency and robustness of the proposed independence sampler method.
We find that the FTG prior has better performance than the TG for the detail information in the unknowns, especially for recovering textures of image.
A natural extension of the present work is to use the FTG prior in other applications such as the reconstruction of rough surfaces.
And the diagonal map-based preconditioner for independence sampler can be applied to the nonlinear infinite-dimensional inverse problems.

\begin{ack}
The work described in this paper was supported by the NSF of China (11301168) and NSF of
Hunan (2020JJ4166).
\end{ack}


\end{document}